\documentclass{article}

\newtheorem{thm}{Theorem}[section] %

\newtheorem{corollary}[thm]{Corollary} %

\newtheorem{prop}[thm]{Proposition} %

\newtheorem{lemma}[thm]{Lemma} %

\newcommand{\eqn}{\begin{eqnarray}}

\newcommand{\eeqn}{\end{eqnarray}}

\begin{document}

\title{Complex Valued Analytic Torsion for Flat Bundles and for Holomorphic Bundles with (1,1)
 Connections}

\author{Sylvain E. Cappell and Edward Y. Miller  }

\date{}

\maketitle

\begin{flushleft}

\footnotesize

\centerline{ABSTRACT}

 The work of Ray and Singer which introduced analytic torsion,
a kind of determinant of the Laplacian operator in topological and
holomorphic settings, is naturally generalized in both settings. The
couplings are extended in a direct way in the topological setting to
general flat bundles and in the holomorphic setting to bundles with
(1,1) connections, which using the Newlander-Nirenberg Theorem are seen to be the bundles with both holomorphic
and anti-holomorphic structures. The resulting natural generalizations of Laplacians are not always self-adjoint
and the corresponding generalizations of analytic torsions are thus not
always real-valued. The Cheeger-Muller theorem, on equivalence in a topological setting of
analytic torsion to classical topological torsion, generalizes to this complex-valued torsion.
On the  algebraic side the methods introduced include a notion of torsion
associated to a complex equipped with both boundary and coboundry maps.

\vspace{.1in}

\S  1: Introduction.

\S 2: Coupling operators to a flat  bundle.

\S 3: Holomorphic  bundles with   type $(1,1)$ connections and a
$\bar{\partial}$--Laplacian.

\S 4: Holomorphic torsion  and Hermitian metric variation.

\S 5: An algebraic torsion for complexes with boundary and coboundary.

\S 6: Variation of algebraic torsion
with changing Hermitian inner product on $TW$.

\S 7: Definition of the flat  complex analytic torsion $\tau(M,F)$.

\S 8: Combinatorial torsion for general flat bundles over compact smooth manifolds.

\S 9: Comparison of Combinatorial and Analytic Torsions: Generalization of the Cheeger--M\"{u}ller theorem.

\S 10: Zeta functions   for d-bar setting

\end{flushleft}

\normalsize

\section{Introduction: }  \label{section1}

In  this paper     the celebrated
work of Ray and Singer \cite{RaySinger1, RaySinger2} on analytic torsion
for Riemannian and
complex Hermitian manifolds is generalized to two
natural geometrical setting, one topological and the other
holomorphic.
In their works,  coupling  geometric operators, the Riemannian
Laplacian  and $\bar{\partial}$--Laplacian, respectively,
with flat unitary bundles $E$ yielded self--adjoint operators
with  real eigenvalues whose spectral properties were encoded by these two
analytic torsions. These are hence
real numbers, in the acyclic cases considered by
Ray and Singer  and are
expressed as elements of real determinant line bundles.

Here the analytic torsion for the Riemannian case
is extended to coupling with a general flat bundle; in this general setting  the torsion has
a complex valued character.
For the complex Hermitian case, the holomorphic  torsion is extended to coupling
with an arbitrary holomorphic bundle with compatible
connection of type $(1,1)$. This includes both unitary and flat (not necessarily
unitary) bundles couplngs as  special cases.
As will be seen, by the Newlander-Nirenberg theorem \cite{NewlanderNirenberg}, and discussed in
\S 3 below,  such bundles $E$
are characterized by being endowed with   both a holomorphic structure and
also an anti-holomorphic structure.
However, in the present general setting
the associated operators are
not necessarily self--adjoint and the torsion is complex valued.
The present developments involves introducing algebraically a torion
invariant associated to a complex equipped with both boundary and coboundary maps.

In the flat unitary  Riemannian case, Ray and Singer  proved the invariance of
their analytic torsion under changes of  the Riemannian metric and offered ample evidence for their
foundational  conjecture that their analytic torsion was computable as
 a topological
invariant, the Reidemeister-Franz torsion \cite{RaySinger1} for odd dimensional manifolds.
 This   conjecture
of Ray and Singer was proved by Cheeger and
M\"{u}ller \cite{Cheeger1,Muller1}  by  different  methods of
independent interest.
 In the flat unitary holomorphic case, Ray and Singer
obtained results about the variation of their holomorphic torsion
with respect to the change of Hermitian metric of the underlying  complex manifold
\cite{RaySinger2}.

 Here the analogues of results of  Ray and Singer
 are proved in both our non--self adjoint  settings
 and the extension  of the Cheeger-M\"{u}ller theorem  is also proved.

This paper provides a general (and direct)  answer to   questions raised  by Burghelea and Haller
\cite{BurgheleaHaller4,BurgheleaHaller1} on defining torsion, algebraic and analytic,
for flat bundles.
They have pursued a different approach
\cite{BurgheleaHaller2,BurgheleaHaller3,BurgheleaHaller5} in the topological setting,  under a mild restriction
on the flat bundle, following a suggestion of M\"{u}ller
to define an analytic torsion using a different generalization
of the Laplace operator
from that utilized  here. Their method  employs a choice of
non--singular bilinear form on the flat bundle $E$
and then an anomaly cancelation term to correct for
this choice. Their  conjectured  relation to topological
invariants was proved by Su and Zhang \cite{SuZhang}
by extending previous methods of
Bismut and Zhang \cite{BismutZhang1}. The
seminal paper of Bismut and Zhang \cite{BismutZhang1}
provided the first direct
 proof of the Cheeger-M\"{u}ller
theorem and its extension to some non-unitary
contexts. Underlying this approach was the
remarkable  deformation of the
Laplacian, invented by Witten \cite{Witten1}.

The method of proof of theorem
\ref{thm.mainCheegerMuller} below generalizing the Cheeger--M\"{u}ller
theorem
proceeds via an
adaptation  of these  techniques used by
Su and Zhang \cite{SuZhang}  in  the Burghelea and Haller  setting
 to  our slightly different
operator.

More closely related to the topological part of the present work,  is
deep  work of
Braverman and Kappeller \cite{BravermanKappeller1,BravermanKappeller2,
BravermanKappeller3,BravermanKappeller4,BravermanKappeller5} which also addresses a  definition of
torsion for non-unitary flat complex bundles over
a smooth odd dimensional,  oriented manifold. They achieved this, under a mild restriction
on the flat bundle, using
a  sophisticated analysis
via the odd signature operator; the general treatment  in this paper goes directly via the
de Rham operator.
In particular, the invariant here relates to the square of
the Braverman-Kappeller invariant. In  recent papers \cite{BravermanKappeller4,BravermanKappeller5}, they related the
various  approaches  and prove under some mild restriction a version
of theorem \ref{thm.mainCheegerMuller} and so proved with this mild restriction
what was  a conjecture in an earlier circulated
version of the present work.

\vspace{.1in}
 As recalled in section \ref{section.FlatC}, there is a quite
general way of coupling any differential operator $D$ with any flat vector bundle $E$,
yielding a differential operator $D^\flat_E$. No metric on $E$ is necessary for
this procedure.
Applied to self--adjoint operators, the resulting operators, although in general not self--adjoint,
retain  most of the desirable properties of that simpler case; so natural results
about geometric operators should have corresponding interesting extensions
to this case.

This method is applied in this paper to the Riemannian Laplacian and the result
is called the flat Laplacian ${\Delta^{\flat}}_E$.
For flat
unitary bundles $E$ , one reacquires  the coupled self-adjoint
operators
considered by Ray and Singer \cite{ RaySinger1}.

As explained in section \ref{section.analyticholo1}, in the setting of complex Hermitian manifolds, there is
a further extension coupling the $\bar{\partial}$--Laplacian
to any holomorphic bundle $E$ endowed with a compatible
type $(1,1)$ connection $D$, that is a connection which is compatible with
the holomorphic structure and whose curvature has type $(1,1)$. Equivalently, by the Newlander
Nirenberg theorem \cite{NewlanderNirenberg}, the bundle $E$
is endowed with the structure of both a holomorphic structure and
also an anti-holomorphic structure. The result is called
$\Delta_{E,\bar{\partial}}$; its definition does not use any
choice of Hermitian metric on $E$.  For $E$ flat holomorphic it reduces to the flat
coupling above. For $E$ with connection compatible with a given
Hermitian inner product on $E$, $\Delta_{E,\bar{\partial}}$ recovers the self--adjoint operators
considered by Bismut, Gillet, Lebeau, Soul\'{e} \cite{BismutGillet1, BismutGillet2,
BismutLebeau1, BismutLebeau2}.
For $E$ flat unitary, one reacquires the self-adjoint d-bar operator
of Ray and Singer \cite{RaySinger2}.

In section \ref{section.analyticholo2}, the definition of the
holomorphic torsion of a holomorphic bundle $E$ with compatible type $(1,1)$ connection $D$ is given
and the local variation formula for changing Hermitian
metrics is stated  together with consequences parallel
to the theorems of Ray and Singer for holomorphic
torsion. The analogous definitions and results for the flat Riemannian
case appear in section \ref{section.analyticRiem1}.

In the complex Hermitian setting of
a holomorphic bundle with a compatible type $(1,1)$ connection $D$ defined over
a complex Hermitian manifold $W$ of complex dimension $n$, for each $p$,
$1 \le p \le n$,
the $p$ dimensional holomorphic torsion  $\tau_{holo,p}(W,E) $ will
be a non-vanishing element of a product of complex determinant line
bundles:
\begin{eqnarray}
\tau_{holo,p}(W,E) \  \in \{ \ det \ H_{\bar{\partial}_E}^{p,*}(W,E)  \otimes
\ [  det \ H_{D''}^{*,n-p}(W,E)]^{(-1)^{n+1}} \}. \nonumber  \\
 \nonumber
\end{eqnarray}
Note that these Dolbeault type cohomologies do not depend
on a choice of Hermitian inner product $g$ on $TW$.
In the case that these cohomologies vanish, the acyclic case,
this defines a complex number.

Quillen  had introduced the determinant line
bundles into the study of torsion
and interpreted the Ray-Singer torsion as a metric on this line bundle \cite{Quillen1, Quillen2}.
Here, in this  general setting, we employ  the complex determinant  line bundles  but dispense with
the metric on them. The torsion is a section of the determinant line bundle, and in particular
in the acyclic case is a complex number. From the present point of view, a Quillen metric
arises by declaring this section to have length one.

Correspondingly, the total torsion $\tau(M,F)$ for a smooth flat bundle $F$ over a
closed smooth manifold $M$ is in this approach a non--zero element
of the determinant line bundle:
$$
\tau(M,F) \in det(H^*(M, F \oplus F^*))
$$
where $F^*= Hom(F,C)$ is the dual of $F$.

As explained in section \ref{section.comb}, for any compact
smooth,  oriented manifold $M$, possibly with boundary, with flat complex bundle $F$ over $M$,
there is a combinatorial Reidemeister--Franz torsion
that takes its values in precisely this
determinant line for any flat bundle $F$.
It utilizes directly a subdivision and
dual subdivision of the manifold $M$.

\vspace{.1in}
A key  technical difficulty of this paper is that
it is \textbf{not} sufficient to just use the generalized eigenvalues
of $\Delta^{\flat}_{F}$, respectively $\Delta_{E,\bar{\partial}}$,
  to define the analytic
or holomorphic torsion.
The zero--modes have additional structure:  in each case there
are two differentials, one going up and one going down, that commute with the operator
$\Delta^{\flat}_{F}$, $\Delta_{E,\bar{\partial}}$ respectively.
One must encode this information in the definition of analytic
or holomorphic torsion
to get control of the changes under metric variations.
The added term is an algebraic torsion associated to the
zero modes with these two differentials \ref{section.algebraictorsions};  it is defined
and related to the standard Reidemeister--Franz
torsion of a complex with bases for the cochains.

In section \ref{section.algebraictorsions},
a description of the algebraic torsion of a finite complex
with two differentials, one going up and one going down,
is given. It is related to the more classical approaches
and expressed in terms of eigenvalues in the acyclic
cases.

More explicitly, if $F$ is a flat bundle over a smooth manifold $M$, then the flat extensions $d_F,d^*_F$ of the
exterior derivative $d$ and its adjoint $d^*$ via Riemannian metric
are operators on the smooth forms with values in $F$, $A^*(M,F)$. These two each have
square zero and commute with the flat Laplacian $\Delta^{\flat}_F$.
Hence, the generalized zero modes of $\Delta^{\flat}_F$ is endowed with
two differentials. To such a  structure we associate  an algebraic torsion
and it is incorporated into the definition of the analytic torsion.
Similarly  for the complex holomorphic cases, we use the added information
contained in the zero modes to define the holomorphic torsion.

Since signs are a key  issue, to avoid a proliferation
of different conventions and facilitate comparisons, we have decided to largely
adopt  conventions of Braverman and Kappleller \cite{BravermanKappeller3}
and  Bismut and Zhang \cite{BismutZhang1}; these are in agreement but differ slightly
from those of Turaev \cite{Turaev1} and
Milnor \cite{Milnor1}.
For example, in accordance with  Braverman and Kappeller,
 we will take our basic complex to be
a co-chain complex  with differential of degree $+1$, $(C^*,d)$,
$d:C^q \rightarrow C^{q+1}, \  d^2 =0$.
In section 4,  the basic algebra concerns the situation where
in addition there is  a differential
of degree $-1$, $d^*$, on this same $C^*$.

\vspace{.1in}

In section \ref{section5}  the variation of the holomorphic torsion is computed
via an analysis of the variation of the algebraic torsion associated to the
generalized zero-modes.
The parallel development for the flat Laplacian case appears in section \ref{section.analyticRiem1}

The combinatorial torsion for a general flat bundle $F$ over a closed Riemannian manifold $M$
is described in section \ref{section.comb}.

 The proof of the generalized  Cheeger-M\"{u}ller
theorem, the equality of analytic and combinatorial torsions,  is carred out in section \ref{section.mainCheegerMuller}.
The method of proof is modeled on the work of Su and Zhang \cite{SuZhang} with the slight
modifications needed for our slightly different operator. As mentioned above, an analogue
 has been previously proved by Braverman and Kappeller \cite{BravermanKappeller4,BravermanKappeller5}
in a slightly restrictive setting.

The idea of regarding the original Ray--Singer invariant, a real
type invariant, \cite{RaySinger2},
as giving a real norm on the determinant bundles
originated in Quillen's work \cite{Quillen1},\cite{Quillen2}
and has been extensively utilized
in the seminal work of Bismut and his collaborators, \cite{Bismut1, BismutGillet1, BismutGillet2,
BismutLebeau1, BismutLebeau2, BismutZhang1}. It is natural to conjecture
analogues to the results in those papers and others  in the present settings, i.e., dropping
the Hermitian metrics and systematically obtaining complex valued formulae.
Indeed, the naturality and good fit of the analytic details in the present use
of $(1,1)$ connections, suggests that it may provide natural settings for some extensions of notions
of global analysis, symplectic geometry, and gauge theory.

The authors are grateful for the encouragement of I. Singer, J. Cheeger, W. M\"{u}ller,
and J.-M. Bismut, and their comments on an early version circulated several years ago.

\section{Coupling operators to a flat  bundle:}
\label{section.FlatC}

There is a simple general principle: If $D$ is a linear differential
operator mapping smooth sections of a complex bundle $E_1$ to smooth sections of
a complex bundle $E_2$ over a smooth manifold $M$, i.e.,
\eqn
D : \ \Gamma(E_1) \rightarrow \Gamma(E_2) ,
\eeqn
then for any complex flat bundle $E$ there is a canonically
determined associated differential operator $D_E^{\flat}$ from smooth sections
of the tensor product $E_1 \otimes E$ to smooth sections of  $E_2 \otimes E$
\eqn
D^{\flat}_E : \Gamma(E_1 \otimes E) \rightarrow \Gamma(E_2 \otimes E)
\eeqn
This is specified as follows:

\vspace{.3in}
Since a differential operator is defined  locally, we need but
specify $D_E^{\flat}$ over an open set, say $U \subset M$ over which
the flat bundle $E$ is trivial. Let the restriction $E|U$ have
a basis of flat smooth sections $s_1,s_2,\cdots, s_k$, with $k$ the dimension
of the bundle $E$. Then define
\eqn
(D_E^{\flat}|U) :\  \Gamma(E_1 \otimes E|U) \rightarrow \Gamma(E_2 \otimes E|U)
\eeqn
via $\Sigma_{i=1}^k \ f_i \otimes s_i \mapsto \Sigma_{i=1}^k \ D(f_i) \otimes s_i$
for any smooth sections $f_i,i=1,\cdots,k$ of $E_1|U$.
Since the transition mappings comparing flat sections are
constant and $D$ is linear, these local definitions are compatible
and define the desired operator $D_E^{\flat}$.

Proceeding thusly for any elliptic operator, say $D$,  the resultant operator
has symbol $\sigma(D_E^{\flat})$ equal to the symbol of $D$ tensor the identity
on $E$, so
the operator $D_E^{\flat}$ is elliptic.
The primary  difference is that for general flat $E$, one
is led from self-adjoint elliptic operators to
non-self adjoint ones.

Applied to the Laplace operator, $ \Delta$, of a Riemannian manifold $M$, this
construction produces the flat Laplacian, $\Delta_E^{\flat} : A^*(M,E) \rightarrow A^*(M,E)$,
acting on smooth forms with values in the flat bundle $E$. In the case that $E$ is unitary,
this flat Laplacian will be the standard self-adjoint extension of the Laplacian $\Delta$
where the unitary structure is used to define the inner-product. However, without any
choice of Hermitian metric, using just the flat structure on $E$,
the flat Laplacian $\Delta_E^\flat$ is well-defined.

\section{
Holomorphic  bundles with   type $(1,1)$ connections and a
$\bar{\partial}$--Laplacian } \label{section.analyticholo1}

There is a natural parallel construction for the $\bar{\partial}$--Laplacian. It is described
in this section.

Let $W$ be a complex manifold
with Hermitian inner product, say $g = <..>$   and $E \rightarrow W$ be a complex
holomorphic, bundle over $W$ endowed with a linear connection $D$.

\vspace{.2in}
Let $\bar{\partial}: A^{p,q}(W,C) \rightarrow A^{p,q+1}(W,C), \partial
: A^{p,q}(W,C) \rightarrow A^{p+1,q}(W,C)$ be the standard
operators obtained by decomposing by type the exterior derivative
$$
d = \bar{\partial} + \partial
$$
acting on complex valued smooth forms of type $(p,q)$. By $d^2=0$,
one has $(\bar{\partial})^2 =0 , \ (\partial)^2 = 0$.
Let $<.,.>$ be the induced Hermitian metric on these $(p,q)$ forms and
$\bar{\partial}^*$ denote the adjoint of $\bar{\partial}$ under this Hermitian
metric on $TW$
$$
<\bar{\partial} s_1, \ s_2> = < s_1 , \ \bar{\partial}^*  s_2> \ for \ s_1 \in A^{p,q-1}(W,C) \
and \ s_2 \in A^{p,q}(W,C)
$$
here $\bar{\partial}^* : A^{p,q}(W,C) \rightarrow A^{p,q-1}(W,C)$.

\vspace{.2in}
By $E$ holomorphic,
 the $\bar{\partial}$ operator
on forms on $W$ with complex values
has a unique natural extension to the smooth forms with values in $E$,
$$
\bar{\partial}_E : A^{p,q}(W,E) \rightarrow A^{p,q+1}(W,E).
$$
This  first order differential operator is uniquely characterized by: \newline
1) $\bar{\partial}_E ( f \ s) = (\bar{\partial}_E f) \wedge  \ s + f (\bar{\partial}_E \ s)$
for any smooth function $f$; and,  \newline 2) over any open set $U$ for which $t$ is
holomorphic section of $E$ and any smooth $(p,q)$ form $a$ with complex values
$\bar{\partial}_E ( a \otimes  t)|_U = (\bar{\partial} a)\otimes \ s|_U$.

\vspace{.2in}

A linear connection $D$ on $E$ is by definition a linear mapping of smooth sections:
$$
D :  \Gamma(E) \rightarrow \Gamma(T^*W \otimes_R E) = \Gamma((T^*W\otimes_R C) \otimes_C E)
$$
satisfying the Leibnitz formula $D(f \ s) = df \otimes  s  + f \ D(s)$.

The complex structure on $TW$ provides the natural splitting
$$
T^*W \otimes C \cong {T^*}'W \oplus {T^*}''W
$$
with ${T^*}'W$, respectively ${T^*}''W$, the $\mp \ i$ eigenspaces of
the operator $(J \otimes Id)$ where $J :TW \rightarrow TW, \ J^2 = -Id$
specifies the complex structure. Hence,  the smooth one forms with values in $E$,
$ \Gamma((T^*W\otimes_R C) \otimes_C E)$ splits as a direct sum:
$$
\begin{array}{l}
 \Gamma((T^*W \otimes_R C) \otimes_C E) \\
\cong
  A^{0,1}(W,E) \oplus A^{1,0}(W,E) \cong  \Gamma({T^*}'W \otimes_C E) \oplus  \Gamma({T^*}''W \otimes_C E)
\end{array}
$$

Under this decomposition, the connection $D$ decomposes as a sum:
$D = D' \oplus D''$ with $D' : A^{0,0}(W,E)= \Gamma(E) \rightarrow A^{0,1}(W,E)$ and
$D'' : A^{0,0}(W,E)= \Gamma(E)  \rightarrow A^{1,0}(W,E)$.

The connection $D$ is said to be \textbf{compatible} with the holomorphic structure
on $E$ if $D' = \bar{\partial}_E$. In particular, in this case
any holomorphic section, $t$ of $E$ over an open set $U$ of $W$,
satisfies $D'(t) = \bar{\partial}_E t = 0$. The natural extension
of $D''$ to all $(p,q)$ forms is then
$D'' := \bar{\partial }_E : A^{p,q}(W,E) \rightarrow A^{p,q+1}(W,E)$
in this case.

Now as above, the operator $D''$ has a unique extension:
$$
D'' : A^{p,q}(W,E) \rightarrow A^{p+1,q}(W,E)
$$
$D''$  is  uniquely characterized by: \newline
1) $D'' ( f \ s) = (\partial f) \wedge  \ s + f (D'' \ s)$
for any smooth function $f$; and, \newline 2)  for any smooth section $t$
of $E$ and any smooth $(p,q)$ form $a$ with complex values
$D'' ( a \otimes t)|_U = (\partial  \otimes  a)\ t + (-1)^{p+q} \ a \ \wedge D''(t)$.

In the standard manner one defines the curvature of the linear connection
$D$ via
$$
curvature = D^2 = (D')^2 + (D' D'' + D'' D') + (D'')^2
$$
This curvature is a two form with values in $End(E)$, the self maps of $E$,
and its type decomposition $A^{0,2}(W,End(E)) \oplus A^{1,1}(W,End(E)) \oplus A^{2,0}(W,End(E))$
is $(D')^2,(D' D'' + D'' D') , (D'')^2$. We say the connection is
of \textbf{type $(1,1)$} if its  curvature is of type $(1,1)$, that is,
$$
(D')^2=0, \ and \  (D'')^2 = 0
$$

The first condition follows from $D' = \bar{\partial}_E$ when $D$
is compatible with the holomorphic structure on $E$. By the Newlander--Nirenberg
theorem \cite{NewlanderNirenberg, KobayashiNomitzu} the condition $(D'')^2=0$ is equivalent
to having a holomorphic structure on the conjugate $\bar{E}$ of $E$,
or equivalently an anti--holomorphic structure on $E$.

\vspace{.2in}
Conversely,  if $E$ is endowed with a holomorphic structure, there is a unique
linear operator $D'  : \Gamma(E) \rightarrow \Gamma( T^{0,1}W \otimes E)$ with
$D' ( f \ s) = \overline{\partial} f \otimes  \ s + f (D' \ s)$
for any smooth function $f$ and $D'(s) = 0$ for any holomorphic
section $s$. The condition that any point has an open neighborhood $O$
for which $E$ has a basis, say $\{s_i \}$ of holomorphic sections,
implies that $(D')^2 = 0$, since for any smooth section $s$, one has
$s| 0 = \Sigma_i \ h_i \ s_i$ and consequently, $(D')^2(s)|0 =
D'( \overline{\partial} h_i \otimes s_i ) = \overline{\partial} \ \overline{\partial} h_i  \otimes s_i
=0$. If in addition, $E$ is endowed with an anti-holomorphic structure,
then there is a unique linear operator $D'' :  \Gamma(E) \rightarrow \Gamma( T^{1,0}W \otimes E)
$ with $D'' ( f \ s) = (\partial f) \otimes   \ s + f (D'' \ s)$
for any smooth function $f$ and as before $(D'')^2 = 0$. Together
the sum $D := D' + D''$ enjoys the property
$D(fs) =  ( (\overline{\partial} f ) \otimes   \ s + f (D' \ s)) +((\partial f) \otimes  \ s + f (D'' \ s))
= df \otimes s + f ( D(s))$. That is, $D$ is a standard linear connection. It is of type $(1,1)$
since $(D'')^2 = (D')^2 = 0$.

These observations prove the proposition:

\begin{prop}
A complex bundle $E$ endowed with a holomorphic and anti-holomorphic
structure defines a unique connection $D$  of type $(1,1)$  compatible
with these two structures. Here compatible means that the type decomposition
$D = D' + D''$ yields operators with $D' s=0$ for holomorphic sections
and $D''s = 0$ for anti-holomorphic sections.

Conversely, a complex bundle $E$ endowed with a type $(1,1)$
connection has a uniquely determined holomorphic and
anti-holomorphic structure.

\end{prop}
The last statement uses  the  Newlander--Nirenberg
theorem \cite{NewlanderNirenberg, KobayashiNomitzu}.

\vspace{.2in}
Let $D$ be assumed to be a linear connection which is
compatible with the given holomorphic structure and
having curvature of type $(1,1)$.

To specify the desired generalization of the standard
d-bar Laplacian on smooth forms
$$
\Delta_{\bar{\partial}} = \bar{\partial} \partial+\partial \bar{\partial}
: A^{p,q}(W,C) \rightarrow A^{p,q}(W,C)
$$
to the smooth $(p,q)$ forms with values in $E$, $A^{p,q}(W,E)$,
it is helpful give  more  explicit details.

\vspace{.2in}

 Somewhat  confusingly, the star operator
$\star $ acting on forms  is a complex conjugate  linear mapping
\eqn
\star : A^{p,q}(W,C) \rightarrow A^{n-p,n-q}(W,C)
\eeqn
induced by a conjugate linear bundle isomorphism. It is characterized by
 the property that for $a,b \in   A^{p,q}(W,C)$,
\eqn
<a,b> = \int_W \ a \wedge \star b
\eeqn
where $<a,b> = \int_W  <a(x),b(x)> \  dvol_W$ with
$<a(x),b(x)>$ the pairing over $x \in W$ induced by
the Hermitian inner product on $T_xW$.

Now $\star$ composed with the natural conjugation mapping
\eqn
conj : A^{p,q}(W,C) \rightarrow A^{q,p}(W,C)
\eeqn
is a complex linear mapping, induced by a bundle
isomorphism. Here $conj$  is induced by the
bundle automorphism  $T^*W \otimes_R  C \rightarrow T^*W\otimes_R  C,\
v \otimes \lambda \mapsto v \otimes \bar{\lambda}$ of the complexified
cotangent bundle.  Denote this induced,  composite, complex linear mapping   by $\hat{\star} $. That is,
\eqn
\hat{\star} = conj \ \star : A^{p,q}(W,C) \rightarrow A^{n-q,n-p}(W,C)
\eeqn

This complex linear isomorphism is induced by a complex linear mapping
of bundles, so may be coupled with any bundle $E$ via $\hat{\star} \otimes Id_E$
\eqn
\hat{\star} \otimes Id_E : A^{p,q}(W,E) \rightarrow A^{n-q,n-p}(W,E)
\eeqn

 It should be noted that $\hat{\star}$ being  complex linear
may be coupled to a  complex linear bundle mapping, such as the identity mapping.
The bundle mapping $\star$ can be coupled only to a conjugate complex
linear mapping.

 Recall  again that  the exterior derivative $d $
has the type decomposition
$
d = \bar{\partial} + \partial
$
decomposition, where $\bar{\partial} :A^{p,q}(W,C) \rightarrow A^{p,q+1}(W,C)$
and $\partial : A^{p,q}(W,C) \rightarrow A^{p+1,q}(W,C)$.These  correspond under conjugation
$$
\partial = conj \ \bar{\partial}\  conj : A^{p,q}(W,C) \rightarrow A^{p+1,q}(W,C).
$$
Recall also that  the adjoint of $\bar{\partial}$ under the chosen
Hermitian inner product on $TW$ has the  well known explicit description
$$
{\bar{\partial}}^* = - \ \star \bar{\partial} \star.
$$

In particular,
$$
{\bar{\partial}}^*= - \hat{\star}\  \  conj \ \bar{\partial} \ conj \ \hat{\star}
 = - \hat{\star} \ \partial \  \hat{\star}.
$$

Motivated by this formula, since $D''$ extends $\partial $ and $\hat{\star} \otimes Id_E$
extends $\hat{\star}$ to the sections of $E$, one defines
$$
\bar{\partial}^{*}_{E,D''} = - ( \hat{\star} \otimes Id_E) \cdot
      D'' \cdot ( \hat{\star}  \otimes Id_E)
$$
and introduces the desired generalization of the
d-bar Laplacian for $E$ a holomorphic bundle
with compatible type $(1,1)$ connection by:
$$
\Delta_{E,\bar{\partial}} = \bar{\partial}_E \bar{\partial}^{*}_{E,D''}
+ \bar{\partial}^{*}_{E,D''} \bar{\partial}_E
$$
Note that $(\bar{\partial}^{*}_{E,D''})^2 = 0$ since $(D'')^2 = 0$ for
the given type $(1,1)$ connection $D$ and $(\bar{\partial}_E)^2 = 0$.
Also, both the differentials $\bar{\partial}^{*}_{E,D''}$ and $\bar{\partial}_E$ commute
with $\Delta_{E,\bar{\partial}}$.

\vspace{.2in}
In the case that $E$  has a flat connection, called $D$, utilizing the method of section
\ref{section.FlatC}, $\bar{\partial}^{*}_{E,D''}$ above
is just the coupled version of ${\bar{\partial}}^*
= - \hat{\star} \ \partial \  \hat{\star}$ and $\Delta_{E,\bar{\partial}}$
is the coupled  version of the d-bar Laplacian $\Delta_{\bar{\partial}}$.

\vspace{.3in}
It is important to note that if a smooth Hermitian inner product
is chosen on $E$,
say $<.,.>_E$, then the adjoint of $\bar{\partial}_E$
is \textbf{not} in general the above mapping.  For example, its
definition does not utilize any such choice.
Rather, as specified by M\"{u}ller \cite{Muller2},
the adjoint is defined in terms of the induced conjugate linear
bundle isomorphism $\mu : E \rightarrow E^*$, of $E$ and its dual $E^*$ as:
\begin{eqnarray}
adjoint(\bar{\partial}_E) =
- (\star^{-1} \otimes \mu)^{-1}  \ \bar{\partial}_{E^*} \
   \ (\star \otimes \mu) \label{eqn.operadj1}
\end{eqnarray}
With this definition $<\bar{\partial}_E s, t>_E =
 <s, adjoint(\bar{\partial}_E)\  t>_E$  for the induced
inner product on forms with values in $E$.
In most of the literature, this adjoint is (rather sloppily) written
as $-\star \bar{\partial} \ \star$.
Note that the conjugate complex linear operators, $\star$ and
$\mu$ appear here, not the complex linear one $\hat{\star}$.
[Since both $\star$ and $\mu$ are both conjugate linear bundle isomorphisms, this makes sense.]

\vspace{.2in}
In the case that the holomorphic bundle $E$ is endowed with a Hermitian inner product, say $<.,.>_E$, which
is compatible with the connection $D$, in the sense that
$$
d \ <s_1,s_2>_E = < D \ s_1, s_2>_E + <s_1, D \ s_2>_E,
$$
that is, the connection $D$ is unitary,
then the conjugate linear mapping $\mu$ identifies the dual bundle $E^*$ equipped with its dual connection $D^*$
with $E$ equipped with its connection $D$. Consequently in  this unitary situation,
the adjoint $adjoint(\bar{\partial}_E)$ equals precisely $\bar{\partial}^*_{E,D''}$
and so the standard self-adjoint d-bar Laplacian equals the above one,
$$
\Delta_{E,\mu} = (\ \bar{\partial}_E +  adjoint(\bar{\partial}_E) )^2 \ .
= \Delta_{E,\bar{\partial}}
$$
Hence, the above operator restricts to the standard self--adjoint
one considered by Ray and Singer \cite{RaySinger2}, Bismut, Gillet,  Lebeau, Soul\'{e}  in \cite{Bismut1, BismutGillet1, BismutGillet2,
BismutLebeau1, BismutLebeau2}
for $E$ flat unitary.

\vspace{.2in}

Nonetheless, in all type $(1,1)$ compatible holomorphic cases, it is apparent from these formulas that the first order
differential operators
$\bar{\partial}^{*}_{E,D''}$ and the adjoint
$adjoint(\bar{\partial}_E)$ have identical symbols and
thus  differ by a smooth bundle automorphism, say $\alpha$,
$$
\bar{\partial}^{*}_{E,D''} = adjoint(\bar{\partial}_E) + \alpha \ .
$$
Thus in all cases,
\begin{eqnarray}
\Delta_{E, \bar{\partial}} = (\ \bar{\partial}_E +  adjoint(\bar{\partial}_E) )^2
 + \alpha \bar{\partial}_E +  \bar{\partial}_E \alpha
\eeqn
 is an elliptic second order partial differential equation with
scalar symbol.

\vspace{.3in}
Since $D$ is a type $(1,1)$ compatible connection,  it follows that
\eqn
\bar{\partial}_E : A^{p,*}(W,E) \rightarrow A^{p,*}(W,E) \ has \ (\bar{\partial}_E)^2 = 0 \\
D'': A^{*,n-p}(W,E) \rightarrow A^{*,n-p}(W,E)  \ has \ (D'')^2 = 0,
\eeqn
so two Dolbeault type homologies
$ \ H_{\bar{\partial}_E}^{p,*}(W,E)$ and
$\ H_{D''}^{*,n-p}(W,E)$ are well defined.
These  cohomologies  appear in the  main theorem
for the holomorphic torsions.
The proposed invariant is to lie in the product of determinants.
\begin{eqnarray}
\tau_{holo,p}(W,E) \  \in \{ \ det \ H_{\bar{\partial}_E}^{p,*}(W,E)  \otimes
\ [  det \ H_{D''}^{*,n-p}(W,E)]^{(-1)^{n+1}} \}. \nonumber  \\
\end{eqnarray}
It will extend the Ray-Singer holomorphic torsion \cite{RaySinger2}  from their case
of $E$ unitarily flat, and so holomorphic, to this more general and flexible setting of
coupling with any holomorphic bundle with compatible type $(1,1)$
connection. [For $E$ flat holomorphic,
one can replace $D''$ by $\partial^\flat_E$, the flat extension
of $\partial$.]

Its variation with the choice of Hermitian metric $g$ on $TW$
is specified in theorem \ref{thm.MAIN} . The cohomologies
are independent of that choice. This theorem is a generalization of
the Ray-Singer theorem of \cite{RaySinger2}.

\vspace{.1in}

\textbf{Foundational Example:}

\vspace{.1in}

An important universal example of a holomorphic $(1,1)$ connection is described as follows:
Let $V$ be a complex vector space of dimension $n$, and $Gr_k(V)$ denote the Grassmannian
of $k$-dimensional subspaces of $V$. That is, $Gr_k(V) = \{ E_1 \subset V\ | \ dim_C \ E_1 = k \}$.
Similarly, let $Gr_{n-k}(V) =\{ E_2 \subset V\ | \ dim_C \ E_2 =(n- k) \}$,  the $n-k$ dimensional
subspaces.
Let
$$
{ \cal{O}} = \{ (E_1,E_2) \ | \ (E_1, E_2 ) \in \ Gr_{k}(V) \times Gr_{n-k}(V) \ with \ E_1 \cap E_2 = \{0\} \}
$$
That is, ${ \cal{O}}$ is the open subset of those pairs that are transverse, or equivalently
$E_1  +  E_2 \cong V$.

Over ${ \cal{O}}$ one has  the natural subbundle
$$
{\cal{S}}_1 = \{ p \in E_1 \ | \ (E_1,E_2 ) \in U \}
$$
and a natural inclusion of ${\cal{S}}_1 \subset V \times Gr_k(V)$ into the trivial bundle
and a natural projection defined by $V \mapsto V /E_2 \cong E_1$ of
this trivial bundle $V \times Gr_k(V)$ back onto ${\cal{S}}_1$. Let these bundle mappings be denoted by
$$
i : {\cal{S}}_1 \subset V \times Gr_k(V)  \ and  \ p:   V \times Gr_k(V)\rightarrow {\cal{S}}_1
$$

Hence, one  may define a canonical connection $D$ on smooth sections of ${\cal{S}}_1$ by the following procedure.
For
$f $ a smooth section of ${\cal{S}}_1$, that is $ f \in \Gamma({\cal{S}}_1)$ and
$X$ a vector field on ${ \cal{O}}$, use the standard derivative $\nabla_X$ on the product
bundle $V \times Gr_k(V)$ applied to $i(f)$ and then project via $p$ to ${\cal{S}}_1$
again. That is, set
$$
D_X(f) = p \ \nabla_X ( i(f))
$$
By construction for $g$ a smooth function on $U$, $D_X(g \ f) = p \ \nabla_X( i (g\ f))
=p \ \nabla_X(g  \ i (f)) = p \ ( g \  \nabla_X( i (f)) + (X\ g) \ ( i (f)) )
= g \ p \ \nabla_X( i (f)) + (X\ g) \ f = g \ D_X (f) + (X\ g) \ f$, so
$D$ is a connection on the bundle ${\cal{S}}_1$ as desired.

It is claimed that the connection $D$  is naturally a holomorphic connection of type $(1,1)$
for a suitable complex structure on ${ \cal{O}}$. This is the complex structure on ${ \cal{O}}$
induced from that on the product of Grassmanians $G_k(V) \times Gr_{n-k}(V)$
where the complex structure on the first factor $Gr_k(V)$ comes from the standard
complex structure $J = i$ on $V$ and the complex structure on the second factor
$Gr_{n-k}(V)$ comes from the conjugate complex structure $-J$ acting on $V$.

\vspace{.1in}
To see this  concretely fix a pair $(E_1,E_2)$ in $U$ and
take linearly independent vectors $e_1,e_2,\cdots, e_n$ where $E_1$ is the span of
$\{e_1,\cdots, e_k\}$ and $E_2$ is the span of $\{e_{k+1},\cdots, e_n\}$.
Now for complex coordinates $z_{p,q}, \ 1 \le p \le k, 1 \le q \le (n-k)\}$ specify the
$k$ vectors
$$
s_p = e_p \ + \ \Sigma_{q=1}^{n-k} \ z_{p,q} \ e_{k+q}
$$
and let $E(\{z_{p,q}\})$ be the span of these $k$ vectors. This provides an open
set, say $O_1$, in $Gr_k(V)$ and a complex coordinate system for the complex structure induced from $J$
acting on $V$. Here
$E_1$ is the $k$-plane with vanishing coordinates. Let $Z$ denote
the complex $k \times (n-k)$ matrix with entries $z_{p,q}$.
Let $\vec{s} = \left( \begin{array}{l} s_1 \\ \vdots \\ s_k  \end{array} \right)$
the column vector with entries the $s_p$'s.

Similarly,  for complex coordinates $u_{r,s}, \ 1 \le r \le (n-k), 1 \le s \le k\}$ specify the
$n-k$ vectors
$$
t_r =  \Sigma_{s=1}^k \ \overline{u_{r,s}} \ e_s \ + \ e_{k+r}
$$
and let $E(\{\overline{u_{r,s}} \})$ be the span of these $n-k$ vectors. Note the conjugates
in the definition of $t_r$. This provides an open
set, say $O_2$,  in $Gr_{n-k}(V)$ and a complex coordinate system for the complex structure induced from $-J$.
Thus $E_2$ is the $n-k$-plane with vanishing coordinates. Let $\overline{U}$ denote
the complex $(n-k) \times k $ matrix with entries $\overline{u_{r,s}}$.
Let $\vec{t} = \left( \begin{array}{l} t_1 \\ \vdots \\ t_{n-k}  \end{array} \right)$
the column vector with entries the $t_r$'s.

Together the complex coordinates $z_{p,q}, u_{r,s}$ specify the desired complex structure on the
open set $O_1 \times O_2$ in the product $G_k(V) \times Gr_{n-k}(V)$.
The $s_p$'s provide a basis of sections of the subbundle ${\cal{S}}_1$ while the $t_r$'s provide
as basis of sections of the subbundle ${\cal{S}}_2$.

The relation between the vectors $s_p,t_r$ and the standard vectors basis vectors $\{e_1,\cdots, e_n\}$
is
$$
\left( \begin{array}{l} \vec{s} \\ \vec{t} \end{array} \right)
= \left( \begin{array}{ll} I_k &  Z \\ \overline{U} & I_{n-k} \end{array} \right)   \
          \left( \begin{array}{l} e_1 \\ \vdots \\ e_n \end{array} \right)
$$
where $I_l$ is the identity matrix with $l$ rows and columns. Let $W$ denote the above $n \times n$ matrix.
The condition that the subspaces $E(\{z_{p,q}\}),E(\{\overline{u_{r,s}} \})$ be transverse is
precisely that $det \ W  \neq 0$. That is, on ${ \cal{O}} \cap ( O_1 \times O_2) $
these provide coordinates. Let the inverse of the matrix $W$, which is defined on  ${ \cal{O}} \cap ( O_1 \times O_2)$,
be written also in block form:
$$
W
= \left( \begin{array}{ll} I_k &  Z \\ \overline{U} & I_{n-k} \end{array} \right)
          \ and \
W^{-1} = \left( \begin{array}{cc} \ \alpha & \beta \\ \gamma & \delta \end{array} \right)
$$

In these terms one computes,
$$
\begin{array}{l}
D(\vec{s}) = p \ ( \nabla \ \ i \left(  \vec{s}  \right))  = p\ ( \left( \begin{array}{ll} 0 ,   &  dZ \end{array} \right)
 \left( \begin{array}{l} e_1 \\ \vdots \\ e_n \end{array} \right) )  = \\
 p \ (
\left( \begin{array}{ll} 0 ,   &  dZ \end{array} \right) \
\left( \begin{array}{cc} \ \alpha & \beta \\ \gamma & \delta \end{array} \right) \
\left( \begin{array}{l} \vec{s} \\ \vec{t} \end{array} \right) ) \\

= p \ ( dZ \cdot \gamma \cdot \vec{s} + dZ \cdot \delta \cdot \vec{t}) =  dZ \cdot \gamma \cdot \vec{s}
\end{array}
$$
That is, the connection form for this basis $\vec{s}$ is $\Theta = dZ \ \cdot \   \gamma$, a 1-form of type $(1,0)$.

Correspondingly, the curvature $K$ is $K = d\Theta - \Theta \wedge \Theta$ via
$D^2(\vec{s}) = D( \Theta \otimes \vec{s}) = d\ \Theta  \otimes \vec{s} - \Theta \wedge \ D(\vec{s})=
(d\Theta - \Theta \wedge \Theta) \otimes \vec{s}$.

In this instance $K = d( dZ \cdot \gamma) - dZ \cdot \gamma \wedge dZ \cdot \gamma
= - dZ \ \wedge d\gamma - dZ \cdot \gamma \wedge dZ \cdot \gamma
= - dZ \wedge  (d\gamma + \gamma  \cdot dZ \cdot \gamma)$.

However, in view of the identity
$$
\begin{array}{l}
d ( W^{-1}) = \left( \begin{array}{cc} \ d\alpha & d\beta \\ d \gamma & d \delta \end{array} \right)
= - W^{-1} \cdot dW \cdot W^{-1} \\
= - \left( \begin{array}{cc} \ \alpha & \beta \\  \gamma &  \delta \end{array} \right) \cdot
\left( \begin{array}{ll} 0 &  d\ Z \\ \overline{d\ U} & 0 \end{array} \right) \cdot
\left( \begin{array}{cc} \ \alpha & \beta \\  \gamma &  \delta \end{array} \right)  \\
= - \left( \begin{array}{cc} \alpha & \beta \\  \gamma & \delta \end{array} \right) \cdot
\left( \begin{array}{cc} \ dZ \cdot \gamma & dZ \cdot \delta \\ \overline{d\ U} \ \cdot \alpha & \overline{d \ U} \cdot
\beta \end{array} \right) \\
= - \left( \begin{array}{ll} \star & \star \\ \gamma \cdot dZ \cdot \gamma
+ \delta \ \overline{d\ U} \cdot \alpha & \star \end{array} \right) \ ,
\end{array}
$$
one obtains the equality $d\gamma = - \gamma \cdot dZ \cdot \gamma - \delta  \cdot \overline{d\ U} \cdot \alpha$, so
the curvature $K$ is
$$
K = + \  dZ \wedge (\delta \ \overline{d\ U} \cdot \alpha)
$$
a curvature of type $(1,1)$, as expected.

Note that the sections $\vec{s}$ provide a holomorphic basis of sections since $D(\vec{s})= dZ \cdot \gamma \cdot \vec{s}
$ is of type $(1,0)$. This defines the required holomorphic structure.

Similarly,  $(I- Z \cdot \gamma) \vec{s}$ provides a basis of anti-holomorphic sections, since
$D( (I- Z \cdot \gamma) \vec{s})$ is of type $(0,1)$ as is seen from the computation:
$$
\begin{array}{l}
D( (I- Z \cdot \gamma) \vec{s}) = - (dZ \cdot \gamma + Z \cdot d\gamma)  \cdot \vec{s}+(I- Z \cdot \gamma)  \cdot D(\vec{s})\\
= - (dZ \cdot \gamma + Z \cdot d\gamma) \cdot \vec{s} +  (I- Z \cdot \gamma) \cdot ( dZ \cdot \gamma \cdot \vec{s}) =
- Z \cdot d \gamma \cdot \vec{s} - Z \cdot \gamma \cdot ( dZ \cdot \gamma \cdot \vec{s}) \\
= - Z \cdot ( d\gamma + \gamma \cdot ( dZ \cdot \gamma) \cdot \vec{s} =
+ \  Z \cdot (\delta  \cdot \overline{d\ U} \cdot \alpha ) \cdot \vec{s}
\end{array}
$$
This defines the corresponding anti-holomorphic structure.

\vspace{.1in}
More conventionally, one chooses a Hermitian metric, say $<.,.>$ on $V$, then the subbundle, say
${\cal{E}}_1$, over the Grassmanian $Gr_k(V)$ acquires an induced Hermitian metric, so there is
a unique connection $D$ which is compatible with the metric $<.,.>$, in that
$$
d<s_1, s_2>  = < Ds_1,s_2> + < s_1, Ds_2>
$$
for any smooth sections $s_1,s_2$ of this subbundle.
This compatible connection  then has curvature of type $(1,1)$.

Using the Hermitian metric, for each $k$-dimensional subspace $E_1$ let
$(E_1)^\perp = \{ w \ | \  <v,w> = 0 \}$ be the orthogonal complement, a $(n-k)$ dimensional
subspace. Hence, the orthogonal projection of $V$ onto $E_1$ is precisely the projection
of $V$ via the mapping $p$, $V \rightarrow V/ (E_1)^\perp \cong E_1$.
 It is evident that for sections $s_1,s_2$ of ${\cal{E}}_1$
regarded as sections of the trivial bundle, i.e., as $i(s_1),i(s_2)$ for the inclusion
${\cal{E}}_1 \subset V \times Gr_k(V)$, one has the formula:
$$
\begin{array}{l}
X< i(s_1),i(s_2)> = < \nabla_X \  i(s_1), i(s_2)> + < i(s_1),\nabla_X \ i(s_2)> \\
  = <p \  \nabla_X \  i(s_1), s_2> + < s_1 , p \ \nabla_X \ i(s_2)>
  \end{array}
$$
That is, the connection given by $f \mapsto p \  \nabla \  i(f)$ satisfies the defining
property above. Consequently, one has
$$
D(f) = p \  \nabla \  i(f).
$$
That is, the standard compatible connection associated to the Hermitian metric $<.,.>$
is obtained in this natural fashion.

An alternative method to see this is to directly compare the respective connection 1-forms and
see their equality.
Given a Hermitian
metric, one can define a mapping:
$$
\begin{array}{l}
F : Gr_k(V) \subset Gr_k(V) \times Gr_{n-k}(V) \\
via  \\
E_1 \mapsto ( E_1, (E_1)^{\perp} ) . \\
\end{array}
$$
 Note that the image of $F$ lies
in $\cal{O}$.
By the above, the induced connection recovers   precisely the compatible
connection $D$ on $\cal{E}$.

% **************************************

\section{ Holomorphic torsion  and Hermitian metric variation:}
 \label{section.analyticholo2}

For $E$ a holomorphic bundle over the complex manifold $W$ with Hermitian metric
$g = <.,.>$, let $D$ be a linear connection on $E$ which is compatible
with the holomorphic structure and is of type $(1,1)$. Let $W$ have
complex dimension $n$, so real dimension $2n$.

Define the operator $\Delta_{E,\bar{\partial}}$ by
\begin{eqnarray}
\Delta_{E,\bar{\partial}} =( \bar{\partial}_E + \bar{\partial}^*_{E,D''})^2
   = \bar{\partial}_E \bar{\partial}^*_{E,D''} +
          \bar{\partial}^*_{E,D''}  \bar{\partial}_E
\end{eqnarray}
yielding a generally non-self adjoint operator, which commutes with both of
$\bar{\partial}_E ,  \bar{\partial}^*_{E,D''}$ which have squares equal to zero.

Choose a Hermitian metric, say $<.,.>_E$ on the complex bundle $E$.
As noted in the last section, the adjoint of $\partial_E$ differs by a bundle
mapping, say $\alpha$, from $\bar{\partial}^*_{E,D''}$,
\begin{eqnarray}
\bar{\partial}^{*}_{E,D''} = adjoint(\bar{\partial}_E) + \alpha \label{eqn.alpha}
\eeqn
and so in addition
\begin{eqnarray}
\Delta_{E, \bar{\partial}} = (\ \bar{\partial}_E +  adjoint(\bar{\partial}_E) )^2
 + \alpha \bar{\partial}_E +  \bar{\partial}_E \alpha
\eeqn
is an elliptic second order partial differential equation with
scalar symbol.

Consequently,  standard methods of elliptic theory apply,
e.g., the methods of Atiyah, Patodi, Singer and
Seeley, \cite{AtiyahPatodiSingerIII,Seeley1,Seeley2} to
prove the following basic results summarized here,
see section \ref{section.Analysis}. In modified form they appear
in the papers of Ray and Singer \cite{RaySinger1,RaySinger2}.

Firstly, the spectrum of $\Delta_{E,\bar{\partial}}$
is discrete and has generalized eigenspaces of
finite multiplicity.

\vspace{.2in}

As a second step towards understanding the spectral properties
of these operators,  recall an observation from  Atiyah-Potodi-Singer's
part III  \cite{AtiyahPatodiSingerIII}. There they  also study the eta invariant of the operator
$B_{ev}$ coupled to a flat non-unitary bundle. The following is just a more explicit version.

\begin{lemma} \label{lemma.AN1}
Fix a Hermitian inner product $<.,.>$ and use this to
define an inner product on forms with values in $E$ via
$<f,g> = \int_M \  <f(x) \wedge \star g(x)> \ dvol_x$.
Then in terms of the smooth bundle mapping $\alpha$ with
$\bar{\partial}^{*}_{E,D''} = adjoint(\bar{\partial}_E) + \alpha$, the spectrum of
$(\bar{\partial}_E + \bar{\partial}^*_{E,D''}) $ lies in the strip of $\{\lambda\}$ satisfying
$$
-||\alpha|| \le Im \ \lambda \ \le +||\alpha||.
$$

If $\alpha  \neq 0$,
 the spectrum of $\Delta_{E,\bar{\partial}}$ lies on or inside  the parabola enclosing
the positive x-axis: \  $\{ \lambda\}$
with
$$
Re \ \lambda \ge \frac{(Im \ \lambda)^2}{4 ||\alpha||^2} - ||\alpha||^2
$$
Of course, if $||\alpha||=0$, then the spectrum of  the operator
$\Delta_{E, \bar{\partial}} = (\ \bar{\partial}_E +  adjoint(\bar{\partial}_E) )^2$
is real non-negative.

\end{lemma}

Let $K>0$ be a real number which is not the real part of
any eigenvalue.
Let  $S(p,q,K)$
be a complete enumeration of all the generalized eigenvalues counted with
multiplicities with real part greater than $K$ of $\Delta_{E,\bar{\partial}}$
acting on $A^{p,q}(W,E)$. That is, $\Re (\lambda_j) >K$.

Then the zeta function
\begin{eqnarray}
\zeta_{p,q,E,K,g}(s) = \Sigma_{\lambda_j \in S(p,q,K)} \ \frac{1}{\lambda_j^s} \ ,
\end{eqnarray}
defined using the principle values  for  the complex powers,  converges
for $Re(s) > n/2$ with $n = dim_R  \ W = 2n$. If $\Pi_{E,K,g}$ denotes the spectral projection
on the span of the generalized eigenvectors with generalized eigenvalues
with real part less than $K$ and $Q_{E,K,g} = (1 - \Pi_{E,K,g})$, then the
elliptic operator $Q_{E,K,g} \ \Delta_{E,\bar{\partial}}$
fits the setting of Seeley \cite{Seeley1,Seeley2}. In particular, its complex
powers $[Q_{E,K,g} \Delta_{E,\bar{\partial}}]^{-s}$ are well defined, as in that
paper. Also as in \cite{Seeley2}, the ``heat kernel''
$e^{- t \ Q_{E,K,g} \ \Delta_{E,\bar{\partial}} }= Q_{E,K,g} e^{-t \ \Delta_{E,\bar{\partial}}}$ is well defined for
$t $  real positive, is of trace class, and for $Re(s) >N$ the
formula
\begin{eqnarray}
[ Q_{E,K,g} \ \Delta_{E,\bar{\partial}}]^{-s}
 = \frac{1}{\Gamma(s)}\int_0^{+\infty} \ t^{s-1} \ e^{- t \ Q_{E,K,g}
 \ \Delta_{E,\bar{\partial}}} \ dt \quad \quad
\end{eqnarray}
holds and its trace gives the zeta function $\zeta_{p,q,E,K,g}(s)$
for the bundle $E$.

Moreover, the methods of Seeley imply that $\zeta_{p,q,E,K,g}(s)$
has a meromorphic extension to the whole complex plane and
is analytic at $s=0$. Consequently, the derivative
at $s=0$, \ $\zeta'_{p,q,E,K,g}(0)$ is meaningful.

Following Ray and Singer \cite{RaySinger2}, one sets
$$
\begin{array}{l}
Ray-Singer-Term(p,E, K,g) \\
 \hspace{1in} = exp((1/2) \  \Sigma_{q=0}^n \ (-1)^q \ q \ \zeta'(p,q,E,K,g)(0))
\end{array}
$$

The methods of Ray and Singer again apply to prove the following lemma
about  the Ray--Singer terms  for $E$,
for details  see section \ref{Analysis}. Here the metric
dependence is explicitly recorded in the operators.

\begin{lemma} \label{lemma.3A}
Fix $K >0$. Let $g(u), \  a \le u \le b$,  be a smooth family
of Hermitian metrics on $TW$ such that the operators
$\Delta_{E,\bar{\partial},g(u)}$
 have no generalized
eigenvalues with real parts equal
to $K$.

Let the star operator $\star_{u,p,q} $ be the bundle
isomorphism associated to the Hermitian inner product
$g(u)$,
\begin{eqnarray}
\star_{u,p,q} :  \Lambda^p T^{1,0}W \otimes
\Lambda^q T^{0,1}W \stackrel{\star_u}{\rightarrow}  \Lambda^{n-p} T^{1,0}W \otimes
\Lambda^{n-q} T^{0,1}W\ ,
\end{eqnarray}
and similarly for $\hat{\star}_u$.
and set $\alpha_{u,p,q} = (\star_{u,p,q})^{-1} \ d/du \ ( \star_{u,p,q})= \newline
(\hat{\star}_{u,p,q})^{-1} \ d/du \ ( \hat{\star}_{u,p,q}) $,
 a the self mapping
of the bundle $ \Lambda^p T^{1,0}W \otimes
\Lambda^q T^{0,1}W$ depending on $u$.

As is the case,  by the results of Seeley,
let the trace $Tr(\alpha_{u,p,q} e^{-t \Delta_{E,\bar{\partial}}})$ have an asymptotic expansion for small $t$
given by
$$
Tr(\alpha_{u,p,q} e^{-t \Delta_{E,\bar{\partial}}}) = \Sigma_{j=0}^M a_{j,u,p,q} \ t^{-N/2+j} + 0 ( t^{-N/2+ M+1})
$$
with $N = dim \ W = 2n$ and $M > N/2 + 1$. Here  $a_{u,p,q}$ is given by explicit local formulas $a_{u,p,q} = \int_W \ b_{u,p,q}$
with the forms $b_{u,p,q}$ expressed directly, locally, in  terms of the Hermitian metric
$g(u)$, the connection form of $D$, and their covariant derivatives. In particular,  the coefficient
of $t^0$ is $a_{N/2,u,p,q}$.

Let $\Pi_{E,K,p,q,u}$ be the spectral projection onto  the span
of the generalized eigenvectors   of $\Delta_{E,\bar{\partial},g(u)}$
acting on $(p,q)$ forms
with generalized eigenvalues of real part less than $K$.

Then
$$
\begin{array}{l}
Ray-Singer-Term(p,E, K,g(u)) \\
= exp((1/2) \ \Sigma_{q=0}^n \ (-1)^q \ q \ \zeta'(p,q,E,K,g(u))(0))
\end{array}
$$
  varies smoothly with $u$, for $ a \le u \le b$

 Moreover,
$$
\begin{array}{ll}
d/du [
log \  ( \ Ray-Singer-Term(p,E,K,g(u) ) \ ]  \\
 = (1/2) \Sigma_{q=0}^n \ (-1)^q \ [Tr (\Pi_{E,K,p,q, g(u)} \  (\alpha_{u,p,q} \otimes Id_E))
+ \int_W \ b_{n,u,p,q}]
\end{array}
$$

\end{lemma}

\vspace{.3in}
Similarly, the methods of Ray and Singer straight-forwardly give the expected dependence
of the $Ray-Singer-Term(p,E, K,g)$
on changing $K>0$. It is:

\begin{lemma} \label{lemma.3BB}
Fix $L > K >0$. Let $g$  be a Hermitian metric
on $TW$ such that the  operator $\Delta_{E,\bar{\partial},g}$ acting
on $A^{p,*}(W,E)$ has  no  eigenvalues with real part equal to
$K$ or $L$. Let $\{ \lambda_{p,q,j}| j = 1\cdots n_q\}$ be a complete
enumeration counting with multiplicities
of the generalized eigenvalues of $\Delta_{E,g}$ acting on $A^{p,q}(W,E)$
which have real part  in the range $K$ to $L$.
Then the formula below holds:
\begin{eqnarray}
\left[ \frac{Ray-Singer-Term(p,E, K,g)}{Ray-Singer-Term(p,E, L,g)} \right]^{\, 2}
= \left[ \prod_{q=0}^n \ (\prod_{j=1}^{n_q} \  \lambda_{p,q,j})^{(-1)^q  \, q} \right]^{-1} \nonumber
\end{eqnarray}
\end{lemma}
For details see  section \ref{Analysis}.

\vspace{.3in}
Let $C^{p,q}(E,K,g)$ denote the span in $A^{p,q}(W,E)$
of the generalized eigensolutions of $\Delta_{E,\bar{\partial},g}$
 with generalized eigenvalues with real part
less than $K$. By elliptic theory, $C^{p,q}(E,K,g)$ is a finite dimensional
subspace of smooth sections.

Since $\Delta_{E,\bar{\partial}}$ commutes with $\bar{\partial_E}, \bar{\partial}^*_{E,D''}$
which each have square zero, by $D$ a type $(1,1)$ connection,
the graded  complex
\begin{eqnarray}
C^{p,\star}(E,K,g) := \bigoplus_{q=0}^n  \ C^{p,q}(E,K,g)
\nonumber
\end{eqnarray}
 has  two differentials, $ d, d^*$,  $d= \bar{\partial}_E$ and $d^* = \bar{\partial}^*_{E,D''} =
-(\hat{\star} \otimes Id_E)
 D'' (\hat{\star} \otimes Id_E)$ with $ d :
C^{p,q}(E,K,g) \stackrel{\bar{\partial}_E}{\rightarrow } C^{p,q+1}(E,K,g) $  increasing the  grading
by one  and $
d^* :
C^{p,q}(E,K,g) \rightarrow  C^{p,q-1}(E,K,g) $
 decreasing  degree by one.
That is, one has the same complex equipped with two differentials:
$C^{p,\star}(E,K,g), d, d^*$
of degrees $+1,-1, $ respectively.

\vspace{.3in}
In the next section,  \S \ref{section.algebraictorsions}, it will be proved that
in the general  algebraic setting of a  graded complex of length $n$,\ $C^*$
equipped with two differentials, $d$ of degree $+1$ and $d^*$ of degree $-1$,
i.e, $(C^*,d, d^*)$, there is
a natural non-vanishing algebraic torsion invariant
$$
torsion( C^*,d,d^*) \in (det
\ H^*(C^*,d) ) \otimes (det \ H_*(C^*,d^*))^{(-1)} \ .
$$
Here $H^q(C^*,d) = (ker \ d:C^q \rightarrow C^{q+1})/(im \ d: C^{q-1} \rightarrow C^q)$,
the cohomology, and
$H_q(C^*,d^*) = (ker \ d^*:C^q \rightarrow C^{q-1})/(im \ d^*: C^{q+1} \rightarrow C^q)$,
the homology.

 Hence, one gets a non-vanishing torsion invariant
$$
\begin{array}{l}
torsion((C^{p,*}(E,K,g),\bar{\partial},  \bar{\partial}^*_{E,D''}) ) \\
\hspace{.4in} \in (det
\ H^*( C^{p,*}(E,K,g),\bar{\partial}_E ) \otimes (det \
  H_*(C^{p,*}(E,K,g), \bar{\partial}^*_{E,D''})^{(-1)}
\end{array}
$$
These cohomologies and homologies can be identified as follows:

Let $M$ be the multiplicity of the generalized eigenvalue $0$
for $\Delta_{E,\bar{\partial}}$
and $\cal{H}$ denote the span of these generalized eigenvectors
with generalized eigenvalue $0$. By spectral theory $\cal{H}$ is finite
dimensional and consists of smooth sections.  Then, by adapting standard Hodge--theoretic methods,
 one has the
direct sum decompositions:
$$
\begin{array}{l}
A^{*,*}(W,E) =   {\cal{H}}
 \oplus \bar{\partial}_E  [(\Delta_{E,\bar{\partial}_E})^M \
A^{*,*}(W,E)] \\ \hspace{1.5in} \oplus  ( (\hat{\star} \otimes Id_E) \partial_E
(\hat{\star} \otimes Id_E) ) [ (\Delta_{E,\bar{\partial}})^M \ A^{*,*}(W,E)]
\end{array}
$$
and $\bar{\partial }_E$ maps the second summand isomorphically to
the first. By this means one can easily show that the inclusion
of co-chain complexes
$$
(C^{p,*}(E,K,g),\bar{\partial}_E )
\subset (A^{p,*}(W,E),\bar{\partial}_E)
$$
induces an isomorphism on cohomology.
That is,
$H^q( C^{p,*}(E,K,g),\bar{\partial}_E )  \cong H_{\bar{\partial}}^{p,q}(W,E)$.
 In particular, $det \ H^*( C^{p,*}(E,K,g),\bar{\partial}_E )$  $
 \cong det \ H_{\bar{\partial}_E}^{p,*}(W,E)$.

Similarly, by $(\hat{\star} \otimes Id_E)^2|A^{p,q}(W,E) = (-1)^{(p+q)(2n-(p+q))}$, it follows that
$(\hat{\star} \otimes Id_E) (-(\hat{\star} \otimes Id_E) D'' (\hat{\star} \otimes Id_E))
 (\hat{\star} \otimes Id_E)^{-1} = \pm \ D''$;
so $(\hat{\star} \otimes Id_E)$ induces a  complex linear isomorphism
of the graded complex $(C^{n-*,n-p}(E,K,g),  \pm D'')$
to the complex $(C^{p,*}(E,K,g),d^*)$ sending $C^{n-q,n-p}(W,E) \rightarrow C^{p,q}(W,E)$
 This identifies
$ H_q(C^{p,*}(E,K,g), d^* ) ) = H_{D''}^{n-q,n-p}(W,E)$
and thus the determinants
$ det \ H_q(C^{p,*}(E,K,g), \bar{\partial}^*_{E,D''} ) ) = det \ H_{D''}^{n-q,n-p}(W,E) $.

Hence, in this special situation the algebraic invariants  may be regarded
as an element of the same complex line,
\begin{eqnarray}
torsion((C^{p,*}(E,K,g),d,
   d^*), \nonumber  \\
 \in \  det  \ H_{\bar{\partial}}^{p,*}(W,E)) \otimes
         [det \ H_{D''}^{n-*,N-p}(W,E) ]^{-1} \nonumber \\
\cong  [det  \ H_{\bar{\partial}}^{p,*}(W,E))] \otimes
         [det \ H_{D''}^{*,n-p}(W,E) ]^{(-1)^{n+1}}
\end{eqnarray}
for any real number $K >0$ and any choice of Hermitian metric $g$ on $TW$.
This last line is independent of the choice of Hermitian metric $g$.

Now let $L>K>0$ be real and positive.
It will follow from the algebraic lemmas B,C of \S \ref{section.algebraictorsions},
that under the assumptions of lemma \ref{lemma.3BB},  the algebraic
torsion of
$( C^{p,*}(E,K,g), d, d^*)$
and that for $L$ with  $K>L>0$ are precisely related by the eigenvalues $\lambda_{q,j}$
of lemma \ref{lemma.3BB} as follows:
\begin{eqnarray}
 torsion( C^{p,*}(E,L,g), d, d^*)) \nonumber \\
=   torsion ( C^{p,*}(E,K,g), d, d^*))  \label{eqn.stabilize1}\\
\hspace{.5in} \bullet \left[ \prod_{q=0}^n \ (\prod_{j=1}^{n_q} \  \lambda_{q,j})^{(-1)^q  \, q}
\right]^{-1}
\nonumber
\end{eqnarray}
 as  elements of  $\  [det  \ H_{\bar{\partial}_E}^{p,*}(W,E)] \otimes
         [det \ H_{D''}^{*,n-p}(W,E) ]^{(-1)^{N+1}}$.

Similarly, it will follow from lemma \ref{lemma.5vary} of \S \ref{section.Analysis}
 for a smooth
family of Hermitian metrics $g(u)$ satisfying the  assumptions of lemma \ref{lemma.3BB}
that one gets the equality
\begin{eqnarray}
d/du \left[ \ log  \ torsion(C^*_p(E,K,g(u)), d, d^*)
\right] \nonumber \\
  = - \ Tr (\Pi_{E,K,g(u)} \  (\alpha_{g(u)} \otimes Id_E)) \quad \quad .
\end{eqnarray}

As an immediate consequence of the above two equations ,
the main theorem below follows:

\begin{thm}[Variation theorem of holomorphic torsion ] \label{thm.MAIN}

For any choice of Hermitian inner product $g$ on
$TW$, pick a real number $K>0$ for which the
operators $ \Delta_{E,\bar{\partial},g}$ acting on $A^{*,*}(W,E)$ have
no eigenvalue with real part equal to $K$. Define
the graded complex $C^{p,*}(E,K,g) = \bigoplus_{q=0}^N  \ C^{p,q}(E,K,g)$ with its
two differentials, $d,d^*$ as above. Then
the combination
$$
\begin{array}{lll}
\tau_{holo,p}( W,E) \  := & \
\left[ \ torsion((C^{p,*}(E,K,g),d,d^*)
 \right] \\
& \hspace{.5in} \bullet
\left[  Ray-Singer-Term(p,E,K,g) \right]^{\, 2}\\
\end{array}
$$
is independent of the choice of $K>0$.

By definition $\tau_{holo,p}( W,E)$ is a non-vanishing element
of the determinant line bundle
\begin{eqnarray}
\tau_{holo,p}(W,E) \  \in  [\ det \ H_{\bar{\partial}}^{p,*}(W,E)]  \otimes
\ [  det \ H_{D''}^{*,N-p}(W,E)]^{(-1)^{N+1}} . \nonumber  \\
\end{eqnarray}

For a smooth variation of Hermitian metrics $g(u), a \le u \le b$, \newline
$\tau_{holo,p}(W,E,g(u))$ varies smoothly. Here the explicit dependence
on the metric is recorded. Also with the notation of lemma
\ref{lemma.3A}, $d/du \ \tau_{holo,p}(W,E,g(u))$ is given by
the local formula:
$$
d/du \ \tau_{holo,p}(W,E,g(u))
= \Sigma_{q=0}^N \ (-1)^q \ \int_W \ b_{n/2,u,p,q}
$$
with $n = dim_R \ W = 2N$.

In the acyclic cases it is a complex number. In the case that
$E$ is acyclic and flat unitary, this torsion is a real number.
Since the operator $\Delta_{E,\bar{\partial}}$ is self-adjoint
in this case, one may take $K>0$ less than the smallest non-vanishing
real eigenvalue; for this choice of $K>0$, the algebraic correction term
is just $+1$ and the above reduces to the formula implicit in  Ray and Singer
\cite{RaySinger2}.

\end{thm}

This theorem has a corollary which generalizes the explicit theorem
of Ray and Singer on holomorphic torsion \cite{RaySinger2}.

\begin{corollary}
If $E$ is a holomorphic bundle with  connection $D$ which is compatible and type $(1,1)$
over a Hermitian complex manifold $W$ and $F_1,F_2$ are two flat complex bundles
over $E$ of the same dimension,  then the ``quotient'' of torsions
$$
\tau_{holo,p}(W,E \otimes F_1,g(u))\otimes [\tau_{holo,p}(W,E \otimes F_2,g(u))]^{-1}
$$
in the tensor product of determinant line bundles
$$
\begin{array}{l}
 ([ \ det \ H_{\bar{\partial}}^{p,*}(W,E\otimes F_1)]  \otimes
\ [  det \ H_{D''}^{*,N-p}(W,E \otimes F_1)]^{-1} ) \\
 \otimes ([ \ det \ H_{\bar{\partial}}^{p,*}(W,E\otimes F_2)]  \otimes
\ [  det \ H_{D''}^{*,N-p}(W,E\otimes F_2)]^{-1})^{-1}\\
\end{array}
$$
is independent of the Hermitian metric $g(u)$ chosen.

\end{corollary}

This corollary follows since the two bundles $E \otimes F_1, E \otimes F_2$
are locally identical as bundles with connections,  so the local corrections
$b_{n/2,u,p,q}$ are the same and  cancel in computing the variation
of the holomorphic torsions.

\section{An algebraic torsion for complexes with boundary and coboundary:}  \label{section.algebraictorsions}

In this section the torsion of a chain complex $C^*$ over the complex numbers
equipped with two differentials, $d$ going up with  $d^2=0$, and $d^*$ going down  with $(d^*)^2 = 0$, is defined.
On the one hand $(C^*,d)$ is a cochain complex,$d:C^q \rightarrow C^{q+1}, d^2=0$, and also $(C^*,d^*)$ is a chain
complex, $d^*: C^q \rightarrow C^{q-1}, (d^*)^2=0$, so the cohomologies of $d$, $ H^q(C^*,d)$,
and the homologies of $d^*$, $H_q(C^*,d^*)$, can be formed.
The torsion of $(C^*,d,d^*)$  is to be an element of the product of
determinant lines
$$
torsion(C^*,d,d^*) \in det(H^*(C^*,d)) \otimes det(H_*(C^*,d^*))^{-1}
$$
formed from the cohomology and homology of $C^*$.

Here a complex one dimensional vector space $L$ is called a complex line.
The dual of a complex line $L$ is denoted by $L^{-1}$.
This torsion will be  shown to
have good algebraic properties and the relation to other
standard algebraic definitions is explained. It is crucial
that certain sign difficulties are dealt with carefully,
so the definition involves a suitable choice of signs.

\vspace{.2in}

For $E$ a k  dimensional vector space over the complex numbers, $C$,
set
$$
det(E) = \Lambda^{k}(E)
$$
, the complex line given by the highest exterior
product, $k= dim \ E$.
By definition, if ${\bf{v}} = \{v_1,\cdots,v_k\}$ is an ordered  basis of $E$, then
the ordered wedge product ${\Lambda \bf{v}} := v_1 \wedge v_2 \wedge \cdots v_k$ is a basis for $det(E)$ and
if  ${\bf{w}} = \{w_1,\cdots,w_k\}$ is another ordered basis, then
$\Lambda \bf{w} = [\bf{w}/ \bf{v}] \cdot \Lambda \bf{v}$ where the complex number $[\bf{w}/ \bf{v}]$
 is the determinant of the $k$ by $k$ matrix expressing the ordered basis
 $\bf{w}$ in terms of the ordered basis $\bf{v}$.

By convention set $det(E) = C$
for the 0-dimensional vector space $E=\{0\}$
and take its standard generator  to be $1 \in C$.

\vspace{.2in}
Similarly, if $(C^*,d)$ is a finite dimensional cochain complex of length $N$
\begin{eqnarray}
0 \rightarrow C^0 \stackrel{d}{\rightarrow} C^1   \stackrel{d}{ \rightarrow} \cdots
 \rightarrow \cdots  \stackrel{d}{\rightarrow} C^N \rightarrow 0
\end{eqnarray}
so $d^2 = 0$,
set
$$
det(C^*) = \bigotimes_{q=0}^{q=N} (\ det \ C^q)^{(-1)^q} \ .
$$
In this context, one defines the q-coboundries, $B^q = dC^{q-1}$, the q-cocyles,
$Z^q = ker \ d:C^{q} \rightarrow C^{q+1}$, and the q-cohomologies
$H^q(C^*,d) = Z^q/B^q$ and the also associated determinant line
$$
det(H^*(C^*,d)) = \bigotimes_{q=0}^{q=N} (\ det \ H^q(C^*,d))^{(-1)^q} \ .
$$

\vspace{.3in}

In this terminology, following the conventions of Bismut and Zhang \cite{BismutZhang1},
the torsion isomorphism of the graded cochain complex
$(C^*,d)$ is a standard  isomorphism of determinant
lines:
\begin{eqnarray}
\tau : det(C^*)  \stackrel{\cong}{\rightarrow} \ det( H^*(C^*,d))
\end{eqnarray}
recalled  below.

Now if in addition, $C^*$ has another differential $d^*: C^q \rightarrow C^{q-1}$
decreasing dimension and $(d^*)^2 = 0$,
\begin{eqnarray}
0 \leftarrow C^0 \stackrel{d^*}{\leftarrow} C^1   \stackrel{d^*}{ \leftarrow} \cdots
 \leftarrow \cdots  \stackrel{d^*}{\leftarrow} C^N \leftarrow 0,
\end{eqnarray}
then correspondingly one  has a standard
isomorphism of determinant lines:
\begin{eqnarray}
\tau' : det(C^*)  \stackrel{\cong}{\rightarrow} \ det( H_*(C^*,d^*))
\end{eqnarray}
defined also below. Here $H^q(C^*,d^*) = Z_q(C^*)/B_q(C^*)$ is the q-homology
groups of $d^*$. By definition the q-cycles are $Z_q(C^*) = ker \ d^*: C ^q \rightarrow C^{q-1}$
and the q-boundaries are $B_q(C^*) = d^*C^{q+1}$.

\vspace{.3in}
For the chain complex with these two differentials,$(C^*,d,d^*)$, there is then defined the
induced isomorphism of complex lines:
$$
det(H_*(C^*,d)) \stackrel{(\tau')^{-1}}{\rightarrow} det(C^*) \stackrel{\tau}{\rightarrow} det(H^*(C^*,d))
$$
which then defines a element of the product of determinant lines:
$$
\tau(C^*,d,d^*) \in  det(H^*(C^*,d)) \otimes  (det(H_*(C^*,d^*)))^{-1}
$$

\begin{flushleft}

  \textbf{Definition of the algebraic torsion of  the complex with two differentials: $ (C^*,d,d^*) $   }

\end{flushleft}

One sets
$$
\begin{array}{l}
torsion(C^*,d,d^*) \\
 = (-1)^{S(C)}  \tau(C^*,d,d^*) \in  det(H^*(C^*,d)) \otimes  (det(H_*(C^*,d^*)))^{-1}
\end{array}
$$
Here  the sign $(-1)^{S(C^*)}$ is defined by setting
\begin{eqnarray} \label{eqn.sign}
 S(C^*) = &
\Sigma_q \ [dim_C \ B_{q-1}(C^*) \cdot dim_C \ B^{q+1}(C^* ) \nonumber \\
&  + dim_C \ B^{q+1}(C^*)  \cdot dim_C H_q(C^* ,\partial)   \\
&  +
dim_C \ B_{q-1}(C^*)  \cdot dim_C H^{q}(C^* ,\partial^*)]. \nonumber
\end{eqnarray}
The mysterious sign is chosen so that this algebraic invariant
has the desirable properties A,B,C below.

\vspace{.2in}
To specify the isomorphism $\tau$ for the cochain complex $C^*,d$,
fix a  direct sum decompositions of  vector spaces
\begin{eqnarray}
C^q = B^q \oplus H^q \oplus A^q
\end{eqnarray}
where $B^q$ are the coboundaries, $B^q \oplus H^q$ are
the cocyles, and $A^q$  fills out the rest of $C^q$. Note $A^N = \{0\}$.
 Then by choice $H^q$ is naturally
isomorphic to the cohomology, $H^q \cong H^q(C^*,d)$ under
the projection of $Z^q$ onto $H^q(C^*,d)$.

Taking the wedge product in this order establishes a
natural isomorphism
\begin{eqnarray}
\Psi: (det \ B^q) \otimes (det \ H^q) \otimes (det \ A^q) \cong det \ C^q
\end{eqnarray}

Fix generators $c(q) \in det \ C^q, a(q) \in det \ A^q$ for each
$q=0,\cdots ,N$. Use
 the convention that for the 0-dimensional subspace,
$\{0\}$,  $det \ \{0\} =C$ by definition and the chosen generator
is  $1 \in C$.

Let $d(a(q))$  denote the generator of $det \ B^{q+1}$ determined
by the isomorphism $d : A^q \rightarrow B^{q+1}$. Then for each
$q=0,\cdots,N$ there is a unique element $h(q) \in det \ H^q$ such that
\begin{eqnarray}
c(q) =  \Psi(\  d(a(q-1)) \ \otimes  \ h(q) \ \otimes \ a(q) \ )
\end{eqnarray}

For $x \in L$ a non-vanishing element of the one dimensional
vector space $L$,
let  $x^{-1}$ be the unique element of $L^{-1} = L^*$ which
maps $x$ to $+1$.

With these choices, the algebraic torsion mapping
\begin{eqnarray}
\tau : det(C^*)  \stackrel{\cong}{\rightarrow} \ det H^*(C^*,d) \label{eqn.tau}
\end{eqnarray}
is defined as the isomorphism which maps
$$
 c(0)\otimes c(1)^{-1} \otimes \cdots \otimes c(N)^{(-1)^N}
\mapsto h(0) \otimes h(1)^{-1} \otimes \cdots \otimes h(N)^{(-1)^N}
$$
It is a foundational theorem that the resultant mapping
is independent of all choices \cite{Milnor1, BismutZhang1}.

[Note that Braverman and Kappeller \cite{BravermanKappeller2}
introduce an additional sign at this point in their definition
of a related  isomorphism. For
simplicity this is  not done here. The cost of this simplicity is that one must   introduce
the total sign correction $(-1)^{S(C)}$ in the above
definition of $torsion(C^*,d,d^*)$ to reestablish  desirable functional properties.]

\vspace{.3in}
More concretely and explicitly,
let $\{h^{q,k}| k=1,\cdots, v_q \}$, \ $ v_q = dim_C \ H^q$ be a chosen ordered basis for
$H^q(C^*,d)$ for each $q=0,\cdots, N$ and
${\bf{c}}^q = \{c^{q,1},c^{q,2},\cdots, c^{q,n_q}\}, n_q = dim_C \ C^q$, be a chosen
ordered basis for $C^q$ in each degree, $q=0,\cdots, N$.
Let $c(q)$ be the ordered wedge product, $\wedge_j \ c^{q,j} \in det \ C^q$,
and $h(q)$ be the ordered wedge product $\wedge_k \ h^{q,k} \in det \ H^q$.

Following the notation of Ray and Singer \cite{RaySinger1}, chose an ordered basis, say
${\bf{b}}^q = \{ b^{q,1}, b^{q,2},\cdots, b^{q,s_q} \}$,
for the coboundaries $B^q = Im \ d: C^{q-1} \rightarrow C^q$  in $C^q$.
Now chose
$\tilde{ {\bf{b}} }^q = \{ \tilde{b}^{q,1}, \tilde{b}^{q,2},\cdots,
\tilde{b}^{q,s_q}\}$  in $C^{q-1}$ so that $d \tilde{b}^{q,i} = b^{q,i}$ in $B^q$
for $i = 1,\cdots, s_q$. Thus taking the ordered elements $\tilde{ {\bf{b}} }^q$
as a basis for $A^{q-1}$ and $ {\bf{b}}^q$ as a ordered basis for $B^q$, one has
$d( \wedge_j \tilde{b}^{q,j})
= (\wedge_j b^{q,j}) \in det \ B^q$.

Now $Z^q/B^q \cong H^q$, so pick elements $lh^{q,k} \in Z^q \subset C^q$
projecting to the basis elements $h^{q,k}$ of $H^q$. Taking
the ordered set ${\bf{lk}}^q := \{lh^{q,k}|k=1,\cdots, v_k\}$ as a basis for
$H^q \subset C^q$ gives the above direct sum decomposition of $C^q$,
since with these choices $\{ {\bf{b}}^q,{\bf{lh}}^q,\tilde{\bf{b}}^{q+1}\}$
is an ordered basis for $ C^q$. This may  be compared to the
chosen basis ${\bf{c}}^q$.

In this notation one has the translation of invariants:
\begin{eqnarray}
\tau(c(0)\otimes c(1)^{-1} \otimes \cdots \otimes c(N)^{(-1)^N}) &
  \\
  = \left[ \prod_{q=0}^N \ [ {\bf{b}}^q,{\bf{lh}}^q ,
\tilde{ {\bf{b}} }^{q+1} / {\bf{c}}^q  ]^{(-1)^q  } \right]^{-1} &  \nonumber  \\
 \times ( h(0) \otimes h(1)^{-1} \otimes \cdots \otimes h(N)^{(-1)^N}) &
 \nonumber
\end{eqnarray}

\vspace{.3in}
In the below  this concrete point of view is often taken. For a choice
of bases for $C^*$ and for $H^*$ one  forms  the above complex number \newline
$\prod_{q=0}^N \ [ {\bf{b}}^q,\textbf{lh}^q ,
\tilde{ {\bf{b}} }^{q+1} / {\bf{c}}^q  ]^{(-1)^q  }$. For the above choices
one defines
\begin{eqnarray}
\tau(C^*,d, \{ {\bf{h}}_q\}, \{ {\bf{c}}_q\})
:= \left[ \prod_{q=0}^N \ [ {\bf{b}}^q,{\bf{lh}}^q ,
 \tilde{ {\bf{b}} }^{q+1} / {\bf{c}}^q  ]^{(-1)^q  } \right]^{-1} \label{eqn.tau4}
\end{eqnarray}
regarding it as specifying the isomorphism
$\tau: det(C^*) \rightarrow det(H^*(C^*,d))$  above. A simple check shows
$ \tau(C^*,d, \{ {\bf{h}}_q\}, \{ {\bf{c}}_q\})$
is independent of all choices but those
of the ordered bases $\{ {\bf{c}}^q\}, \{ {\bf{h}}^q\}$.

\vspace{.2in}

In analogy to the isomorphism  \ref{eqn.tau}
there is a natural isomorphism of determinant
lines defined by the chain complex $(C^*,d^*)$:
\begin{eqnarray}
\tau' : det \ C^*  \rightarrow det \ H_*(C^*,d^*) \label{eqn.tau'} \ .
\end{eqnarray}
It is defined as follows:

To specify the isomorphism $\tau'$ for the cochain complex $C^*,d$,
fix a  direct sum decompositions of  vector spaces
\begin{eqnarray}
C^q = B_q \oplus H_q \oplus A_q
\end{eqnarray}
where $B_q= d^* C^{q+1}$ are the boundaries, $B_q \oplus H_q$ are
the cocyles, and $A_q$  fills out the rest of $C_q$. Note $A_0 = \{0\}$.
 Then by choice $H_q$ is naturally
isomorphic to the homology, $H_q \cong H_q(C^*,d^*)$ under
the projection of $Z_q$ onto $H_q(C^*,d^*)$.

Taking the wedge product in this order establishes a
natural isomorphism
\begin{eqnarray}
\Psi': (det \ B_q) \otimes (det \ H_q) \otimes (det \ A_q) \cong det \ C^q
\end{eqnarray}

Fix generators $c(q) \in det \ C^q, a'(q) \in det \ A_q$ for each
$q=0,\cdots ,N$. Use
 the convention that for the 0-dimensional subspace,
$\{0\}$,  $det \ \{0\} =C$ by definition and the chosen generator
is  $1 \in C$.

Let $d^*(a'(q))$  denote the generator of $det \ B_{q-1}$ determined
by the isomorphism $d^* : A_q \rightarrow B_{q-1}$. Then for each
$q=0,\cdots,N$ there is a unique element $h'(q) \in det \ H_q$ such that
\begin{eqnarray}
c(q) =  \Psi'(\  d^*(a'(q+1)) \ \otimes  \ h'(q) \ \otimes \ a'(q) \ )
\end{eqnarray}

With these choices, the algebraic torsion mapping
\begin{eqnarray}
\tau' : det(C^*)  \stackrel{\cong}{\rightarrow} \ det( H_*(C^*,d)) \label{eqn.tau'}
\end{eqnarray}
is defined as the isomorphism which maps
$$
\begin{array}{l}
 c(0)\otimes c(1)^{-1} \otimes \cdots \otimes c(N)^{(-1)^N} \\
\mapsto h'(0) \otimes h'(1)^{-1} \otimes \cdots \otimes h'(N)^{(-1)^N}
\end{array}
$$

\vspace{.2in}

Concretely, having already chosen bases ${\bf{c}}^q$ for
$C^q$, let ${\bf{b}}_q = \{ b_{q,1}, b_{q,2},\cdots, b_{q,r_q} \}$
be an ordered basis
for the boundaries $B_q = Im \ d^*: C^{q+1} \rightarrow C^q$  in $C^q$.
Now chose
$\tilde{ {\bf{b}} }_q = \{ \tilde{b}_{q,1}, \tilde{b}_{q,2},\cdots,
\tilde{b}_{q,r_q}\}$  in $C^{q+1}$ so that $d^* \tilde{b}_{q,i} = b_{q,i}$ in $B_q$
for $i = 1,\cdots, r_q$. Thus taking the ordered elements $\tilde{ {\bf{b}} }_q$
as a basis for $A_{q+1}$ and $ {\bf{b}}_q$ as a ordered basis for $B_q$, one has
$d^* ( \wedge_j \tilde{b}_{q,j})
= (\wedge_j b_{q,j}) \in det \ B_q$ under the isomorphism $d^* :A_{q+1} \rightarrow B_q$.

Now $Z_q/B_q \cong H_q$, where $Z_q$ are the q-cycles, $Z_q = ker \ d^*: C^q
\rightarrow C^{q-1}$.  Pick elements $lh_{q,k} \in Z_q\subset C^q $
projecting to the chosen ordered basis elements $h'_{q,k}$ of $H_q$. Taking
the ordered set ${\bf{lk}}_q := \{lh_{q,k}|k=1,\cdots, u_k\}$ with $u_q = dim_C \ H_q(C^*,d^*)$
 as a basis for
$H_q \subset C^q $ gives the  direct sum decomposition, $C^q = B_q \oplus H_q \oplus A_{q-1}$,
since with these choices $\{ {\bf{b}}_q,{\bf{lh}}_q,\tilde{\bf{b}}_{q-1}\}$
is an ordered basis for $C^q$. This can be compared to the
chosen basis ${\bf{c}}^q$.

For the above choices
one sets
\begin{eqnarray}
\tau'(C^*,d^*, \{ {\bf{c}}^q\}, \{ {\bf{h'}}_q\})
= \left[ \prod_{q=0}^N \ [ {\bf{b}}_q,{\bf{lh}}_q ,
 \tilde{ {\bf{b}} }_{q-1} / {\bf{c}}^q  ]^{(-1)^q  }\right]^{-1} \label{eqn.tau5}
\end{eqnarray}
. A simple check shows it
is independent of all choices but those
of $\{ {\bf{c}}_q\}, \{ {\bf{h}}_q\}$.

The desired isomorphism $\tau'$ is then the element of \newline
$Hom( det(C^*), det( H_*(C^*,d^*)))$ defined by
$$
\begin{array}{l}
\tau'(c(0)\otimes c(1)^{-1} \otimes \cdots \otimes c(N)^{(-1)^N})
  \\
 = \tau'(C^*,d^*, \{ {\bf{c}}^q\}, \{ {\bf{h'}}_q\})
\times ( h'(0) \otimes h'(1)^{-1} \otimes \cdots \otimes h'(N)^{(-1)^N})
\end{array}
$$

In terms of the above notation, if $(C^*,d,d^*)$ has the
two differentials $d,d^*$, \ then for choices of ordered
bases $\{ {\bf{c}}^q\}$, $\{ {\bf{h}}^q\}$, and
$\{ {\bf{h}}_q\}$ for $C^q$, $H^q(C^*,d)$, and $H_q(C^*,d^*)$,
respectively,
one may form the quotient
$$
\tau(C^*,d, \{ {\bf{c}}^q\}, \{ {\bf{h}}^q\})/ \tau'(C^*,d^*, \{ {\bf{c}}^q\}, \{ {\bf{h'}}_q\}).
$$
Consequently one has a natural invariant
$torsion(C^*,d,d^*, \{ {\bf{h}}^q\}, \{ {\bf{h'}}_q\})$
 defined in several
equivalent manners:

\begin{eqnarray}
torsion(C^*,d,d^*, \{ {\bf{h}}^q\}, \{ {\bf{h'}}_q\}) \nonumber \\
:= (-1)^{S(C^*)} \ \tau(C^*,d, \{ {\bf{c}}^q\}, \{ {\bf{h}}^q\})/ \tau'(C^*,d^*, \{ {\bf{c}}^q\}, \{ {\bf{h}}_q\})
\nonumber  \\
= (-1)^{S(C^*)} \  \prod_q \left[ \frac{ [ {\bf{b}}^q,{\bf{lh}}^q ,
\tilde{ {\bf{b}} }^{q+1} / {\bf{c}}^q  ]}
{[ {\bf{b}}_q,{\bf{lh}}_q ,
\tilde{ {\bf{b}} }_{q-1} / {\bf{c}}^q  ]} \right]^{(-1)^{q+1}  }  \nonumber \\
=(-1)^{S(C^*)} \  \ \prod_q  [ {\bf{b}}^q,{\bf{lh}}^q ,
\tilde{ {\bf{b}} }^{q+1} /
 {\bf{b}}_q,{\bf{lh}}_q ,
\tilde{ {\bf{b}} }_{q-1}  ]}{^{(-1)^{q+1}  } \label{eqn.algetor'}  \nonumber
\end{eqnarray}
The last makes the independence of the choice of basis
for $C^q$ clear, so the definition depends only on $(C^*,d,d^*)$
and the choice of bases for the cohomology and homologies.

Here  the sign $(-1)^{S(C^*)}$ is defined as in equation \ref{eqn.sign}.

\vspace{.2in}
As the above makes clear, the non-zero complex number \newline
$torsion(C^*,d,d^*, \{ {\bf{h}}^q\}, \{ {\bf{h}}_q\})$
determines the torsion $torsion(C^*,d,d^*)$
by the formula:
$$
\begin{array}{l}
torsion(C^*,d,d^*)(\Lambda \{ {\bf{h}}_q\}) \\
= torsion(C^*,d,d^*, \{ {\bf{h}}^q\}, \{ {\bf{h'}}_q\})
\cdot (\Lambda \{ {\bf{h}}_q\})
\end{array}
$$
an element of $Hom(det(H_*(C^*,d^*), det(H^*(C^*,d)) )$.

\begin{flushleft}

\textbf{Properties of algebraic torsion, $torsion(C^*,d,d^*)$:}

\end{flushleft}

It will be helpful to clarify the relations between   these
approaches. This is accomplished by the following claims.
The mysterious sign $(-1)^{S(C^*)}$  in \ref{eqn.sign} is chosen so these properties
hold exactly.

\vspace{.2in}
\textbf{Claim A:}{ [\textit{Relation to Bilinear pairings}]}
 Let $(C^*,d)$ be a finite cochain
complex  let $\{c^{q,j}| j=1,\cdots, n_q\}$ be
the specified basis for $C^q$ and $ \{\ h^{q,j} | j=1 \cdots, v_q\}$ be
the specified basis for  $H^q(C^*,d )$. Thus $d$ is a differential
of degree $+1$.

Introduce non-degenerate complex bilinear pairings
$$
(.\ ,\ .) : \ C^q \times C^q  \rightarrow C
$$
defined by $(c^{q,i}, c^{q,j}) = \delta_{i,j}$. That is,
$( \Sigma_i \ \lambda_i \ c^{q,i}, \Sigma_j \ \mu_j \ c^{q,j}) = \Sigma_j \ \lambda_j \ \mu_j$.
Note $(.,.)$ is not Hermitian, rather $(\lambda a, b) = \lambda (a,b) = (a,\lambda b)$
for a complex number $\lambda$.

Let $d^* $ be unique differential of degree $-1$ with
 $(d a,b) = (a,d^*  b)$, i.e.,  the dual of $d$ under $(.,.)$.
Let ${\bf{h'}}_q  = \{ h'_{q,j} | j=1 \cdots, v_q \}$
be the dual basis in $H_q(C^*,d^*)$ to the chosen basis
 $ \{\ h^{q,j} | j=1 \cdots, v_q\}$ in $H^q(C^*,d)$ under this
pairing, $(.,.)$.
It is claimed:
$$
\begin{array}{l}
torsion(C^*,d,d^*,\{{\bf{h}}^q\}, \{ {\bf{h'}}_q\})
= \tau(C^*,d ,  \{{\bf{h}}^q\}, \{ {\bf{c}}_q\})^{2}  \label{eqn.chaim}
\end{array}
$$
The significance of this equation  is that the
square of the traditional torsion appears naturally in the bi-complex
setting.

\vspace{.2in}
\textbf{Claim B:}{[\textit{Relation to Eigenvalues}]
Suppose $(C^*,d,d^*)$ is
a finite bi-complex as above. Assume moreover
that the combinatorial Laplacian, $\Delta:= (d\ d^* + d^* \  d )
$, has no eigenvalue $0$. Record
the eigenvalues counted with multiplicities
of $\Delta_q = (\Delta|C^q) : C^q \rightarrow C^q$ as
$\{ \lambda_{q,i}| i=1, \cdots, n_q \} $.
By assumption they are all non-zero.

It is  claimed that the cohomologies   and homologies
$H_q(C^*,d)$ and $H^q(C^*,d^*)$ vanish,
so no choice of basis is needed. Moreover, the algebraic torsion
of $(C^*,d,d^*)$, a complex number,  is expressible as:
\begin{eqnarray}
torsion(C^* ,d , d^* ) \hspace{2in} \\ =
 [\prod_q \ (det(  \Delta_q))^{(-1)^q}]^{-1}
= [\prod_{q=0}^N \
 ( \prod_{j=1}^{n_q} \ \ \lambda_{q,j} )^{(-1)^{q} q}]^{-1} \nonumber \label{eqn.acyclic1}
\end{eqnarray}

This allows   the right   hand side to be extended
from the setting of non-zero eigenvalues to the
general case. This is precisely what is needed
in the situation of the zero-mode correction.

\vspace{.2in}
\textbf{Claim C:}{[ \textit{Stability Property} ] }
 Let $(C^*,d, d^*)$ be a finite chain complex
with differentials, $d, d^*$ of degrees $+1,-1$ respectively.

Suppose
specified bases
of the cohomologies, $H^q( C^*,d)$,
 say ${\bf{h}}^q = \{h^{q,j} | j=1 \cdots, v_q \}$
with $v_q = dim_C H^q(C^*, d)$,
and there are  also given  specified bases
of the homologies, $H_q( C^*,d^*)$, say ${\bf{h}}_q = \{h_{q,j} | j=1 \cdots, u_q \}$
with \newline $u_q = dim_C \ H_q(C^*, d^*)$.

Moreover, let $(D^*,d', d'^*)$ be another finite chain complex
with differentials $d',d'^*$ of degrees $+1,-1$ respectively
and suppose that the cohomology $H^*(D^*,d')$ and homology $H_*(D^*,d'^*)$
vanish. That is, $D^*$ is acyclic under both  $d'$ and $d'^*$.

In this circumstance, one has the natural isomorphisms,
$H^*(C^* \oplus D^*, (d \oplus d')) = H^*(C^*,d)$ and
$H^*(C^* \oplus D^*,  (d^* \oplus d'^*)) = H_*(C^*,d^*)$,
so the choices above may be used in computing the torsion
of the bi-complex, $(C^* \oplus D^*, d + d', d^* + d'^*)$.
It is claimed:
$$
\begin{array}{l}
 torsion(C^* \oplus D^*,(d \oplus d'),(d^* \oplus d'^*),\{{\bf{h}}^q\}, \{ {\bf{h}}_q\})
\hspace{1in} \nonumber \\
 = torsion(C^*,d,d^*,\{{\bf{h}}^q\}, \{ {\bf{h}}_q\}) \cdot
            torsion(D^*,d',d'^*) \nonumber  \\
\end{array}
$$
Here by acyclicity, $torsion(D^*,d',d'^*)$ is just a complex number.

\vspace{.3in}
This stability property is particularly desirable, and  depends on the choice of
signs in \ref{eqn.sign}. In the case that $\Delta|D_* $ has
all non-zero eigenvalues, the term $ torsion'(D^*,d,d^*)$
is expressible by claim B in terms of those eigenvalues.

In view of these formulas, one
deduces the equation \ref{eqn.stabilize1} of section 3. This
proves that the invariant defined in theorem \ref{thm.MAIN} is independent
of the choice of $K>0$.

\centerline{ \bf{Details of Proofs of Claims A,B,C: }}

\vspace{.1in}

{\bf{A:}} Here since $c^{q,i}, i=1 , \cdots , n_q$ is a basis for
$C^q$, and $B^{q+1} = d C^q$ has rank $s_{q+1}$, one  may reorder this
basis so the last $s_{q+1}$ elements, that is, \newline
$\{c^{q,n_q-s_{q+1}+1} , \cdots , c^{n_q}\}$  project by $d$ to a basis
of $B^{q+1}$. Set $b^{q+1,i} = d c^{q,i}$, for $n_q-s_{q+1}+1 \le i \le n_q$,
giving a  basis ${\bf{b}}^{q+1}$  for $B^{q+1}$, and set
 $\tilde{b}^{q+1,i} = c^{q,i}$. Let $\tilde{{\bf{b}}}^{q+1} = \{ \tilde{b}^{q+1,i} \}
=\{ c^{q,i}| n_q-s_{q+1}+1 \le i \le n_q \} \subset C^q$.
The span of these last elements maps bijectively to $B^{q+1} = dC^q$.

Now for each $i$, \ $1 \le i \le n_q-s_{q+1}$, there are unique numbers
$U_{q,i,j},\  n_q-s_{q+1}+1 \le j \le n_q$ with
$d  c^{q,i} = \Sigma_{j=n_q-s_{q+1}+1}^{n_q} \ U_{q,i,j} \ d  c^{q,j}$.
Thus, $z^{q,i}  = c^{q,i}   - \Sigma_{j=n_q-s_{q+1}+1}^{n_q} \ U_{q,i,j}  \ c^{q,j}$,
for $1 \le i \le n_q-s_{q+1}$,
 provides a basis for the cycles $Z^q = (ker \ d  : C^q \rightarrow C^{q+1} )$.
Since $Z^q $ projects onto $H^q(C^*,d )$ which is of dimension
$v_q$, we can and do reorder these $n_q-s_{q+1}$ entries
so that the $v_q$ elements $\{z^{q,i} | n_q-s_{q+1}-v_q+1 \le i \le n_q-s_{q+1}\}$
project to a basis of $H^q(C^*,d)$.

Since $C^q/Z^q \cong B^{q+1}$ under $d$, one  has $n_q - (s_q+v_q) = s_{q+1}$ or
equivalently $n_q = s_q + s_{q+1}+v_q$.

Set $lh^{q,k} =  z^{q, n_q-s_{q+1}-v_q+k}= z^{q,s_q +k}$ for $1 \le k \le v_q$,
and set ${\bf{lh}}^q = \{z^{q, s_q+k} | 1 \le k \le v_q \}
= \{ c^{q,s_q+ k} - \Sigma_{j=n_q-s_{q+1} +1}^{n_q} \ U_{q,s_q+ k,j}  \ c^{q,j}
| 1 \le k \le v_q\}$. These $q$-cocyles  project to a basis for $H^q(C^*,d)$.

With these choices one has  the ordered basis
 $$
\begin{array}{l}
 {\bf{b}}^q ,{\bf{lh}}^{q},  \tilde{ {\bf{b}} }^{q+1} \\
 = \{ d c^{q-1,i}|n_{q-1}-s_q+1 \le i \le n_{q-1}\} \\
\hspace{.5in} \sqcup \{  c^{q,s_q+k} - \Sigma_{j=n_q-s_{q+1}+1}^{n_q}
\ U_{q,s_q+k,j} \ c^{q,j}| 1 \le k \le v_q \} \\
\hspace{.5in} \sqcup \{  c^{q,j}| n_q-s_{q+1}+1 \le j \le n_q  \}
\end{array}
$$
of $C^q$ to compare to the basis  $\{c^{q,i}| 1 \le i \le n_q \}$.
Here $n_q-s_{q+1} = s_q + v_q$.

If $d c^{q-1,i} = \Sigma_{j=1}^{s_q} \ R_{q,i,j}\ c^{q,j}
+  \Sigma_{j=s_q+1}^{s_q+v_q} \ S_{q,i,j} \ c^{q,j}
+ \Sigma_{j=s_q+v_q+1}^{n_q} \ T_{q,i,j}\ c^{q,j}
$, for
$n_{q-1}-s_q+1 \le i \le n_{q-1}$.    Then the change of basis matrix from ${\bf{c}}^q$ to
${\bf{b}}^q ,{\bf{lh}}^{q},  \tilde{ {\bf{b}} }^{q+1}$ is
$$
M_q = \left(
 \begin{array}{ccc}
 R_q & 0 & 0 \\
S_q & I & 0 \\
T_q & -U_q & I \end{array} \right)
$$

In particular, $det(R_q ) \neq 0$ and for the above
explicit bases ${\bf{h}}^q = \{[z^{q,j}]|s_q+1 \le j \le s_q+v_q  \}$ for $H^q(C^*,d)$

$$
\begin{array}{l}
\tau(C^*,d, \{{\bf{h}}^q\}, \{ {\bf{c}}^q \} )
 := \prod_{q=0}^N \ [ {\bf{b}}^q, {\bf{lh}}^q,
\tilde{ {\bf{b}} }^{q+1} / {\bf{c}}^q  ]^{(-1)^{q+1}  } \\
= \prod_q \ (det(R_q))^{(-1)^{q+1}}
\end{array}
$$

Since reordering the basis elements of $\{ {\bf{c}}_q \}$ leaves the
square unchanged, one  may  use this specific choice of ordering to compute
the square
$\tau(C^*,d, \{{\bf{h}}^q\}, \{ {\bf{c}}^q \} )^2
= (\prod_q \ (det(R_q)^{(-1)^{q+1}})^2$.

\vspace{.2in}
Now consider elements $z^*_{q,i} \in C^q$  for $1 \le i \le v_q$
of the special form = $c^{q,s_q+i} + \Sigma_{j=1}^{s_q} \ V_{q,i,j} \ c^{q,j}$.
It is claimed  that  there is a unique choice of the $V_{q,i,j}$'s  which makes
these elements lie in  the q-cycles for $d^*$, the dual
of $d$ under the pairing $(.,.)$. Here
the q-cycles are $Z_q = (ker \ d^*|C^q : C^q \rightarrow C^{q-1})$.

The non-degenerate bilinear pairing $(.,.)$ on $C_q$ is
completely specified by $(c^{q,i},c^{q,j}) = \delta_{i,j}$,
the Kronecker delta function. This induces the identification
of $C^q$ with its dual, $(C^q)^* = Hom(C^q,C)$, given by the specified
choices of bases.

The constraint $z^*_{q,i} \in Z_q$ is the same as $(z^*_{q,i},B^q) =0$.
Since $\{ d\ c^{q-1,n_{q-1}- s_q+k}| \newline 1 \le k \le s_q \}$
is a basis for $B^q$, these constraints are the same as
the $s_q$ equalities:
$0  = (z^*_{q,i},d\ c^{q-1,n_{q-1}-s_q+k})
=  ( c^{q,s_q+i} + \Sigma_{j=1}^{s_q} \ V_{q,i,j} \ c^{q,j},
\Sigma_{m=1}^{s_q} \ R_{q,k,m}\ c^{q,m}
+  \Sigma_{m=s_q+1}^{s_q+v_q} \ S_{q,k,m} \ c^{q,m}
+ \Sigma_{m=s_q+v_q+1}^{n_q} \ T_{q,k,m}\ c^{q,m})
= (\Sigma_{j=1}^{s_q} \ V_{q,i,j} \ c^{q,j}, \newline \Sigma_{m=1}^{s_q} \ R_{q,k,m}\ c^{q,m})
+ (c^{q,s_q+i},  \Sigma_{m=s_q+1}^{s_q+v_q} \ S_{q,k,m} \ c^{q,m})
= \Sigma_{j=1}^{s_q} \  \ V_{q,i,j}  \ R_{q,k,j} + S_{q,k,s_q+i }
= (V \cdot R^t)_{i,k} + (S^t)_{s_q+i,k}$. Since $R$ is invertible,
the equations $(V \cdot R^t)_{i,k} = - (S^t)_{s_q+i,k}$ determine
$V$ uniquely.
The claim follows.

Here the rank of $d:C^q \rightarrow C^{q+1}$ is the same
as the rank of the adjoint via $(.,.)$ which is
$d^* :C^q \rightarrow C^{q+1}$. That is,
$s_{q+1} = dim \ B^{q+1} = dim  \ B_q$.

Let $lh'_{q,i} =c^{q,s_q+i} + \Sigma_{j=1}^{s_q} \ V_{q,i,j} \ c^{q,j}$
for $ 1\le i \le v_q$ for these unique choices. Then $lh'_{q,i} \in Z_q$
represent elements of $H_q(C^*,d^*)$
and moreover the pairing with the cycle representatives $lh^{q,j}$
specifying a basis for $H^q(C^*,d)$ is:
$(lh'_{q,i}, lh^{q,j}) =  (c^{q,s_q+i} + \Sigma_{j=1}^{s_q} \ V_{q,i,j} \ c^{q,j},
 c^{q,s_q+k}  - \Sigma_{j=n_q-s_{q+1}+1}^{n_q}
\ U_{q,s_q+k,j} \ c^{q,j} ) = \delta_{i,j}$.

Hence, the cycles $lk'_{q,i}$ represent in $H_q(C^*,d^*)$
the dual basis to the basis represented by
the cocycles $lh^{q,j}$ in $H^q(C^*,d)$.
Let ${\bf{lh}}'_q = \{ lh'_{q,i}| 1 \le i \le v_q\}$.

Since $R$ has non-zero determinant, and
$\{ d c^{q,n_q -s_{q+1} +k}|1 \le k \le s_{q+1}\}$
form a basis for $B^{q+1}$,  the
elements $\{ d^* c^{q+1,j}| 1 \le j \le s_{q+1}\}$
form a basis for $B_q = d^* C^{q+1}$.
So one may take ${\bf{b}'}_q = \{ d^* c^{q+1,j}|1 \le j \le s_{q}\}$
as a basis for $B_q$ and set
 $ \tilde{ {\bf{b}}' }_q
= \{ c^{q+1,j}|1 \le j \le s_{q+1}\}$.

Indeed, for $1 \le j \le s_{q+1}, 1 \le k \le s_{q+1}$, one computes
$( d^* c^{q+1,j}, c^{q,n_q -s_{q+1}+k}) \newline
= (  c^{q+1,j}, d c^{q,n_q -s_{q+1}+k})
= ( c^{q+1,j}, \Sigma_{m=1}^{s_{q+1}} \ R_{q+1,k,m} \ c^{q+1,m})
= R_{q+1,k,j}  = \newline( R_{q+1,k,j} \ c^{q,n_q -s_{q+1}+k},  c^{q,n_q -s_{q+1}+k})$,
so  $d^* c^{q+1,j} =
\Sigma_{m=1}^{s_{q+1}} \ R_{q+1,m,j}\ c^{q,s_q + v_q+ m}
+  \Sigma_{m=1}^{s_q} \ X_{q,k,m} \ c^{q,m}
+ \Sigma_{m=s_q+1}^{s_q+v_q} \ Y_{q,k,m}\ c^{q,m}
$
for some constants $X_{q,i,j},Y_{q,i,j}$.

Consequently, one gets an ordered  basis
$$
\begin{array}{l}
 {\bf{b}'}_q ,{\bf{lh}}'_{q},  \tilde{ {\bf{b}}' }_{q-1}  \\
 = \{ d^* c^{q+1,i}|1 \le i \le s_{q+1}\}
\sqcup \{ c^{q,s_q+i} + \Sigma_{j=1}^{s_q} \ V_{q,i,j} \ c^{q,j}| 1 \le i \le v_q\} \\
\hspace{.5in} \sqcup \{ c^{q,i}|1 \le i \le s_q\}\\
\end{array}
$$
for $C^q$. The change of basis matrix to  ${\bf{c}}^q$ is:
$$
N_q = \left(
 \begin{array}{ccc}
X & V & I \\
Y & I & 0 \\
(R_{q+1})^t & 0& 0 \end{array} \right)
$$

Hence, one  sees that for this ordering of the basis $\{c^{q,i}\}$,
$$
\begin{array}{l}
\tau'(C^*,d^*, \{{\bf{h}}_q\}, \{ {\bf{c}}_q \} )
=  \prod_q \ (det(N_q)^{(-1)^{q+1}} \\
=  \ \prod_q \ ( (-1)^{s_{q+1}s_q + s_q v_q + s_{q+1}v_q}  det(R_{q+1}))^{(-1)^{q+1}} \\
= (-1)^R \  \prod_q \ [(det(R_q))^{(-1)^{q+1}}]^{-1}
\end{array}
$$
with the sign  $(-1)^R = (-1)^{s_{q+1}s_q + s_q v_q + s_{q+1}v_q}$ coming
from permuting the blocks $I,I,(R_{q+1})^t$ which have ranks
$s_{q+1},s_q,v_q$ respectively.

Recall that $s_{q+1} = dim \ B^{q+1} = dim \ B_q$
and $v_q = dim \ H^q = dim \ H_q$.

But by  equation \ref{eqn.sign},   $S(C^*) = \Sigma_q \ S_q$ with
$S_q =  [dim\ B_{q-1}(C^*) \cdot dim_C \ B^{q+1}(C^* )
 \newline + dim_C \ B^{q+1}(C^*)  \cdot dim_C H_q(C^* ,d) +
dim_C \ B_{q-1}(C^*)  \cdot dim_C H^{q}(C^* ,d^*)]$.
So this special case, $S_q = s_q s_{q+1} + s_{q+1} v_q
+ s_q v_{q} $.
Consequently,
$(-1)^R = (-1)^{S(C^*)} $ with $S(C^*)$ as in \ref{eqn.sign}

Hence,
$$
\begin{array}{l}
torsion(C^*,d,d^*,\{ {\bf{h}}^q \},\{ {\bf{h}}'_q \}) \\
 = (-1)^{S(C^*)}
 \ \prod_{q=0}^N \ [ {\bf{b}}^q, {\bf{h}}^q,  \tilde{ {\bf{b}} }^{q+1}
 /  {\bf{b'}}_q, {\bf{h}}'_q, \tilde{ {\bf{b'}} }_{q-1} ]^{(-1)^{q+1}  } \\
= (-1)^{S(C^*)}  \  [\prod_q \ (det(R_q))^{(-1)^{q+1}}]/
 [(-1)^R  [\prod_q \ (det(R_q))^{(-1)^{q+1}}]^{-1}] \\
= \  [\prod_q \ (det(R_q))^{(-1)^{q+1}}]^{2} \\
=  \ [\tau(C^*,d,\{ {\bf{h}}_q \},\{ {\bf{c}}^q \} )]^2
 \end{array}
$$

This proves claim A for a special choice of bases
${\bf{h}}_q$. A change of basis  of $H_q(C^*,d)$ modifies both
sides by the same non-zero multiple so
it suffices to prove this result for any
specific choice. Changing the ordering of the $\bf{c_q}$
also leaves both sides unchanged, so claim A is established.

\vspace{.5in}
{\bf{B:}} If $\Delta$ has no vanishing eigenvalues, then
$\Delta$ is an isomorphism, so in particular, for each $x \in C^q$ there
is a $y \in C^q$ with $x = \Delta y = (d^* d y + d d^* y)$. Hence,
$C^q= d^* C^{q+1} + d C^{q-1}$. On the other hand, if $x \in d^* C^{q+1} \cap
d C^{q-1}$, then $d^*x = 0, d x = 0$ by $(d^*)^2 =0, (d )^2=0$, so
one gets
$\Delta  x= (d^*d + d  d^*) x =0$ also
which implies $x=0$. Hence, one has  a direct sum splitting
of $C^q$ as:
$$
C^q= d^* C^{q+1} \oplus d  C^{q-1}
$$

By $d^*\Delta = d^* d  d^* = \Delta d^*,
d  \Delta = d d^* d  = \Delta d $
the isomorphism $\Delta$ preserves the above splitting. Hence,
$\Delta$ maps $d^* C^{q+1}$  and  $d C^{q-1}$
 isomorphically to themselves. Since $\Delta|d^*C^* = d^* d $
and $\Delta| d  C^* = d  d^*$, these maps are isomorphisms
on $d^*C^{q+1}, d  C^{q-1}$ respectively. In particular, $d :
d^*C^{q+1}\rightarrow  d C^q$
is an injection and moreover $d^*: d  C^q\rightarrow d^*C^{q+1}$ is an injection
with composite an isomorphism. Hence, one   finds  all these maps
are necessarily isomorphisms. In summary, one has  induced isomorphisms:
$$
d^*C^{q+1} \stackrel{ d }{\rightarrow  } d  C^q\quad and
\quad d  C^q \stackrel{ d^*}{\rightarrow }
d^*C^{q+1} \ .
$$

In particular, $H_*(C^*,d) =0$ and $H^*(C^*,d^*) =0$,
as claimed. Also, by $H^*(C^*,d) = 0$, $n_q = s_q + s_{q+1}$.

Now take a basis, say $d y_{q,i}, 1 \le i \le s_q$ for $B^q  = d C^{q-1}$ for which
the matrix $\Delta$ has matrix representative
in Jordan block form. Then
$$
\Delta \ d y_{q,i}= d d^*   \ d y_{q,i}  = \lambda_{q,i}\
 d y_{q,i}  + \Sigma_{j<i} \ T_{q,ij} \
d y_{q,j} \quad  1 \le i \le s_q
$$
for some constants $T_{ij}$. Here the $\lambda_{q,i} \neq 0$ are the
generalized eigenvalues of $ \Delta|d C^{q-1}$ counted with their multiplicities.
Also $y_{q,i}  \in C^{q-1}$.

Since $\lambda_{q,i} \neq 0$ for any $i$, one  may invert the above to get
unique constants $U_{q,i,j}$ with
$$
 d d^* ( d  y_{q,i} + \Sigma_{j<i}\  U_{q,i,j} \
d y_{q,j} ) = \lambda_{q,i} \  d  y_{q,i}
$$

Let
$$
x_{q,i} = (1/\lambda_{q,i})(d  y_{q,i} + \Sigma_{j<i}\  U_{q,i,j} \ d y_{q,j})
 \quad so \quad d d^* \  x_{q,i} = d y_{q,i}
$$
with $x_{q,i} \in C^q$.

Since the $\{d y_{q,i}| 1 \le i \le s_q\}$ are a basis for $B^q = d C^{q-1}$
and $d d^* $ is an isomorphism of
$d  C^{q-1}$ onto itself, then necessarily $x_{q,i}$ is also a basis for $B^q = d C^{q-1}$.
By $d^* |d C^{q-1}$ an isomorphism
to $d^* C^q$, it follows that $\{d^*  x_{q,i}|1 \le i \le s_q\}$ is a basis for $d^* C^q$.
In particular, in view of $C^q=   d  C^{q-1} \oplus d^* C^{q+1}$
one obtains that
$$
{\bf{c}}_q =     \{ d  y_{q,j} | 1 \le j \le s_q\}
 \sqcup  \{ d^* \ x_{q+1,i}| 1 \le i \le s_{q+1} \}
$$
is a basis for $C^q$. In particular,
\begin{eqnarray}
dim \ B_q = dim \ B^{q+1} = s_{q+1} \label{eqn.help1}
\end{eqnarray}

The  computation of  the torsion $\tau ( C^*, d, \{   {\bf{c}}^q   \}  )$
proceeds as follows for these choices.

For this choice of ordered bases $ {\bf{c}}^q $, one  chooses ${\bf{b}}^q$, the ordered
 elements
$b^{q,i} = d y_{q,i}$, as basis for $B^q = d C^{q-1}$ inside $C^q$.
One specifies  ${\tilde{\bf{b}}^q}$ by $\tilde{b}_{q,i} = d^*  x_{q+1,i}$
since $d d^*  x_{q+1,i} = d y_{q+1,i}$. With these choices
$$
[ {\bf{b}}^q, \tilde{ {\bf{b}} }^{q+1} / {\bf{c}}^q  ]
= [\{ d  y_{q,j} \} \sqcup  \{ d^*  x_{q+1,i} \}
/ \{ d  y_{q,j} \} \sqcup  \{ d^*  x_{q+1,i} \}
 ] = +1
$$
Hence, $\tau ( C^*, d , \{   {\bf{c}}^q   \}  ) = +1$ for these choices.

\vspace{.2in}
In order to compute $\tau' ( C^*, d^* , \{   {\bf{c}}^q   \}  )$
one  may take ${\bf{b}'}_q$ to be the
 basis \newline $\{ d^*  x_{q+1,i}|1 \le i \le s_{q+1}\}$ of $B_q = d^* C^{q+1}$
and ${\tilde{\bf{b}'}_q}$ to be the elements $\{ x_{q,j}|1 \le j \le s_q \}$.
Hence,
$$
[ {\bf{b}'}_q, \tilde{ {\bf{b}'} }_{q-1} / {\bf{c}}_q  ]
 = [\{ d^*  x_{q+1,i}\} \sqcup  \{  x_{q,j}\}
  /  \{ d  y_{q,j} \} \sqcup  \{ d^*  x_{q+1,i} \}]
$$

In view of the equalities  $x_{q,i} = (1/\lambda_{q,i})(d y_{q,i}
 + \Sigma_{j<i}\  U_{q,ij} \ d y_{q,j})
$, this implies that
$$
[ {\bf{b}'}_q, \tilde{ {\bf{b}'} }_{q-1} / {\bf{c}}_q  ] =
 (-1)^{s_q s_{q+1}} (\prod_i \ (\lambda_{q,i})^{-1})
$$
since ${\flat}\{ d y_{q,j}\} =   s_q$ and
${\flat} \{d^* x_{q+1,i}\} = s_{q+1}$

In toto one gets:
\begin{eqnarray}
torsion(C^*,d^*,d ) \nonumber \\=  (-1)^{S(C^*)} \ \tau(C^*,d
,\{{\bf{c}}_q\})/\tau(C^*,d^* ,\{{\bf{c}}_q\})  \nonumber \\
 = (-1)^{S(C^*)} \ (-1)^R   \prod_q \ (\prod_i \lambda_{q,i})^{(-1)^{q+1}} \nonumber
\end{eqnarray}
with $R = \Sigma_q \ s_q s_{q+1}$.

Here by  equation \ref{eqn.sign},  $S(C^*) = \Sigma_q \ S_q$ with
$S_q =  [dim\ B_{q-1}(C^*) \cdot dim_C \ B^{q+1}(C^* )
 + dim_C \ B^{q+1}(C^*)  \cdot dim_C H_q(C^* ,d) +
dim_C \ B_{q-1}(C^*)  \cdot dim_C H^{q}(C^* ,d^*)]$.

In this doubly acyclic case $S_q = (dim \ B_{q-1})(dim B^{q+1})
= s_q s_{q+1}$ via the equality $s_q = (dim \ B_{q-1})$ observed
in equation \ref{eqn.help1}. Hence, $(-1)^R = (-1)^{S(C^*)}$
and so
\begin{eqnarray}
torsion(C^*,d^*,d )
=  \prod_q \ (\prod_i \lambda_{q,i})^{(-1)^{q+1}}
\end{eqnarray}

But the eigenvalues of $\Delta_q$ are those of $\Delta | d C^{q-1}$, i.e.,
$\{\lambda_{q,i}| 1\le i \le s_q \}$ and those of $\Delta |d^* C^{q+1}$. These last
under the isomorphism $d^* $ from  $d C^{q}$ to $d^*C^{q+1}$
 are the eigenvalues $\{\lambda_{q+1,j}|1 \le j \le s_{q+1}\}$.

That is, $det(\Delta_q) = (\prod_i \ \lambda_{q,i})(\prod_j \ \lambda_{q+1,j})$.
 In particular,  as in  Ray and Singer \cite{RaySinger1}
$$
\prod_q \ (det(\Delta_q))^{(-1)^q q} = \prod_{q>0} \ (\prod_i \lambda_{q,i})^{(-1)^q}.
$$

Substituting this into the above gives:
$$
torsion(C^*,d,d^* ) =  [\prod_q \ (det(\Delta_q))^{(-1)^q q}]^{-1}
$$
as claimed.

\vspace{.3in}
\textbf{C:}
Make  choices for $(C^*,d,d^*)$ and  make corresponding
choices for $(D^*,d',d'^*)$.
Comparing terms one finds the  same basis entries occurring
except for the ordering of the terms. For example, in $C^* \oplus D^*$
the determinants
$$
\begin{array}{l}
[ { \bf{b}}(C^*)^q \sqcup {\bf{ b}}(D^*)^q, {\bf{h}}(C^*)^q,
  \tilde{ {\bf{b}}}(C^*)^{q+1} \sqcup \tilde{ {\bf{b}}}(D^*)^{q+1}  \\
 \hspace{.5in} /  {\bf{b}}(C^*)_q \sqcup {\bf{b}}(D^*)_q, {\bf{h}}(C^*)_q,
\tilde{ {\bf{b}}}(C^*)_{q-1} \sqcup \tilde{ {\bf{b}}}(D^*)_{q-1} ]
\end{array}
$$
with the evident notation occurs, while the corresponding term
for $C^*$ is \newline $[ { \bf{b}}(C^*)^q , {\bf{h}}(C^*)^q,
  \tilde{ {\bf{b}}}(C^*)^{q+1}
  /  {\bf{b}}(C^*)_q , {\bf{h}}(C^*)_q,
\tilde{ {\bf{b}}}(C^*)_{q-1} ]$ and that for
$D_*$ is \newline $[ { \bf{b}}(D^*)^q ,
  \tilde{ {\bf{b}}}(D^*)^{q+1}
  /  {\bf{b}}(D^*)_q ,
\tilde{ {\bf{b}}}(D^*)_{q-1} ]$.  The first of these
is the product of the next two times a sign
$(-1)^{T_q}$ which accounts for the change of ordering.

Let $r_q = dim_C\  B_q(C^*) = dim \ im \ d^*_{q+1}[C^{q+1}]$, \
$s_q = dim_C \ B^q(C^*) = dim \ im \ d_{q-1} \ [C^{q-1}]$,\
$u_q = dim_C \ H_q(C^*,d^*), v_q = dim_C \ H^q(C^*,d)$.
Let $x_q = dim_C\  B_q(D^*) = dim \ im \ d^*_{q+1}[D^{q+1}]$,\
$y_q = dim_C\  B^q(D^*) = dim \ im \ d_{q-1}[D^{q-1}]$.
With the above notation, the number of elements in ${\bf{ b}}(D^*)_q$
is $x_q$ and each  must be passed by $ {\bf{h}}(C^*)_q,\tilde{ {\bf{b}}}(C^*)_{q-1}$
which have $u_q + r_{q-1}$ elements to get the desired order.
Similarly, the number of elements in ${\bf{b}}(D_*)^q$ is $y_q$
and each must be passed by $\bf{h}(C^*)^q, \tilde{ {\bf{b}}}(C^*)^{q+1}$
which has $v_q + s_{q+1}$  to get the correct order.
Consequently, multiplying over $q$ gives the relation
$$
\begin{array}{l}
torsion'(C^*\oplus D^*,(d+ d'),(d^*+ d'^*),\{ {\bf{h}}^q \},\{ {\bf{h}}_q \}) \\
= (-1)^{S(C^* \oplus D^*)} (-1)^{S(C^*)}(-1)^{S(D^*)} (-1)^T \\
\hspace{.2in} \cdot   torsion'(C^*,d,d^*,
\{ \textbf{h}_q \},
\{ \textbf{h}_q \}) \cdot
torsion'(D_*,d,d^*)
\end{array}
$$
where by the above $T =\Sigma_q \ T_q$ with
$T_q =  x_q(u_q + r_{q-1}) + y_q(v_q + s_{q+1})$.

Now by $B_q \cong C_{q+1}/Z_{q+1}, H_q \cong Z_q/B_q$ one has for $n_q =
dim_C \ C^q$, the equations $r_0 = n_0-u_0, r_1=n_1-u_1-r_0,r_2 = n_2-u_2-r_1, \cdots$,
which inductively gives
\begin{eqnarray}
r_q = dim_C \ B^q(C^*)
= (\Sigma_{i=0}^q \ (-1)^i \ n_{q-i}) - (\Sigma_{i=0}^q \ (-1)^i \ u_{q-i}).
\label{eqn.rs1}
\end{eqnarray}
   Similarly from
$B^q \cong C_{q-1}/Z^{q-1}, H^q \cong Z^q/B^q$, one has the equations
$s_1 = n_0 -v_0, s_2 = n_1-v_1-s_1, \cdots$ which  inductively give
\begin{eqnarray}
s_{q+1} = dim_C \ B^{q+1}(C^*) = (\Sigma_{i=0}^q  \ (-1)^i \ n_{q-i}) -
(\Sigma_{i=0}^q (-1)^i \ v_{q-i}). \hspace{.4in} \label{eqn.rs2}
\end{eqnarray}
In particular, one has
\begin{eqnarray}s_{q+1} = r_q + \Sigma_{i=0}^q \ (-1)^i \
[u_{q-i}- v_{q-i}].  \label{eqn.rs3}
\end{eqnarray}

Applied to the doubly acyclic complex, $(D^*,d ,d^*)$, with
$\hat{n}_q = dim_C \ D^q$, this yields the formulas:
\begin{equation} \label{eqn.rs4}
y_{q+1} = x_{q} = \Sigma_{i=1}^q \ (-1)^i \ \hat{n}_{q-i}
\end{equation}

As an example of the usefulness of these relations, consider the
difference
$$
\begin{array}{l}
D_1 := dim\ B_{q-1}(C^*\oplus D^*) \cdot dim_C \ B^{q+1}(C^* \otimes D^*)\\
-dim\ B_{q-1}( C^*) \cdot dim_C \ B^{q+1}( C^*)
-dim\ B_{q-1}( D^*) \cdot dim_C \ B^{q+1}( D^*)\\
\end{array}
$$ which in the above notation equals
$r_{q-1} y_{q+1} + x_{q-1} s_{q+1}$. By equation \ref{eqn.rs4} the last
is precisely $r_{q-1} x_q + y_{q} s_{q+1} $ which is two of the terms in
$T_q$.

The difference
$$
\begin{array}{l}
D_2 := dim_C \ B^{q+1}(C^*\oplus D^*) \cdot dim_C H_q(C^* \oplus D^*) \\
- B^{q+1}(C^*)  \cdot dim_C H_q(C^*,d)
- B^{q+1}(  D^*)  \cdot dim_C H_q( D^*)
\end{array}
$$
via $H_q(D^*, d' ) =0$ is in the above notation
$B^{q+1}( D^*)  \cdot u_q  = y_{q+1} u_q = x_q u_q$ which is another of the
terms in $T_q$.

Finally, the difference
$$
\begin{array}{l}
D_3 := B_{q-1}(C^* \oplus D^*) \cdot dim_C H^{q}(C^* \oplus D^*) \\
-B_{q-1}(C^* ) \cdot dim_C H^{q}(C^*,d^*)
-B_{q-1}(D^*) \cdot dim_C H^{q}( D^*)$ via $H^{q}(D^*,d'^*) =0
\end{array}
$$
is in the above notation $B_{q-1}(D^*) \cdot v_q= x_{q-1} v_q = y_q v_q$ which is
the last term in $T_q$.

These last three formulas  show that the sign \newline
$(-1)^{S(C^* \oplus D^*)} (-1)^{S(C^*)}(-1)^{S(D^*)} (-1)^T $ is always
$+1$ proving claim C since $S(C^*)$ is defined precisely by
$S(C^*) =  \Sigma_q \ [dim\ B_{q-1}(C^*) \cdot dim_C \ B^{q+1}(C^* ) \newline
 + dim_C \ B^{q+1}(C^*)  \cdot dim_C H_q(C^* ,d^*)  +
dim_C \ B_{q-1}(C^*)  \cdot dim_C H^{q}(C^*,d )]$.

\section{Variation of algebraic torsion
with changing Hermitian inner product on $TW$:} \label{section5}

Let $g(u), a \le u \le b$, be a smoothly varying Hermitian  inner product on
the bundle $TW$. Suppose that
$K>0$ is not a generalized eigenvalue for any
of the operators $\Delta_{E,\bar{\partial},g(u)}$ associated
to  these Hermitian metrics, $g(u)$ acting on $A^{p,q}(W,E)$.

The spectral projections $\Pi_{E,K,g(u)}$ onto the finite dimensional
subspaces $C^{p,*}(t) = \bigoplus_{q=0}^N \ C^{p,q}(E,K,g(u))$, the
spans of the generalized eigenvalues with real part
less than $K$,
by standard spectral theory, are smoothly varying operators in  $u$, for
$a \le u \le b$.

Hence, for any fixed $u_0, a \le u_0 \le b$,
 the fixed spectral projection $\Pi_{E,K,g(u_0)}$ applied
to $C^{p,q}(E,K,g(u))$ must be an isomorphism onto
$C^{p,q}(E,K,g(u_0))$ for $u$ sufficiently close to $u_0$.
 Pick
an $\epsilon>0$ so this is true for all $p,q$ for all $|u - u_0| \le \epsilon$
with
$a \le u \le b$.

Let $P_{p,q}(u_0)$ be the kernel of the reference spectral projection
$\Pi_{p,q,E,K}(u_0)$. Then one has  a direct sum decomposition of $(L^{p,q})^2(W,E)$:
$$
(L^{p,q})^2 (W,E) =C^{p,q}(E,K,g(t_0)) \  \oplus \ P_{p,q}(t_0)
$$
with the first summand of finite dimension and the second closed
in $(L^{p,q})^2(W,E)$.

Recall that for two ordered  bases, say
$Y = \{y_1,\cdots,y_m\}$ and $X = \{x_1,\cdots,x_m\}$
of an $m$ dimensional vector space $V$, the determinant of change of basis
matrix  $[Y/X]$ is the complex number with
$$
y_1 \wedge \cdots \wedge y_m = [Y/X] \ x_1 \wedge \cdots  \wedge x_m
$$

If $V$ is a subspace  $j : V \subset W$ of a complex vector space, $W$,
which is a direct sum  $W = U_1 \oplus U_2$ and the projection,
$\Pi: V \rightarrow U_1$,
 of $W $  onto $U_1$
along  $U_2$ is defines an  isomorphism $V \cong U_1$,
 then one  may compute the number $[Y/X]$ as follows:

Write out the elements, $x_j,y_j$ in this decomposition,
$x_j = \Pi(x_j) + (Id-\Pi)x_j, y_k = \Pi(y_j) + (Id-\Pi)y_k$ and expand
the terms of the wedge product with respect to the direct sum
decomposition
$$
\wedge^m \ W \cong  \oplus_{k=0}^m \ \wedge^{k} \ U_1 \otimes \wedge^{m-k} \ U_2
$$
where $m = dim_C \ U_1$.

Comparing the entries in the summand $\ \bigwedge^m \ U_1$ gives  the formula:
$$
\Pi(y_1) \wedge \cdots \wedge \Pi(y_m) = [Y/X] \ \Pi(x_1) \wedge \cdots  \wedge \Pi(x_m)
$$
this computation being carried out in $U_1$. This method will be
utilized for the above fixed decomposition at $u=u_0$.

\vspace{.3in}

Now to compute the change of algebraic torsion for varying $u$,
for the family of complexes $(C^{p,*}(u),d,d^*) = (\oplus_q \ C^{p,q}(E,K,g(u)), d,d^*)$,
one  may first specify choices for $u=u_0$ and then
for $u$ nearby.

 Choose
an ordered  basis, say ${\bf{h}}^q = \{ h^{q,j}| j =1 \cdots v_q\}$
for $H_{\bar{\partial}_E}^{p,q}(W,E)$. As noted  above,
$$
\begin{array}{l}
H^q(C^{p,*}(u_0),d) = H^q(\oplus_{q'} \  C^{p,q'}(E,K,g(t_0),\bar{\partial}_E)  \\
 \cong H^q(A^{p,*}(W,E),\bar{\partial}_E) \cong H_{\bar{\partial}_E}^{p,q}(W,E)
\end{array}
$$
so one   may chose $lh^{q,j}(u_0) \in \ ker \ \bar{\partial}_E:
C^{p,q}(E,K,g(u_0) ) \rightarrow C^{p,q+1}(E,K,g(u_0) )$ which project to
the ordered basis $\{h^{q,j}| 1 \le j \le v_q\}$ chosen.

Let ${\bf{lh}}^q(u_0)$
be the ordered set $\{ lh^{q,j}(u_0)|j=1,\cdots v_q\}$.

Since $\hat{\star}_{g(u)} \otimes Id_E$ commutes with $\Delta_{E,\bar{\partial},g(u)}$,
it maps $C^{p,q}(E,K,g(u))$ isomorphically to $C^{N-q,N-p}(E,K,g(u))$
and since $\hat{\star}_{g(u)} \ \hat{\star}_{g(u)} = \pm Id$, it interpolates the coboundary map
$D''$
up to sign with the differential $d^* = - (\hat{\star}_{g(u)} \otimes Id_E)
D''(\hat{\star}_{g(u)} \otimes  Id_E)$. Hence, it suffices to make
choices for the complex \newline $(C^{N-*,N-p}(E,K,g(u)), D'')$ and transplant them via
$\hat{\star}_{g(u)}$ to the complex \newline $(C^{p,\star}(E,K,g(u)), d^*)$.

In a parallel fashion,
choose
an ordered  basis, say ${\bf{h'}}^q = \{ {h'}^{q,j}| j =1 \cdots u_q\}$
for $H_{D''}^{q,N-p}(W,E)$. Again,
$
H^q(\oplus_{q'} C^{q',N-p}(E,K,g(u_0),D'')
  \cong H_{D''}^{q,N-p}(W,E)
$
so one   may chose ${{lh'}^{q,j}}(u_0) \in \ ker \ D'':
C^{q,N-p}(E,K,g(u_0) ) \rightarrow \newline C^{q+1,N-p}(E,K,g(u_0) )$ which project to
the ordered basis $\{{h'}^{q,j}| 1 \le j \le u_q\}$ chosen.
Let ${\bf{lh'}}^q(u_0)$
be the ordered set $\{ {lh'}^{q,j}(u_0)|j=1,\cdots u_q\}$.

Next pick  for each $q$ an an ordered basis, say
${\bf{b}}^q(u_0) \newline
= \{ b^{q,1}(u_0), \cdots, b^{q,s_q}(u_0)\}$ for the coboundaries,
$B^q(u_0) = image \ \bar{\partial}_E: \newline C^{p,q-1}(E,K,g(u_0))
\rightarrow C^{p,q}(E,K,g(u_0)$. Finally pick
for each $b^{q,j}(u_0)$ an element $\tilde{b}^{q,j}(u_0) \in C^{p,q-1}(E,K,g(u_0))$
with $\bar{\partial}_E \ \tilde{b}^{q,j}(u_0) = b^{q,j}(u_0)$.
 Let $\tilde{{\bf{b}}}^q(u_0) =
\{ \tilde{b}^{q,1}(u_0), \cdots , \tilde{b}^{q,s_q}(t_0)\}$ be ordered
by increasing last index. Here $\tilde{{\bf{b}}}^q(u_0) \subseteq
C^{p,q-1}(E,K,g(u_0))$.

 In a parallel fashion,
 pick for each $q$ an  ordered basis, say ${\bf{b'}}^q(u_0)
= \newline \{ {b'}^{q,j}(u_0), \cdots, {b'}^{q,s_q}(u_0)\}$ for the coboundaries for $D''$,
$ image \ D'': |\newline C^{q-1,N-p}(E,K,g(u_0))
\rightarrow   C^{q,N-p}(E,K,g(u_0)$. Finally pick
for each ${b'}^{q,j}(u_0)$ an element $\tilde{b'}^{q,j}(u_0) \in C^{q-1,N-p}(E,K,g(u_0))$
with $D'' \ \tilde{b'}^{q,j}(u_0) = {b'}^{q,j}(u_0)$.
 Let $\tilde{{\bf{b'}}}^q(u_0) =
\{ \tilde{b'}^{q,1}(u_0), \cdots , \tilde{b'}^{q,s_q}(u_0)\}$ be ordered
by increasing last index. Here
$\{\tilde{{\bf{b'}}}^q(u_0) \} \subseteq  C^{q,N-p}(E,K,g(u_0))$.

Since $\Delta_{E,\bar{\partial},g(u_0)}$ commutes
with the differential $\bar{\partial}_E$ the spectral
projection $\Pi_{E,K,g(u_0)}$ commutes with $\bar{\partial}_E$
and hence induces a chain mapping
from $(C^*(u),\bar{\partial}_E) $ to $(C^*(u_0),\bar{\partial}_E)$.
In particular, for $|u - u_0|  \le \epsilon$ and $u \in [a,b]$,
the induced map is an isomorphism of chain complexes.
In particular, the coboundaries $B^q(t) $ map isomorphically
to the coboundaries $B^q(u_0)$, the cocyles to cocycles,
and cohomologies to cohomologies.

For  $u$ with $|u-u_0| \le \epsilon$,  let $lh^{q,j}(u),  b^{q,k}(u),\tilde{b}^{q,l}(u)$
be the unique elements of $C^{p,q}(E,K,g(u))$ mapping under the fixed
spectral projection $\Pi_{E,K,g(u_0)}$ to the elements
$lh^{q,j}(u_0),  b^{q,k}(u_0),\tilde{b}^{q,l}(u_0)$. Then necessarily,
$\bar{\partial}_E \tilde{b}^{q,l}(u) = b^{q,l}(u) $,
the $\{b^{q,k}(u)\}$'s form an ordered basis for the coboundaries
$B^q_p(u)$, $\bar{\partial}_E lh^{q,j}(u) =0$
and $\{lh^{q,j}(u)\}$ project to a basis for $H^q(C^*_p(u),\bar{\partial}_E)$.
Since the isomorphism $H^q(C^*_p(u),\bar{\partial}_E) = H^q(A^{p,*}(W,E),\bar{\partial}_E)$
 carries
$ lh^{q,j}(u)$ to the cohomologous \newline
$lh^{q,j}(u_0)\} $ which represents $h^{q,j}
\in H^{p,q}(W,E)$, the element $lh^{q,j}(u)$ under the isomorphism
$H^q(C^*_p(u),d) = H^q(A^{p,*}(W,E),d) = H^{p,q}_{\bar{\partial}}(W,E)$ represents the same
element also.

In particular, the ordered set of cocyles, ${\bf{lh}}^q(u)
= \{ lh_{q,j}(u)|j=1,\cdots v_q\}$ represent cohomology classes
which map to the above basis, ${\bf{h}}^q = \{ h^{q,j}| j =1 \cdots v_q\}$
of $H^{p,q}(W,E)$. Also the ordered coboundaries
${\bf{b}}^q(u)
= \{ b^{q,1}(u), \cdots, b^{q,s_q}(u)\}$ give an ordered
basis for $B^q(u)$ and the ordered set
$\tilde{{\bf{b}}}^q(u) =
\{ \tilde{b}^{q,1}(u), \cdots , \tilde{b}^{q,s_q}(u)\} \newline \subset C^{p,q-1}(u)$
maps under $d = \bar{\partial}_E$ to the ordered
basis ${\bf{b}}^q(u)$ above.

Correspondingly for $C^{*,N-p}(K,g(u))$ and the coboundary $D''$,
for  $u$ with $|u-u_0| \le \epsilon$,  let ${lh'}^{q,j}(u),  {b'}^{q,k}(u),\tilde{b'}^{q,l}(u)$
be the unique elements of $C^{q,N-p}(E,K,g(u))$ mapping under the fixed
spectral projection $\Pi_{E,K,g(u_0)}$ to the elements \newline
${lh'}^{q,j}(u_0),  {b'}^{q,k}(u_0),\tilde{b'}^{q,l}(u_0)$. Then necessarily,
$D'' \tilde{b'}^{q,l}(u) = {b'}^{q,l}(u) $,
the $\{{b'}^{q,k}(u)\}$'s form an ordered basis for the coboundaries
of $D''$, $D'' {lh'}^{q,j}(u) =0$
and $\{{lh'}^{q,j}(u)\}$ project to a basis for $H^q(C^{*,N-p}(K,g(u)),D'')$.
Since the isomorphism \newline $H^q(C^{*,N-p}(K,g(u)),D'') = H^q(A^{*,N-p}(W,E),D'')$
 carries
$ {lh'}^{q,j}(u)$ to the cohomologous ${lh'}^{q,j}(u_0)\} $ which represents ${h'}^{q,j}
\in H_{D''}^{q,N-p}(W,E)$, the element ${lh'}^{q,j}(u)$ under the isomorphism
$H^q(C^{*,N-p}(K,g(u)),D'') =   \newline
H^q(A^{*,N-p}(W,E),D'') = H_{D''}^{q,N-p}(W,E)$ represents the same
element also.

In particular, the ordered set of cocyles, ${\bf{lh'}}^q(u)
= \{ {lh'}_{q,j}(u)|j=1,\cdots u_q\}$ represent cohomology classses
which map to the above basis, ${\bf{h'}}^q = \newline \{ {h'}^{q,j}| j =1 \cdots u_q\}$
of $H_{D''}^{q,N-p}(W,E)$. Also, the ordered coboundaries
${\bf{b'}}^q(u)
= \{ {b'}^{q,1}(u), \cdots, {b'}^{q,r_q}(u)\}$ give an ordered
basis for the coboundaries  and the ordered set
$\tilde{{\bf{b'}}}^q(u) =
\{ \tilde{b'}^{q,1}(u), \cdots , \tilde{b'}^{q,r_q}(u)\} \subseteq C^{q-1,p}(u)$
maps under $D''$ to the ordered
basis ${\bf{b'}}^q(u)$ above.

In particular,
the ordered set $ (\hat{\star}_{g(u)}\otimes Id_E) {\bf{lh'}}^{N-q}(u) \subset C^{N-q,N-p}(E,K,g(u))$
is  in the kernel of $d^*$ and projects to a basis of $H^{N-q}(C^*(u),d^*)$.
Similarly, the ordered set $(\hat{\star}_{g(u)}\otimes Id_E){\bf{b'}}^{N-q}(u)$
forms an ordered basis for these boundaries. Moreover,
$d^*(\hat{\star}_{g(u)}\otimes Id_E)^{-1}\tilde {\bf{b'}}^{N-q}(u)) =
\pm (\hat{\star}_{g(u)}\otimes Id_E){\bf{b'}}^{N-q}(u)$, so these form a suitable
choice to compute the algebraic torsion for the complex $C^{p,\star}(K,g(u)), d^*)$. One sets
 $\{{\bf{h'}}_q(u) \}  = \{[(\hat{\star}_{g(u)}\otimes Id_E) {\bf{lh'}}^{N-q}(u)]\}$
in $H^q(C^*(u),d^*)$. Under the isomorphism
$$
H^q(C^*(u),d^*) \cong H_{D''}^{N-q,N-p}(W,E)
$$
described above,  this basis corresponds to the chosen basis
${\bf{h'}}^{N-q}$ of  \newline $H_{D''}^{N-p,N-q}(W,E)$, which is
independent of $u$.

With these choices one is able to write the algebraic torsion as a function
of $u$ for $|u-u_0| \le \epsilon$. With the short hand $\hat{\star}_{g(u),E}
= (\hat{\star}_{g(u)} \otimes Id_E)$ it is:
$$
\begin{array}{l}
torsion'(C^*_p(t),d,d^*, \{{\bf{h}}^q(u) \},\{{\bf{h'}}_q(u) \} )  \\
= (-1)^S  \ \prod_{q=0}^N \  [ \{{\bf{b}}^q(u)\}, \{{\bf{lh}}^q(u)\}
,\{\tilde{{\bf{b}}}^{q+1}(u)\} /  \\
 \hat{\star}_{g(u),E} \{{\bf{b'}}_{N-q}(u)\}, \hat{\star}_{g(u),E}\{{\bf{lh'}}^{N-q}(u)\}
,\hat{\star}_{g(u),E}\{\tilde{{\bf{b'}}}^{N-(q+1)}(u)\}]^{(-1)^{q+1}}
\end{array}
$$
for a fixed sign $(-1)^S$.

This may be computed via the above device.

For example, if the ordered basis $  \{{\bf{b}}^q(u)\}, \{{\bf{lh}}^q(u)\}, \cdots ,
\{\tilde{{\bf{b}}}^{q+1}(u)\}$ is $[e_1,\cdots,e_L]$ and the ordered
basis
$\hat{\star}_{g(u),E} \{{\bf{b'}}_{N-q}(u)\}, \hat{\star}_{g(u),E}\{{\bf{lh'}}^{N-q}(u)\}
,
\hat{\star}_{g(u),E}\{\tilde{{\bf{b'}}}^{N-(q+1)}(u)\}$ is $[f_1,\cdots, f_L]$, then the required
determinant is $a$ where
$(\wedge e_j) = a ( \wedge \ f_k) $. The required number
$a$ can be computed by just taking the projection onto the fixed
summand $\wedge^L \ C^{p,q}(E,K,g(u_0))$ by means of the spectral
projection $\Pi_{E,K,g(u_0)}$.

By the above  choices one has
$$
\begin{array}{l}
b^{q,j}(u) = b^{q,j}(u_0) + v \ for \ some  \ v \in P_q(u_0)\\
\tilde{b}^{q,k}(u) = \tilde{b}^{q,j}(u_0) + v  \ for \ some  \ v \in P_q(u_0)\\
lk^{q,l}(u) = lk^{q,l}(u_0) +  v  for \ some  \ v \in P_q(u_0)\\
\end{array}
$$

Thus, the projection of  the wedge product of the terms in the ordered
basis $\{{\bf{b}}^q(u)\}, \{{\bf{lh}}^q(u)\}
,\{\tilde{{\bf{b}}}^{q+1}(u)\}$ is precisely the projection wedge product
of the entries in the reference  ordered basis $\{{\bf{b}}^q(u_0)\}, \{{\bf{lh}}^q(u_0)\},
,\{\tilde{{\bf{b}}}^{q+1}(u_0)\}$. In particular, the
derivative of this projection  by $u$  is zero.

By the above  choices,
$$
\begin{array}{l}
{b'}^{q,j}(u) = {b'}^{q,j}(u_0) + v \ for \ some  \ v \in P_{N-p}(u_0)\\
\tilde{b'}^{q,k}(u) = \tilde{b'}^{q,j}(u_0) + v  \ for \ some  \ v \in P_{N-p}q(u_0)\\
{lk'}^{q,l}(u) = {lk'}^{q,l}(u_0) +  v  for \ some  \ v \in  P_{N-p}(u_0)\\
\end{array}
$$

Hence,  since $\hat{\star}_{g(u_0),E}$ commutes with
$\Delta_{E,\bar{\partial},g(u_0)}$  and also the fixed
projection $\Pi_{E,K,g(u_0)}$,
$$
\begin{array}{l}
\hat{\star}_{g(u_0),E} {b'}^{q,j}(u) = \hat{\star}_{g(u_0),E}{b'}^{q,j}(u_0) + w \
for \ some  \ w \in P_q(u_0)\\
\hat{\star}_{g(u_0),E}\tilde{b'}^{q,k}(u) =\hat{\star}_{g(u_0),E}\tilde{b'}^{q,j}(u_0) + w
\ for \ some  \ w \in P_q(u_0)\\
\hat{\star}_{g(u_0),E}{lk'}^{q,l}(u) = \hat{\star}_{g(u_0),E}{lk'}^{q,l}(u_0) +  w  for
 \ some  \ w \in  P_q(u_0)\\
\end{array}
$$

So in particular,
$$
\begin{array}{l}
d/du \ \Pi_{E,K,g(u_0)} \hat{\star}_{g(u),E} {b'}^{q,j}(u)|_{u=u_0}
 = d/du \ \Pi_{E,K,g(u_0)} \hat{\star}_{g(u),E} {b'}^{q,j}(u_0)|_{u=u_0} \\
d/du \ \Pi_{E,K,g(u_0)} \hat{\star}_{g(u),E} \tilde{b'}^{q,k}(u) |_{u=u_0}
 = d/du \ \Pi_{E,K,g(u_0)} \hat{\star}_{g(u),E} \tilde{b'}^{q,k}(u_0) |_{u=u_0} \\
d/du \ \Pi_{E,K,g(u_0)} \hat{\star}_{g(u),E}{lk'}^{q,l}(u) |_{u=u_0}
 = d/du \ \Pi_{E,K,g(u_0)} \hat{\star}_{g(u),E}{lk'}^{q,l}(u_0) |_{u=u_0} \\
\end{array}
$$

Consequently, by these observations, one has the simplification
for $|u - u_0| \le \epsilon$,
$$
\begin{array}{l}
d/du \ torsion'(C^*_p(u),d,d^*, \{{\bf{h}}^q(u) \},\{{\bf{h'}}_q(u) \} )|_{u=u_0}  \\
= (-1)^S \ d/du \ \prod_{q=0}^N \  [ \{{\bf{b}}^q(u_0)\}, \{{\bf{lh}}^q(u_0)\}
,\{\tilde{{\bf{b}}}^{q+1}(u_0)\} /  \\
\hspace{.2in}\Pi_{E,K,g(u_0)} \hat{\star}_{g(u),E} \{{\bf{b'}}_{N-q}(u_0)\}, \\
\hspace{.5in} \Pi_{E,K,g(u_0)} \hat{\star}_{g(u),E}\{{\bf{lh'}}^{N-q}(u_0)\} \\
\hspace{.7in} , \Pi_{E,K,g(u_0)}\hat{\star}_{g(u),E}\{\tilde{{\bf{b'}}}^{N-(q+1)}(u_0)\}]^{(-1)^{q+1}} \ |_{u=u_0}
\end{array}
$$
since the above method of computation yields the
same result for both expressions.

This last gives
$$
\begin{array}{l}
d/du \ \left[ \frac{torsion'(C^{p,*}(u),d,d^*, \{{\bf{h}}^q(u) \},\{{\bf{h'}}_q(u) \} )}
 {torsion'(C^{p,*}(u_0),d,d^*, \{{\bf{h}}^q(u_0) \},\{{\bf{h'}}_q(u_0) \} )} \right]^{-1} \ |_{u=u_0}\\
= d/du \  \prod_{q=0}^N \
[ \Pi_{E,K,g(u_0)}\hat{\star}_{g(u),E} \{{\bf{b'}}_{N-q}(u_0)\}, \\
\hspace{.5in} \Pi_{E,K,g(u_0)}\hat{\star}_{g(u),E}\{{\bf{lh'}}^{N-q}(u_0)\}
,\Pi_{E,K,g(u_0)}\hat{\star}_{g(u),E}\{\tilde{{\bf{b'}}}^{N-(q+1)}(u_0)\} / \\
\hspace{.3in}   \hat{\star}_{g(u_0),E} \{{\bf{b'}}_{N-q}(u_0)\},
           \hat{\star}_{g(u_0),E}\{{\bf{lh'}}^{N-q}(u_0)\} \\
\hspace{.7in} ,    \hat{\star}_{g(u_0),E}\{\tilde{{\bf{b'}}}^{N-(q+1)}(u_0)\}]^{(-1)^{q+1}}  \ |_{u=u_0}
\end{array}
$$

Consequently, one gets
\begin{lemma} \label{lemma.5vary}
Under the above hypotheses:
$$
\begin{array}{l}
d/du \ \log \ {torsion(C^*_p(u),d,d^*, \{{\bf{h}}^q(u) \},\{{\bf{h'}}_q(u) \} )}|_{u=u_0}\\
= -  \ \Sigma_{q=0}^n (-1)^q  Tr ( \Pi_{E,K,g(u_0)} ( \alpha_{g(u_0)} \times Id_E) |V(p,q,E,K,g(u_0)))
\end{array}
$$
with $\alpha_{g(u)}$ the bundle mapping $ \alpha_{g(u)}= (\hat{\star}_{g(u)})^{-1}\  d/du \ (\hat{\star}_{g(u)})
= \newline  - d/du \ (\hat{\star}_{g(u)}) \  (\hat{\star}_{g(u)})^{-1} $
\end{lemma}

\section{
Definition of the flat  complex analytic torsion $\tau(M,F)$}  \label{section.analyticRiem1}

Let $M$ be a closed, Riemannian, smooth manifold of dimension $n$. Let
the metric on $TM$, the tangent bundle of $M$ be denoted by $g: TM \times TM \rightarrow R$.
For convenience we assume $M$ is oriented although this is not necessary.

Suppose that $F \rightarrow M$ is a flat $k$ dimensional complex vector bundle. That is,
the transition matrices  of $F$ are locally constant ( but not necessarily unitary).

The considerations above
carry over almost verbatim to this smooth setting
when the d-bar Laplacian above $\Delta^{\flat}_{E,\bar{\partial}}$ is replaced
by the flat Laplacian
$\Delta^{\flat}_{F}$. The result is the metric independence of the analytic torsion
$\tau(M,E)$ of the metric for flat bundle $E$ over an odd dimensional Riemannian
manifold $M$.
In this section we merely state the parallel results. The proofs are
quite the same in detail.

Utilizing the method of section \ref{section.FlatC}, the exterior derivative
$d$ and its adjoint $d^*$ and the standard Laplacian
$\Delta = (d + d^*)^2 = d d^* + d^* d$ acting on complex valued forms $A^*(M,C)$ have canonical
extensions, called here $d^\flat_F$, $ d^{*,\flat}_F$ and $\Delta^\flat_E$ with
$$
\Delta^\flat_E = (d^\flat_F + d^{*,\flat}_F)^2 =  d^\flat_F  d^{*,\flat}_F+  d^\flat_F d^{*,\flat}_F
$$
and commuting with $(d^\flat_F)^2 = 0,  ( d^{*,\flat}_F)^2 = 0$.

Traditionally, $d^\flat_F$ is denoted by $d_F$, so this convention
is followed here.

\vspace{.2in}
More explicitly,
let $A^p(M,C)$ be the smooth forms on $M$ with complex values,
i.e., smooth sections of the exterior power $\wedge^i (T^*M)$ of the
cotangent bundle $T^*M$, $A^p(M,C) = \Gamma(\wedge^i (T^*M)_R \otimes C)$.
The exterior derivative $d$ defines
a natural first--order operator
\begin{eqnarray}
d: A^p(M,C) \rightarrow A^{p+1}(M,C)
\end{eqnarray}
with  $d^2 = 0$,  giving the deRham complex
$(A^*(M,C),d)$.

Let $A^p(M,F)$ denote $p$-forms with values in the flat bundle
$F$, i.e.,$ A^p(M,F) = \Gamma(\wedge^i (T^*M)_R \otimes F) $.
Since $F$ is a flat bundle, with constant transition matrices,
 there is a unique first--order differential operator
$$
d_F \{= d_F^\flat \} : A^p(M,F) \rightarrow A^{p+1}(M,F)
$$
with the properties $d_F a \otimes s = da \otimes s$ over a open set $U$
for any complex valued $p$-form $a$  and flat section $s$ of $F$ over $U$.
One easily sees that $d_F^2 = 0$. The deRham theorem in this context
asserts that integration of forms gives a natural isomorphism
of cohomologies:
$$
H^p(A^*(M,F) , d_F)  \stackrel{\cong}{\rightarrow} H^p(M,F)
$$
where the last is the topologically defined  cohomology
with coefficients in the  flat coefficient
system $F$.

The chosen Riemannian metric $g$ determines
a unique star operator $\star$ from complex valued $p$-forms
to $n-p$ forms. It is induced by a bundle isomorphism:
\begin{eqnarray}
 \Lambda^p T^*M  \stackrel{\star}{\rightarrow}  \Lambda^{N-p} T^*M
\end{eqnarray}
with $<a(x),b(x)> \ dvol_x = a(x) \wedge \star b(x)$ at any point $x \in W$.
In terms of local coordinates $\{x_1,\cdots x_n\}$ in a neighborhood
of the  point $x \in W$ for which $\{dx_j\}$ at $x$ is an orthonormal
basis of $T^*M$ under $g$, one has  at $x$
$$
\begin{array}{l}
\star \  [( dx_{i_1}\wedge  \cdots \wedge  dx_{i_p} ) ]|_x
=    \ sign  \ [( dx_{k_1}\wedge  \cdots \wedge  dx_{k_{n-p}} )]|_x
\end{array}
$$
where $\{1,2,\cdots, N\} = \{i_1,\cdots, i_p\} \sqcup \{k_1,\cdots, k_{N-p}\}$.
Here the sign is the sign of the permutation
$[1,2,\cdots, n ] \mapsto  [i_1,\cdots, i_p,
k_1,\cdots, k_{n-p}, l_1']$.  Note  $\ \star^2  |_{A^p(M,C))}  = (-1)^{p(n-p)} \ Id $.

In particular,  a  linear star isomorphism of forms with values in $F$ is induced
by $\star \otimes Id_F$:
$$
\begin{array}{l}
A^{p}(M,F) =   \Gamma(\Lambda^p T^M  \otimes F)
  \stackrel{\ (\star \  \otimes \ Id_F)  \ }{\rightarrow}
 \Gamma(\Lambda^{n-p} T^*M  \otimes F) \\
\hspace{1.5in} = A^{n-p}(M,F)
\end{array}
$$

In these terms the  first--order differential operator $d_F^{* ,\flat}:
A^{p}(M,F) \rightarrow \newline  A^{p-1}(M,F)$ is  the composite:
\begin{eqnarray} \label{eqn.oper55}
     (-1)^{np + n+1} \ (\star \otimes Id_F) \
  d_F \ (\star \otimes Id_F) & :
\end{eqnarray}
Here by definition  $(d_F^{* ,\flat})^2 = 0$.

Consequently, one has a complex. $A^*(M,F)$ with two
differential $d_F, d_F^{* ,\flat}$ of degrees $+1$ and $-1$ so
the considerations of sections 3 to 5 apply with but small
changes.

It is important to note that if a smooth Hermitian inner product,
say $<.,.>_F$ is chosen on $F$, then the adjoint of $d_F$
is \textbf{not} in general the above mapping. Rather,  as specified by M\"{u}ller \cite{Muller2},
the adjoint is defined in terms of the induced conjugate linear
bundle isomorphism $\mu : F \rightarrow F^*$, of $F$ and its dual $F^*$ as:
\begin{eqnarray}
 \delta_F=
(-1)^{np+n+1} (\star \otimes Id_F)\ \mu^{-1} \ d_F \ \mu \ (\star \otimes Id_F) \label{eqn.operadj12}
\end{eqnarray}
With this definition $<d_F \  s, t>_F =
 <s, \delta_F \  t>_F$  for the induced
inner product on forms with values in $F$.
In most of the literature, this adjoint is rather sloppily written
as $d^* = (-1)^{np+p+1}\star d_F \star$ for suitable $\star$.

In particular, defining the flat Laplacian $\Delta^{\flat}_{F}$ by
\begin{eqnarray}
\Delta^{\flat}_{E} =( d_F + d_F^{* ,\flat})^2
   = d_F \ d_F^{* ,\flat}  +
         d_F^{* ,\flat} \   d_F  \label{eqn.defflatLap}
\end{eqnarray}
yields a generally non-self adjoint operator, as opposed to the
standard self-adjoint  Riemannian Laplacian coupled to $E$, $\Delta_{E}$, defined by
\begin{eqnarray}
\Delta_{F} = ( d_F +  d_F^{* ,\flat})^2
= d_F \  d_F^{* ,\flat} +  d_F^{* ,\flat}
\ d_F. \nonumber \\  \label{eqn.operadj22}
\end{eqnarray}

It is apparent that $\Delta^{\flat}_{F}$ has a simpler definition than
the classical self-adjoint $\Delta_{F}$ which uses the additional
information of a choice of Hermitian inner product, $<.,.>_F$, on  $F$.  In contrast
the
operator $\Delta^{\flat}_{E}$ utilizes only the Riemannian
inner product $g$ on $TM$.

However, if  $F$ is a unitarily flat
bundle then obviously:
\begin{eqnarray}
\Delta^{\flat}_{F} =\Delta_{F}
\ when \ F \ is \ unitarily \ flat
\end{eqnarray}

\vspace{.1in}

 On the other hand,    the symbols of the two first--order
operators $   \delta^{\flat}_F$ and
$  d^{\star,\flat}_F $ are identical, so they differ by a
smooth bundle mapping, say $B$,
$$
d^{\star,\flat}_F = \delta_F + B.
$$
In particular, $\Delta^\flat_F = \Delta_F + (B d_F + d_F B)$,so
 the generally non-self adjoint operator
$\Delta^{\flat}_{F}$ has the same symbol as the second--order
self adjoint operator $\Delta_{F}$, both acting on $A^{p}(M,F)$.

Consequently,  standard methods of elliptic theory apply,
e.g., the methods of Atiyah, Patodi, Singer and
Seeley, \cite{AtiyahPatodiSingerIII,Seeley1,Seeley2}.

\vspace{.2in}

As before, one has the  theorem
that the generalized eigenvalues of
$\Delta^{\flat}_F$ lie in a fixed parabola
about the positive real axis.

Let $K>0$ be a real number which is not the real part of
any eigenvalue.
Let  $S(p,K)$
be a complete enumeration of all the generalized eigenvalues counted with
multiplicities with real part greater than $K$. That is, $\Re (\lambda_j) >K$.

Then again the zeta function
\begin{eqnarray}
\zeta_{F,K,g}(s) = \Sigma_{\lambda_j \in S(p,K)} \ \frac{1}{\lambda_j^s}
\end{eqnarray}
defined using the principle part of the powers  converges
for $Re(s) > N$. If $\Pi_{F,K,g}$ denotes the spectral projection
on the span of the generalized eigenvectors with generalized eigenvalues
with real part less than $K$, then the
elliptic operator $(1-\Pi_{F,K,g}) \ \Delta^{\flat}_{F}$
fits the setting of Seeley \cite{Seeley1,Seeley2}. In particular, its complex
powers $[(1-\Pi_{F,K,g}) \ \Delta^{\flat}_{F}]^{-s}$ are well defined as in that
paper. Also as in \cite{Seeley2}, the ``heat kernel''
$e^{- t \ (1-\Pi_{F,K,g}) \ \Delta^{\flat}_{F} }$ is well defined for
$t >0$ real, it is of trace class, and for $Re(s) >N$ the
formula
\begin{eqnarray}
[(1-\Pi_{F,K,g}) \ \Delta^{\flat}_{F}]^{-s}
 = \frac{1}{\Gamma(s)}\int_0^{+\infty} \ t^{s-1} \ e^{- t \ (1-\Pi_{F,K,g})
 \ \Delta^{\flat}_{F}} \ dt \quad \quad
\end{eqnarray}
holds and its trace gives the zeta function $\zeta_{F,K,g}(s)$
for the bundle $F$.

Moreover, these  methods of Seeley imply that $\zeta_{F,K,g}(s)$
has a meromorphic extension to the whole plane and
is analytic at $s=0$. Consequently, the derivative
at $s=0$, \ $\zeta'_{F,K,g}(0)$ is meaningful.

Following Ray and Singer \cite{RaySinger2} and the above,  one sets
\begin{eqnarray}
Ray-Singer-Term(F, K,g)  \label{eqn.RaySingerdef2}\\
  = exp((1/2) \  \Sigma_{q=0}^N \ (-1)^q \ q \ \zeta'(F,K,g)(0)) \nonumber \\
\end{eqnarray}

\vspace{.3in}
In view of the fact that $d_F$ and $d_F^{*,\flat}$ commute with $\Delta^{\flat}_F$
the generalized eigensolutions with eigenvalue $\lambda \neq 0$
in the image of $d_F^{*,\flat}$  map by $d_F$ isomorphically
onto those in the image of $d_F$. As in Ray and Singer \cite{RaySinger1},
this implies complete cancelation in equation \ref{eqn.RaySingerdef2}
if the dimension of $M$ is even. Hence,
$
Ray-Singer-Term(F, K,g) =  0
$
in this case.

Consequently, in the remaining part of this section it is
assumed that $M$ is a closed \textbf{odd} dimensional smooth
manifold. For $n$ odd one has the equations:
$$
\begin{array}{l}
(\star \otimes Id_F)^2 = Id \ and  \\
\ d^{*,\flat}_F |_{A^p(M,F)} =  (-1)^p \ (\star \otimes Id_F) d_F
(\star \otimes Id_F) \hspace{.3in}
\end{array}
$$

In this setting
the methods of Ray and Singer \cite{RaySinger1} again apply to prove the following lemma.

\begin{lemma} \label{lemma.3A2} Suppose the smooth closed manifold
$M$ has odd dimension $n$.
Fix $K >0$. Let $g_t, \  a \le t \le b$ be a smooth family
of Riemannian  metrics on $TM$ such that the operators
$\Delta^{\flat}_{F},g_t$
 have no generalized
eigenvalues with real parts equal
to $K$. Here the dependence on $g_t$ is explicitly recorded.

Let the star operator $\star_t: $ be the bundle
isomorphism associated to the Riemannian inner product
$g_t$.
\begin{eqnarray}
 \Lambda^p T^*M  \stackrel{\star_t}{\rightarrow}  \Lambda^{n-p} T^*M
\end{eqnarray}
Let
$\Pi_{F,K,t}$ be the spectral projection onto  the span
of the generalized eigenvectors   of $\Delta^{\flat}_{F,g_t}$
with generalized eigenvalues of real part less than $K$.

Then  $Ray-Singer-Term(F, K,g)
= exp((1/2) \ \Sigma_{q=0}^N \ (-1)^q \ q \ \zeta'(F,K,g_t)(0))$
  varies smoothly with $t$ for $ a \le t \le b$ and similarly for $C^k$ and moreover,
\begin{eqnarray}
d/dt \
log \ det \ Ray-Singer-Term(F,K,g_t)
 \label{lemma.3.12} \\
 = (1/2) \ Tr (\Pi_{F,K,g_t} \  (\alpha_t \otimes Id_F)) \nonumber
\end{eqnarray}
where $\alpha_t =  \star_t^{-1} ( d/dt  \ \star_t ) $,
a bundle mapping acting on $\Lambda^*(T^*M)$.
\end{lemma}

\vspace{.1in}
Similarly, the methods of Ray and Singer straight-forwardly give the expected dependence
of the $Ray-Singer-Term(F, K,g)$
on changing $K>0$. It is:

\begin{lemma} \label{lemma.3B2}
Fix $L > K >0$. Let $g$  be a Riemannian  metric
on $TM$ such that the  operator $\Delta^{\flat}_{F,g}$ acting
on $A^{*}(M,F)$ has  no  eigenvalues with real part equal to
$K$ or $L$. Let $\{ \lambda_{p,j}| j = 1\cdots n_p\}$ be a complete
enumeration counting with multiplicities
of the generalized eigenvalues of $\Delta^{\flat}_{F,g}$ acting on $A^{p}(M,F)$
which have real part  in the range $K$ to $L$.
Then the formula below holds:
\begin{eqnarray}
\left[ \frac{Ray-Singer-Term(F, K,g)}{Ray-Singer-Term(F, L,g)}\right]^{\, 2}
= \left[ \prod_{p=0}^N \ (\prod_{j=1}^{n_p} \  \lambda_{p,j})^{(-1)^p  \, p} \right]^{-1}\quad \quad
\end{eqnarray}
\end{lemma}

\vspace{.1in}
Let $V(F,K,g)$ denote the span in $A^{p}(M,F)$
of the generalized eigensolutions of $\Delta^{\flat}_{F,g}$
of with generalized eigenvalues with real part
less than $K$. By elliptic theory $V(F,K,g)$ is a finite dimensional
subspace of smooth sections.

In view of the commutation relations,
$$
\begin{array}{lll}
\Delta^{\flat}_{F,g} d_F & =   d_F  \
d_F d^{*,\flat}_F \  d_F &  = d_F \
\Delta^{\flat}_{F,g}  \\
and & (\star \otimes Id_F)  \ \Delta^{\flat}_{F,g} &  =
\Delta^{\flat}_{F,g}   \ (\star \otimes Id_F) \nonumber
\end{array}
$$
the graded  complex
\begin{eqnarray}
C^*(F,K,g) := \bigoplus_{p=0}^N  \ V(p,F,K,g)
\nonumber
\end{eqnarray}
 has  two differentials, $d_F$ and $d^{*,\flat}_F = (-1)^{np+n+1} \ (\star \otimes Id_F)
 d_F (\star \otimes Id_F)$. Here $(-1)^{np+n+1} = (-1)^p$,
  $d_F :
V(p, F,K,g)\rightarrow  V(p+1,F,K,g) $  increases grading
by one  and $d^{*,\flat}_F :
V(p,F,K,g) \rightarrow  V(p-1,F,K,g) $
 decreases  degree by one.
That is, one has the same complex equipped with two differentials, a bi-complex:
\begin{eqnarray}
(C^*(F,K,g), d_F, d^{*,\flat}_F  ) \label{eqn.complex12}
\end{eqnarray}
of degrees $+1,-1$ respectively.

\vspace{.1in}
In  \S \ref{section.algebraictorsions}, it was  proved that
in this algebraic setting of a  graded complex of length $n$,\ $C^*$
with two differentials, $d$ of degree $+1$ and $d^*$ of degree $-1$,
i.e, $(C^*,d, d^*)$, there is
a natural non-vanishing algebraic torsion invariant
$$
torsion( C^*,d,d^*) \in (det
\ H^*(C^*,d) ) \otimes (det \ H_*(C^*,d^*))^* \ .
$$

The present example, $(C^*(F,K,g),d_F, \delta^{\flat}_F)$
fits this pattern. Hence, one gets a non-vanishing torsion invariant
$$
\begin{array}{l}
torsion((C^*(F,K,g),d_F, \delta^{\flat}_F) \\
\hspace{.4in} \in (det
\ H^*( C_p^*(F,K,g),d_F ) \otimes (det \
  H_*(C_p^*(F,K,g),d^{*,\flat}_F ) )^*
\end{array}
$$
These cohomologies and homologies can be identified as follows:

Let $N$ be the multiplicity of the generalized eigenvalue $0$
for $\Delta^{\flat}_{F}$
and $\cal{H}$ denote the span of these generalized eigenvectors
with generalized eigenvalue $0$. By spectral theory $\cal{H}$ is finite
dimensional and consists of smooth sections.  Then  one simply has the
direct sum decompositions:
$$
A^{*}(M,F) =  {\cal{H}} \oplus d_F  [(\Delta^{\flat}_{F})^M
A^{*}(M,F)] \oplus  \delta^{\flat}_F   [ (\Delta^{\flat}_{F})^N \  A^{*}(M,F)]
$$
and $d_F$ maps the second summand isomorphically to
the first. By this means one can easily show that the inclusion
of co-chain complexes
$$
(C^*(F,K,g),d_F )
\subset (A^{*}(M,F),d_F)
$$
induces an isomorphism on cohomology. That is,
$H^q( C^*(F,K,g), d_F ) \cong H^{p}(M,F)$.
 In particular, $det \ H^*( C^*(F,K,g), d_F )
 \cong det \ H^{*}(M,F)$ where $H^*(M,F)$ is the
topological cohomology of the local coefficient
system $F$.

Similarly, by $(\star \otimes Id_F)^2|_{A^{p}(M,F)} = \pm Id$, it follows that
$(\star \otimes Id_F) (\delta^{\flat}_F )
 (\star \otimes Id_F)^{-1} = \pm d_F$;
so $(\star \otimes Id_F)$ induces a  complex linear isomorphism
of the graded complex $C^*(F,K,g)$
to the complex $C^*(F,K,g)$ sending $C^p(M,F) \rightarrow C^{n-p}(M,F)$
with the differential $\delta^{\flat}_F $ becoming
$\pm d_F$.

In these terms, the  star operator $(\star \otimes Id_F)$  establishes a complex
linear isomorphism of homologies with cohomologies:
$$
\begin{array}{l}
 H_r(C^*(F,K,g), d^{*,\flat}_F  \ ) \\
\cong  H^{n-r}( C^*(F,K,g), d_F)  \\
\cong  H^{n-r}(M,F)
\end{array}
$$
In particular, one has a complex linear isomorphism of determinants:
$$
det \ H_*(C^*(F,K,g), \delta^{\flat}_F)
= (det H^{*}(M,F))^{(-1)^n}
$$

Hence, in this situation the algebraic invariants  may be regarded
as an element of the same complex line,
\begin{eqnarray}
torsion(C^*(F,K,g),d_F,
 \delta^{\flat}_F ), \nonumber \hspace{2in} \\
 \in \  (det  \ H^{*}(M,F)) \otimes
         [det \ H^{*}(M,F) ]^{(-1)^{n+1}} \hspace{0.2in}
\end{eqnarray}
for any real number $K >0$ and any choice of Riemannian  metric $g$ on $TM$.
Additionally, using Poincar\'{e} duality one has $det \ H^{n-*}(M,F)^* \cong det \ H^*(M,F^*)$.

For the case of interest $n$ is odd, so the algebraic torsion
takes its values in the square of a determinant line bundle:
\begin{eqnarray}
torsion(C^*(F,K,g),d_F,
  \delta^{\flat}_F), \nonumber  \\
 \in \  (det  \ H^{*}(M,F))^2  \cong det H^*( M, F \oplus F^*)
\end{eqnarray}

\vspace{.1in}
Now let $L>K>0$ be real and positive.
It  follows from the algebraic lemmas B,C of \S \ref{section.algebraictorsions},
that under the assumptions of lemma \ref{lemma.3B2},  the algebraic
torsion of
$( C^*(F,K,g), d_F, \delta^{\flat}_F)$
and that for $L$ with  $K>L>0$ are precisely related by the eigenvalues $\lambda_{q,j}$
of lemma \ref{lemma.3B2} as follows:
\begin{eqnarray}
 torsion( C^*(F,L,g), d_F, \delta^{\flat}_F)) \nonumber \\
=   torsion ( C^*(F,K,g), d_F, d^{*,\flat}_F))  \label{eqn.stabilize12}\\
\hspace{.5in} \cdot \left[ \prod_{p=0}^N \ (\prod_{j=1}^{n_p} \  \lambda_{p,j})^{(-1)^p  \, p}
\right]^{-1}
\nonumber
\end{eqnarray}
when regarded as an element of $\  (det  \ H^{p,*}(M,F))^2$. Here $n$ is odd.

Similarly, it  follows just as in the proof of  lemma \ref{lemma.5vary} of \S 6
 for a smooth
family of Riemannian  metrics $g_t$ satisfying the  assumptions of lemma \ref{lemma.3A2}
one gets the equality
\begin{eqnarray}
d/dt \left[ \ log \ det \ torsion(C^*(F,L,g_t), d_F, d^{*,\flat}_F )
\right] \nonumber \\
  = - \ Tr (\Pi_{F,K,g_t} \  (\alpha_{g_t} \otimes Id_F)) \quad \quad
\end{eqnarray}

As an immediate consequence of the above two equations ,
the main theorem below follows:

\begin{thm}[Independence of metric theorem ] \label{thm.MAIN2} Let $M$ be closed, smooth, oriented
and of odd dimension $n$.
For any choice of Riemannian  inner product $g$ on
$TM$, pick a real number $K>0$ for which the
operators $ \Delta^{\flat}_{F,g}$ acting on $A^{*}(M,F)$ have
no eigenvalue with real part equal to $K$ and similarly for the
trivial bundle $C^k$. Define
the graded complex $C^*(F,K,g) = \bigoplus_{p=0}^N  \ V(p,F,K,g)$ with its
two differentials, $d,d^*$ as above and similarly for $C^k$. Then
the combination
$$
\begin{array}{lll}
\tau_{analytic}( M,F)\  := &
 torsion((C^*(F,K,g),d_F,d^{*,\flat}_F)
 \\
& \hspace{.5in} \cdot
Ray-Singer-Term(F,K,g)^{\, 2}\\
\end{array}
$$
is independent of the choice of $K>0$ and also independent of
the choice of Riemannian  inner product $g$ on $TM$ chosen.

By definition $\tau( W,F)$ is a non-vanishing element
of the squared determinant line bundle
\begin{eqnarray}
\tau(M,F) \  \in (det  \ H^{*}(M,F))^2    \cong det H^*( M, F \oplus F^*)
 \nonumber
\end{eqnarray}

\end{thm}

\textbf{Addendum 1):}
\textit{In the acyclic cases the torsion is a complex number.}

\textbf{Addendum 2):} \textit{In the case that
$F$ is acyclic and flat unitary, this torsion is a real number.}

Since the operator $\Delta^{\flat}_{F}$ is self-adjoint
in this case, one may take $K>0$ less than the smallest non-vanishing
real eigenvalue. For this choice of $K>0$, the algebraic correction term
is just $+1$. For $F$ flat unitary and acyclic
$\tau(W,F)$ is identical with  the square of Ray Singer analytic torsion \cite{RaySinger1}
and by the Cheeger-M\"{u}ller theorem \cite{Cheeger1,Muller1} equals
the square of the Riedemeister-Franz topological
torsion in this special flat unitary acyclic case.

As explained, the proof is completely parallel to the d-bar case,
except there is no global anomaly term since the dimension $n$ is odd.

% ******************************************

\section{Combinatorial torsion for general flat bundles over compact smooth manifolds}
  \label{section.comb}

Let $F$ be a flat bundle over a general compact smooth oriented manifold $M$, possibly with boundary.
In this section a ``Reidemiester- Franz'' type of torsion
$$
\tau_{comb}( M,F) \in  det \ H^*(M,F \oplus F^*)
$$
is defined in a combinatorial fashion, an invariant of the pair $(M,F)$.
Here $F^*= Hom(F,C)$ is the dual bundle to $F$.
The method is very straightforward.

Let $F$ be a flat complex bundle over an oriented  n--dimensional manifold with a
cellular  subdivision, say  $W = \{\sigma\}$ with dual cellular complex, say $D(W)= \{ D(\sigma\})$.
 Here  for $\sigma$ an oriented cell, let $D(\sigma)$ is  the dual cell with the
corresponding compatible
orientation. Let $K(\sigma)$ be the
flat sections over $\sigma$ and $K(D(\sigma))$ be the flat sections over $D(\sigma)$.
One  gets a \textbf{canonical} isomorphism, by restricting to the common center:
$$
K(\sigma) \cong K(D(\sigma))
$$

Hence, summing over oriented cells and dually oriented cells, we get a canonical
isomorphism summing over cells of dimension $a$
$$
\Theta : C^a(W,F) = \oplus \ K(\sigma) \cong \oplus \ K(D(\sigma)) = C^{n-a}(D(W),F)
$$
Let $d,d'$ be the coboundaries for $W$, $D(W)$ respectively.

By means of $\Theta$ one  may identify these, and so endow $C^*(W,F)$ with two differentials
on the same complex $C^*(W,F)$:
$$
d,\  going \ up, \ \  and  \ \delta = \Theta^{-1} \ d' \ \Theta, \ \  going \ down
$$  and so define  an algebraic torsion $\tau_{comb.}(C^*(W,F),d, \delta)$
as in section \ref{section.algebraictorsions}.
By  methods similar to those used in  Milnor \cite{Milnor1} $\tau_{comb.}(C^*(W,F),d, \delta)$, can be shown to
be independent of the choice of cellular decomposition.
$$
\tau(M,F) = \tau_{comb.}(C^*(W,F),d, \delta)
$$
The cohomology of the cell complex
$W$ gives $H^*(M,F)$;  the cohomology of the cell complex $D(W)$ gives the relative cohomology
 $H^*(M, \partial M, F)$.

This element, $ \tau(M,F)$
 , by Poincar\'{e} duality  lies in
$$
\begin{array}{l}
\tau(M,F) \in \ det(H^*(W,F))\otimes det( H^{n-*})(D(W),F)^* \\
\cong  det(H^*(M,F)) \otimes det(H^{n-*}(M,\partial M, F))^* \\
 \cong  det(H^*(M,F)) \otimes det(H^{*}(M, F^*))^* \\
\cong det(H^*(M,F \oplus F^*))
\end{array}
$$

 Note on the non-oreinted case: In the case $M$ is not oriented, then the isomorphism $\Theta$
becomes an isomorphism, $\Theta' : C^a(W,F) \cong C^{n-a}(D((W), F \otimes O)$
where $O$ is the orientation sheaf. In this
case the resultant torsion lies in $det(H^*(W,F))\otimes det( H^{n-*})(D(W),F \otimes O)^*
\cong  det(H^*(M,F)) \otimes det(H^{n-*}(M,\partial M,F \otimes O))^*
\cong  det(H^*(M,F)) \otimes \newline det(H^*(M,F^*))
\cong det(H^*(M,F \oplus F^*))$ again.

\vspace{.1in}
A more functorial description  of this element $ \tau(M,F) = \tau_{comb.}(C^*(W,F),d, \delta)$
is as follows: We assume for convenience $M$ is oriented.

The isomorphism $\Theta^{-1}$ defines an element
$$
[\Theta^{-1}] \in (det \ C^*(W,F)) \otimes  ( det \ C^{n-*}(D(W),F))^*
$$
so under the Bismut-Zhang isomorphisms \cite{BismutZhang1}, $\tau \otimes \tau$,  one obtains an element
$$
\begin{array}{l}
(\tau \otimes \tau)[\Theta^{-1}] \in  (det \ H^*(W,F)) \otimes  ( det \ H^{n-*}(D(W),F))^* \\
= (det \ H^*(M,F)) \otimes ( det \ H^{n-*}(M,\partial M, F))^* \cong det \ H^*(M,F \oplus F^*)
\end{array}
$$
where the last step used Poincar\'{e} duality. In these terms one has
$$
\tau_{comb}(M,F) = \tau_{comb.}(C^*(W,F),d, \delta) = (-1)^{S(C)} \ (\tau \otimes \tau)[\Theta^{-1}]
$$
with the sign $(-1)^{S(C)}$ as above.

As a special case of this construction, we may take a Morse function
with its associated downwards and upwards cells.

\section{Comparison of Combinatorial and Analytic Torsions: Generalization of the Cheeger--M\"{u}ller theorem:}
\label{section.mainCheegerMuller}

Let $F$ be a k dimensional flat complex bundle over a smooth oriented odd dimensional Riemannian manifold $M$ of dimension $n$
with Riemannian metric $g$. Let

Chose $K >0$ as in section \ref{section.analyticRiem1} so it is not the real part of any eigenvalue
of the flat Laplacian $\Delta_F^\flat$ and let $V^*(F,K,g)$ the span of the smooth forms which
have eigenvalues with real part less than or equal to $K$.

Now  $d_F$ and $\hat{\star} d_F \hat{\star}$ both commute with  $\Delta_F^\flat$,
and
$V^*(F,K,g)$ acquires by restriction two differentials  and so one may form the
associated algebraic torsion as in \S \ref{section.algebraictorsions}. The torsion $\tau_{analytic}(M,F,g)$ is obtained
from this by multiplying by the square of the Ray-Singer term.

A more functorial way to express this analytic torsion is as follows:
Consider the star operator $\hat{\star} = \star \otimes Id_F$, which gives an isomorphism:
$$
\hat{\star} : V^*(F,K,g) \cong  V^{h-*}(F,K,g)
$$
This isomorphism defines an element
$$
[\hat{\star}] \in ( det \ V^*(F,K,g)) \otimes ( det \  V^{n-*}(F,K,g))^*
$$
Now using the differentials, $d_F, d_F$ on each, by the Bismut-Zhang torsion
isomorphism $\tau$ \cite{BismutZhang1}, one obtains an element:
$$
( \tau \otimes \tau)[\hat{\star}] \in ( det \ H^*(V^*(F,K,g),d_F)) \otimes ( det \ H^{n-*}( V^{*}(F,K,g),d_F)^*
$$
Now integrating these forms over the manifold, by the deRham theorem gives an isomorphism:
$$
\begin{array}{l}
\int \otimes \int : ( det \ H^*(V^*(F,K,g),d_F)) \otimes ( det \ H^{n-*}( V^{*}(F,K,g),d_F)^* \\
\cong  ( det \ H^*(M,F) \otimes ( det \ H^{n-*}(M,F)^* \cong det \ H^*(M,F \oplus F^*)
\end{array}
$$

From this perspective, the analytic torsion $\tau_{analytic}(M,F,g)$ is defined by
$$
\begin{array}{ll}
\tau(M,F,g) &= (-1)^{S(C)} \ [ Ray-Singer-Term(p,E, K,g) ]^2 \  \\
& \times (\int \otimes \int)[ ( \tau \otimes \tau)[\hat{\star}] ] \in det \ H^*(M,F \oplus F^*)
\end{array}
$$
As observed in section \ref{section.analyticRiem1}, it is independent of the metric $g$.

\vspace{.3in}
\begin{thm}[Generalized Cheeger M\"{u}ller Theorem] \label{thm.mainCheegerMuller}
For a closed, oriented, odd dimensional manifold $M$,
the two torsion elements $\tau_{analytic}(M,F,g)$ and $\tau_{comb}(M,F)$ are equal
in $det \ H^*(M,F \oplus F^*)$.
\end{thm}

This theorem was a conjecture in an earlier version of this manuscript.
Following the methods of Su and Zhang,  a direct proof is given below.
Note that for the special case where $F$ admits a smoothly varying
non-degenerate inner product, a mild constraint,  this theorem was  proved by Braverman and
Kappeller \cite{BravermanKappeller3}. However, the following
approach has a certain attractive naturality.

\vspace{.1in}

Chose a Morse function $f$ on $M$ which satisfies the Morse-Smale conditions
with critical points the finite set $B$, downward or descending  cell decomposition  $W^u = \{W_i^u\}$ and
dual upward or ascending  cell decomposition  $D(W_u) = W^s= \{ W_i^s \}$, $W_i^s = D(W_i^u)$. Then the integration
maps and evaluation

$$
\begin{array}{l}
\int \otimes \int : ( det \ H^*(V^*(F,K,g),d_F)) \otimes ( det \ H^{n-*}( V^{*}(F,K,g),d_F)^* \\
\cong  ( det \ H^*(M,F) \otimes ( det \ H^{n-*}(M,F)^* \cong det \ H^*(M,F \oplus F^*)
\end{array}
$$

can be concretely realized by integrating over the cells and dual cells as in:
$$
\begin{array}{l}
\int \otimes \int : ( det \ H^*(V^*(F,K,g),d_F)) \otimes ( det \ H^{n-*}( V^{*}(F,K,g),d_F)^* \\
\cong  ( det \ H^*(C^*(W^u,F) \otimes ( det \ H^{n-*}(C^*(D(W^u) ,F)^* \cong det \ H^*(M,F \oplus F^*)
\end{array}
$$

This allows the reformulation of theorem \ref{thm.mainCheegerMuller} as a comparison
result about the relation of the star operator on forms and the duality mapping
$\Theta$ of cells to dual cells.

\begin{thm}[Analytic - Geometric version:] \label{thm.alternative}
The
element $(-1)^{S(C)} \ [ Ray-Singer-Term(p,E, K,g)]^2 \ [\hat{\star}^{-1}] $
regarded as an element of \newline
$(det \ V^*(K,F,g)) \otimes ( det \ V^{N-*}(V,F,K))^* \stackrel{\tau
\newline
\otimes \tau}{\cong} \newline
(det \ H^*(V^*(K,F,g),d_F) )
\otimes ( det \ H^{n-*}(V^*(V,F,K),d_F))^*$
 passes by integration, $\int \otimes \int$  to the element
of \newline
$( det \ H^*(C^*(W^u,F)) \otimes ( det \ H^{n-*}(C^*(D(W^u),F))^*$ given by
the duality mapping $(-1)^{S(C')} \ [ \Theta^{-1}]$ regarded as an element of \newline
$
( det \ C^*(W^u,F)) \otimes ( det \ (C^{n-*}(D(W^u),F))^*
\stackrel{\tau \otimes \tau}{\cong} \newline ( det \ H^*(C^*(W^u,F))
\otimes ( det \ H^{n-*}(C^*(D(W^u),F))^*$

\end{thm}

\vspace{.1in}
Now introduce a  Witten--type deformation \cite{Witten1} of our present, non--self adjoint,  Laplacian .
That is, introduce for each real number $t$ the operators on $A^*(M,F)$

$$
\begin{array}{l}
\tilde{\Delta}_{T,f,F} = d_{T,f,F} \ \delta_{T,f,F} + \delta_{T,f,F} \ d_{T,f,F} = (\tilde{D}_{T,f,F})^2 \\
with \\
\tilde{D}_{T,f,F} =  d_{T,f,F} + \delta_{T,f,F} \\
d_{T,f,F} = e^{-T  f} d_F e^{+T  f}, and , \ \delta_{f,t,F} = e^{T  f} \delta^{\#}_F  e^{-T  f} \\
\end{array}
$$
and  define the deformed version of $d_F + \delta^{\#}_F$.
$$
D_{T,f,F} := e^{T f} \ \tilde{D}_{T,f,F} \ e^{-T f} = d_F + e^{2 T f}\  \delta^{\#}_F \ e^{-2 T f}
$$
$\tilde{\Delta}_{T,f,F}$ is called a Witten--type  Laplacian.

Note that $D_{T,f,F}$ commutes with $d_F$. Let $\delta^{\#}_{T,f,F} = e^{2 T f} \delta^{\#}_F e^{-2 T f}$,
so the deformed Laplacian $\Delta^\flat_{T,f,F} := (D_{T,f,F})^2 = d_F \ \delta^{\#}_{T,f,F} + d_F \ \delta^{\#}_{T,f,F} $
commutes with $d_F$. It is differs by a first order operator from $\Delta_F^\flat$,  so has all the standard properties.

For $K>0$ not the real part of any eigenvalue of the elliptic operator $\Delta^\flat_{T,f,F}$,
let $V^*(T,f,F,K,g)$ denote the span of the space of generalized eigenforms with real parts less than $K$,
 $\Pi^*(T,f,F,K,g)$ the spectral projection on that finite subspace and
$$
Ray-Singer-Term(T,f,F,K,g)
$$
the Ray and Singer contribution from the eigenvalues with real value greater than $K$ defined
as before via zeta function regularization.

Since this deformation merely changes $\hat{\star} d_F \hat{\star}^{-1} $ into
$ e^{2 T f}\ \hat{\star} d_F \hat{\star}^{-1} \ e^{2 T f}$, the now standard argument
shows that the derivative   $d/dT \ log [ Ray-Singer-Term(T,f,F,K,g)]^2$ is precisely
$$
- \Sigma_{q=0}^n \ (-1)^q \  Tr( \Pi^q(T,f,F,K,g) \  \cdot f )
$$
The absence of anomaly comes from the fact that $M$ is odd dimensional.

\vspace{.1in}
On the other hand, the operator $e^{-2 T f} \hat{\star}$ intertwines the deformed Laplacians
for the values $T,-T$:
$$
[e^{-2 T f} \ \hat{\star}] \  \Delta^\flat_{T,f,F} = \Delta^\flat_{-T,f,F} [ e^{-2Tf} \ \hat{\star}]
$$
So the deformed Laplacians $\Delta^\flat_{T,f,F}$ and $\Delta^\flat_{-T,f,F}$ have the same eigenvalues
and the induced mapping
$$
[e^{-2 T f} \ \hat{\star}] : V^*(T,f,F,K,g) \cong V^{n-*}(-T,f,F,K,g)
$$
identical with the above construction for $T=0$. These subspaces are cochain complexes
under $d_F$ which commutes with $\Delta^\flat_{T,f,F}$ and $\Delta^\flat_{-T,f,F}$.

Consequently, one may again form the associated element of
$ (det H^*(W_u,F) ) \otimes (det \ H^{n-*}(D(W_u),F))^*  \cong \det \ H^*(M,F \oplus F^*)$
as before.
$$
\begin{array}{l}
[e^{-2 T f} \ \hat{\star}  ]^{-1} \in ( det \  V^*(T,f,F,K,g)) \otimes  ( det V^{n-*}(-T,f,F,K,g))^* \\
\stackrel{\tau \otimes \tau}{\cong} ( det \ H^*( V^*(T,f,F,K,g),d_F) \otimes  ( det H^{n-*}(V^*(-T,f,F,K,g),d_F))^* \\
\stackrel{\int \otimes \int}{\cong} (det \ H^*(M,F))\otimes (det \ H^{n-*}(M,F))^* \cong det \ H^*(M,F \oplus F^*)
\end{array}
$$

Since $e^{-2Tf}$ merely deforms the star operator, again the now standard methods
show that the variation of $(-1)^{S(C)} \ (\int \otimes \int) ( \tau \otimes \tau) ([e^{-2 T f} \ \hat{\star}  ]^{-1})$
with $T$ is precisely $
+e^{ 2 \Sigma_{q=0}^n \ (-1)^q \  Tr( \Pi^q(T,f,F,K,g) \  \cdot f )}
$. This together with the above proves the lemma.

\begin{lemma}[Invariance under Witten deformation:] \label{lemma.witteninv}
The element
\newline $[ Ray-Singer_ Term(T,f,F,K,g)]^2 \
 (\int \otimes \int) ( \tau \otimes \tau) ([e^{-2 T f} \ \hat{\star}  ]^{-1})$ \newline
regarded as an element of $det \ H^*(M,F \oplus F^*)$ is independent of $T$.
That is, its equal to $\tau_{analytic}(W,F,g)$ for all $T$.
\end{lemma}

An alternative formulation of this invariance  is that the
element
$$
\begin{array}{l}
[ Ray-Singer_ Term(T,f,F,K,g)]^2 \ \
 (\int \otimes \int) ( \tau \otimes \tau) ([e^{-2 T f} \ \hat{\star}  ]^{-1}) \\
\in (det H^*(W_u,F)) \otimes ( det H^*(D(W_u), F)^*
\end{array}
$$
is independent of $T$. Here the integrals are carried out over the cells of $W_u$
and the cells of $W_s$ respectively.

To analyze the low eigenvalues of the deformed Laplacian $\Delta_{T,f,F}^\flat$
proceed following Bismut and Zhang and Su and Zhang \cite{BismutZhang1}, \cite{SuZhang}.

For each critical point, say $ p \in B$, of index $n_p$,chose local coordinate
$y^1,\cdots,y^n$ so that in this neighborhood $U_p$
$$
f = f(p) - [(1/2) \Sigma_{i=1}^{n_p} (y^i)^2] + (1/2) [\Sigma_{j=n_p+1}^n \ (y^j)^2]
$$
compatible with the orientation of $M$.
Let $W^u_p$ denote the unstable descending cell from $p$, and $W_p^s$ the stable ascending cell
with their inherited orientations from this choice of coordinates.
Take a Hermitain metric on $F$ which is flat on $U_p$ and use the flat sections
to identify $F|U_p \cong U_p \times F_p$.

Set $|Y|^2
= \Sigma_{i=1}^n  \ (y^i)^2$.
Let $\gamma(t)$ be non--negative, smooth compactly supported function identically 1 nearby $0$
for which $\gamma(|Y|)$ has compact support
 in the coordinate chart $U_p$ above.

For an element $e \in F_p$ introduce the following  Gaussian $n_p$--form:
$$
\begin{array}{l}
\rho_{p,T}(e)  = \frac{\gamma(|Y|)^2}{\sqrt{\alpha_{p,t}}}
 \ \times exp(-\frac{T |Y|^2}{2})  \otimes e
 dy^1 \wedge \cdots \wedge dy^{n_p} \\
where \\
\alpha_{p,t} = \int_{U_p} \ \gamma(|Y|)^2 exp( - T |Y|^2) dy^1 dy^2 \cdots dy^n
\end{array}
$$

By sending $[W_p^u] \otimes e $ to $\rho_{p,T}(e) \in A^{n_p}(M,F)$,  define
a mapping of the descending cochain complex to forms:
$$
J_{f,T,u} : C^*(W^u,F)  \rightarrow A^*(M,F)
$$

Similarly, replacing $f$ by $-f$  defines the  mapping for  the ascending
cochain complex:
$$
J_{f,T,s} : C^*(W^s,F) \rightarrow A^*(M,F)
$$

It is a crucial observation that the standard duality isomorphism
$$
\Theta : C^a(W^u,F) \cong C^{n-a}(W^s,F)
$$
sending $[W_i^u]\otimes e$ to $[W_i^s]\otimes e$ and the star operator
$$
\hat{\star} : A^a(M,f) \cong A^{n-a}(M,F)
$$
form the commutative diagram:

\begin{lemma}[$\Theta$ and $\star$ compatibiltity] \label{lemma.compat}
$$
\begin{array}{llll}
J_{T,f,u} :&  C^*(W^u,F) & \rightarrow &  A^*(M,F)  \\
 &\downarrow  \Theta & &  \downarrow \hat{\star} \\
J_{T,f,s} :&  C^*(W^s,F) & \rightarrow &  A^*(M,F)  \\
\end{array}
$$
\end{lemma}

That is, the geometric duality map $\Theta$ becomes
precisely the star operator under these two mappings.
This observation was emphasized by Bismut and Zhang \cite{BismutZhang1}

Recall the from the definitions, the Witten Laplacian is:
$$
\tilde{\Delta}_{T,f,F} := e^{-T f} \ \Delta_{T,f,F} \ e^{+T f}
$$
It has the property that $\hat{\star} \tilde{\Delta}_{f,T,F} = \hat{\star} \ \tilde{\Delta}_{f,-T,F}$.

Let $A_{[0,1],T,f}^*(M,F)$ be the span of the generalized eigenvectors
of $\tilde{\Delta}_{f,T,g}$ with eigenvalues $\lambda$ with
$|\lambda| <1$. Let $\Pi_{f,T}$ denote the spectral projection
of $A^*(M,F)$ onto $A_{[0,1],T,f}^*(M,F)$. Similarly for $f,-T$.

In particular , $a \rightarrow e^{tf} a$ carries $A_{[0,1],T,f}^*(M,F)$ isomorphically
to \newline $V^*(T,f,F,K=1,g)$.

The basic result is the following theorem, see Bismut and Zhang \cite{BismutZhang1}

\begin{thm} \label{thm.approx}
For $T$ sufficiently large there is a
$c>0$ with
$||J_{f,T} - \Pi_{T,f} J_{f,T}|| = 0 (e^{-cT})$
and  $||J_{f,-T} - \Pi_{f,_T} J_{f,-T}|| = 0 (e^{-cT})$.

Moreover, $ \Pi_{f,T} J_{T,f},  \Pi_{f,-T} J_{f,_T}$ are isomorphisms.

Moreover, $\hat{\star} \tilde{\Delta}_{f,T} = \tilde{\Delta}_{t,T}  \hat{\star}$
so $\hat{\star} $ induces an isomorphism $A_{[0,1],f,T}^*(M,F) \rightarrow A_{[0,1],f,-T}(M,F)$.
 Also
 the following commutes
$$
\begin{array}{llll}
J_{f,T,u} :& C^*(W^u,F) & \rightarrow &  A_{[0,1],T,f}^*(M,F)  \\
&\downarrow \Theta  & &   \downarrow \hat{\star} \\
J_{f,T,s} :& C^*(W^s,F) & \rightarrow &  A_{[0,1],T,-f}^*(M,F)  \\
\end{array}
$$

\end{thm}

These relations can be explored using the commutative diagram
with all horizontal arrows isomorphisms for $T$ sufficiently large,

$$
\begin{array}{lllllll}
C^*(W^u,F) &\rightarrow&  A_{[0,1],f,T}^*(M,F) &\stackrel{ \times e^{Tf}}{\rightarrow}&
V^*(f,T,F,K,g) &\stackrel{\int}{\rightarrow}& C^*(W_u,F) \\
\\
\downarrow \Theta & & \downarrow  \hat{\star} & & \downarrow  e^{-2Tf} \hat{\star}& & \downarrow   \\
\\
C^*(W^s,F) &\rightarrow&  A_{[0,1],f,-T}^*(M,F) &\stackrel{ \times e^{-Tf}}{\rightarrow}& V^*(f,-T,F,K,g)
 &\stackrel{\int}{\rightarrow}& C^*(W_s,F) \\
\end{array}
$$
where the first maps are  the isomorphisms:
$\Pi_{f,T} \ J_{f,T,u}: C^*(W^u,F)  \rightarrow  \newline A_{[0,1],f,T}^*(M,F)$ and
$ \Pi_{f,-T} J_{f,T,s}: C^*(W^s,F)  \rightarrow   A_{[0,1],f,-T}^*(M,F)$.

By the above, the  algebraic torsion contribution to $\tau_{analytic}(M,f,T,g)$
of the eigenspace $V^*(T,f,F,K,g)$ and its two differentials is the same as
taking the isomorphism $[e^{-2Tf} \hat{\star}]^{-1}$, regarding it as an element
of $(det \ V^*(f,T,F,K,g)) \otimes ( det \ V^{n-*}(f,-T,F,K,g))^*$ and applying the
integration mapping to send it to $(det \ C^*(W^u,F) ) \otimes ( det \ C^{n-*}(W^s,F))^*$.
By this commutative diagram this may be done by starting with the element $[\Theta]^{-1}$
as an element of $(det \ C^*(W^s,F)) \otimes ( det \ C^{n-*}(W^s,F))^*$  and passing it
along from the left to the right.

Now by theorem \ref{thm.approx} this composite is approximated with error $(e^{-ct})$
by the composite coming from modified diagram:
$$
\begin{array}{lllllll}
C^*(W^u,F) & \stackrel{ \ J_{f,T,u}}{\rightarrow} &  A^*(M,F) & \stackrel{ \times e^{Tf}}{\rightarrow}
  & A^*(M,F) &  \stackrel{\int}{\rightarrow}  & C^*(W_u,F) \\
\\
\downarrow \Theta & & \downarrow  \hat{\star} & & \downarrow  e^{-2Tf} \hat{\star}& & \downarrow   \\
\\
C^*(W^s,F) & \stackrel{ \ J_{f,T,u}}{\rightarrow} &  A^*(M,F) & \stackrel{ \times e^{Tf}}{\rightarrow}
  & A^*(M,F) &  \stackrel{\int}{\rightarrow}  & C^*(W_s,F) \\
\end{array}
$$
This computing these integrals yields the following lemma.
It is the analogue of Theorem 3.3 of the paper of Su and Zhang.

\begin{lemma} \label{lemma.funda}
$$
\begin{array}{l}
\lim_{T \rightarrow +\infty}  \frac{   \tau( e^{T f}  \ V(K=1,F,g),d_F, \hat{\star}  \ d_F \ \hat{\star} )}{
\tau(C^*(W^u,F), d_F, \Theta^{-1} d_F \Theta)}
\\
\cdot ( \frac{T}{\pi})^{((n/2)r - s)}
exp( 2 k \Sigma_p  \ (-1)^{n_p} \ f(p)) \\
\\
= 1
\end{array}
$$

 with $r = k ( \chi(M)) = k ( \Sigma_p \ (-1)^{n_p}), s = k (\Sigma_p \ (-1)^{n_p} f(p))$.

\end{lemma}

Combining the Witten invariance lemma \ref{lemma.witteninv}
with this last lemma, it follows that theorem \ref{thm.alternative}
and so theorem \ref{thm.mainCheegerMuller} will be established once
it is proved that

$$
\begin{array}{l}
lim_{T \rightarrow  \ +\infty} \ [( \frac{T}{\pi})^{((n/2)r - s)}  exp( 2 k \Sigma_p  \ (-1)^{n_p} \ f(p))]^{-1} \\
 \  (Ray-Singer-Term(T,f,F,K=1,g))^2 = 1 \\
\end{array}
$$
Here the eigenvalues can be taken as those of the Witten Laplacian
$\tilde{\Delta}_{f,T,g}$ which has the same eigenvalues as
$\Delta_{f<T,g}$

In the paper of Su and Zhang \cite{SuZhang} the analogous result is proved for the
Witten version of the operator of Haller and Burghelea. Fortunately
 their  methods are very robust and apply without essential change to the Witten--type
Laplacian above. The only major change is the above control on the small eigenvalues.
This concludes the proof of the generalized
Cheeger M\'{u}ller result \ref{thm.mainCheegerMuller}.

\section{Zeta functions  for the d-bar setting:} \label{Analysis} \label{section.Analysis}

The analytical details are given in the context of the $\bar{\partial}$ Laplacian
coupled to a  holomorphic bundle $E$ with compatible type $(1,1)$ connection.
The analogous results for  the Riemannian Laplacian coupled to
a flat complex bundle follow the same outline, so their proofs are omitted.
All this material is relatively standard and is often used without comment  in the literature.
The basic reference is Ray and Singer \cite{ RaySinger1, RaySinger2}.

Let $E$ be a holomorphic  bundle  with compatible connection $D$ of type $(1,1)$
over a complex Hermitian manifold $W$. Let a Hermitian metric $<.,.>_E$ be chosen
for $E$ and the adjoint of the d-bar operator under the induced inner product
on $A^{p,q}(W,E)$ be denote by $\delta_E = adjoint(\bar{\partial}_E)$. Let the
associated self--adjoint d-bar Laplacian be denoted by $\Delta_E$.
$$
\Delta_E = (\bar{\partial}_E+ adjoint(\bar{\partial}_E))^2= ( \bar{\partial}_E + \delta_E)^2
= \bar{\partial}_E \ \delta_E +  \delta_E \ \bar{\partial}_E
$$

In this section fix $p, 1 \le p \le N$ and throughout consider
all operators acting on $A^{p,q}(W,E)$ for $p$ fixed
with possibly varying $q$.

Recall from equation \ref{eqn.alpha},
$
\bar{\partial}^*_{E,D''} = \delta_E + \alpha
$
for a smooth bundle mapping $\alpha :\Lambda^q(T^*M) \rightarrow \Lambda^{q-1}(T^*M)$.
Moreover,
$
\Delta_{E, \bar{\partial}} = (\bar{\partial} + \bar{\partial}^*_{E,D''})^2 =
(\ \bar{\partial}_E +  adjoint(\bar{\partial}_E) )^2
 + \alpha \bar{\partial}_E +  \bar{\partial}_E \alpha
= \Delta_E + \alpha \bar{\partial}_E +  \bar{\partial}_E \alpha$
hence is a elliptic second order partial differential equation with
scalar symbol.

\vspace{.1in}

Proof of lemma \ref{lemma.AN1}:

Use the  reference  choice of Hermitian inner product $<.,.>_E$
on $E$ to define all norms, Solobev spaces, etc.
This is used to construct the adjoint $\delta_E$ and
the self-adjoint d-bar Laplacian $\Delta_E$ as above.

Since $G_t = (\lambda I - (\bar{\partial}_E+\delta_E + t \alpha)), 0 \le t \le 1$,  \ is a first--order
elliptic operator, it is Fredholm regarded
as a map of Sobolev spaces $A^{p,*}(W,E)_{(1)} \rightarrow A^{p,*}(W,E)_{(0)}$,
losing one derivative.

As $G_0$ is self-adjoint, it has index zero. Hence by
homotopy invariance of the index, $G_1 = (\lambda I - (\bar{\partial} + \bar{\partial}^*_{E,D''}))$
is Fredholm of index zero. In particular, if one  shows that
the kernel, $ker \ G_1$, is zero, then $G_1$ is also onto and
so an isomorphism. That is,  $\lambda$ will be not be in the
spectrum of $(\bar{\partial}_E +  \delta_E  + \alpha)  =
\bar{\partial} + \bar{\partial}^*_{E,D''}$ unless there
is an eigenvector with eigenvalue $\lambda$.

Since $\bar{\partial}_E + \delta_E$ is self-adjoint, we may chose an orthonormal
basis, say $\{\phi_j, \  j=1,\cdots\}$, of eigenvectors for
the $L^2$ completion of $A^*(M,E)$. Let $(\bar{\partial}_E + \delta_E) \phi_j =
\mu_j \ \phi_j$ for real non-negative numbers $\mu_j$,

Given  $a = \Sigma_j a_j \phi_j$ with $||a||^2 =\Sigma_j \ ||a_j||^2 = 1$, one  has
the simple estimates $||(\lambda I - (\bar{\partial}_E + \bar{\partial}^*_{E,D''}))a||=
||(\lambda I - (\bar{\partial}_E + \delta_E + \alpha) a||
\ge ||(\lambda I - (\bar{\partial}_E + \delta_E))a|| - ||\alpha a||
= ||\Sigma_j \ (\lambda - \mu_j)a_j \phi_j|| - ||\alpha a||
\ge  \sqrt{ \Sigma_j [(Im(\lambda))^2+ (Re(\lambda - \mu_j))^2]||a_j||^2} - ||\alpha|| \ ||a||
\ge |Im(\lambda)| - ||\alpha||$ because the eigenvalues $\mu_j$ are real.
Hence, if $|Im(\lambda)| >||\alpha||$, then necessarily
$(\lambda I - (\bar{\partial}_E + \bar{\partial}^*_{E,D''}))a \neq 0$. This also implies
the claim since $
\Delta_{E, \bar{\partial}} = (\bar{\partial} + \bar{\partial}^*_{E,D''})^2 $.
The parabola described in lemma \ref{lemma.AN1}
is the image under $z \mapsto z^2$ of the lines $Im \ z = \pm ||\alpha||$.

\vspace{.1in}
The analytic continuation of the zeta function
is addressed as follows. Let $\Pi_{E,K,g}$ denote the
projection onto the span of the generalized eigenspaces
with eigenvalues with real parts less than $K >0$.
Set $Q(E,K,g) = (1-\Pi(E,K,g))$ and $\Delta_{E,K,\bar{\partial},q}
= Q(E,K,g) \Delta_{E,\bar{\partial}}$ acting on $(p,q)$ forms, intuitively this
is the  heat kernel
projected onto the subspace of generalized eigenspaces
with real parts at least $K$.
By the estimates of Seeley \cite{Seeley1}, the associated heat kernel
 is most concretely and explicitly given by the
contour integral:
$$
e^{- \Delta_{E,K,\partial, q }\, t} \stackrel{ \ definition \ }{:=}  \ \frac{1}{2\pi i} \int_{\gamma_K}
 e^{-\lambda t} \ (\lambda - \Delta_{E,K,\partial, q}  )^{-1}
\ d\lambda
$$
 where $\gamma_K$ is is the parabola containing all
the spectrum given by $X = A \ Y^2 - B$ with $B > ||\alpha||^2$ and $0 < A < (1/(4||\alpha||^2)$
traced counterclockwise.

One  assumes that $K>0$ is not the real part of the any
generalized eigenvalue of $\Delta_{E,\partial, q }$, for any $q$,  then
$e^{- \Delta_{E,K,\partial, q }\, t}$  is well defined, trace class
for $t>0$,
and moreover,
$$
Tr(e^{- \Delta_{E,K,\partial, q }, t }) =  \Sigma_{j }
\ e^{- \ \lambda_{q,K,j}\, t } \quad  \mbox{for} \ t>0,
$$
where the sum is over the generalized eigenvalues of $\Delta_{E,\partial}$
acting on $(p,q)$ forms
counted with multiplicities with real part greater than $K$.

\vspace{.1in}

Let $\{ \nu(q,K,1), \cdots , \nu(q,K,N(q,K)) \}$ be a complete
enumeration counted with multiplicities
of the finitely many generalized eigenvalues of $\Delta_{E,\partial}$
having real parts less than $K$. Tautologically one has:
$$
\Sigma_{j=1}^{N(q,K)} \ e^{-\nu(q,K,j)\  t } + Tr(e^{- \Delta_{E,K,\partial, q} \, t })
    =  Tr(e^{- \Delta_{E,\partial, q} \, t })
$$

By Seeley's results the short time asymptotics of the trace of the complete heat kernel
is of the form
$$
Tr(e^{-\Delta_{E,\partial, q} \, t}) = \Sigma_{j=0}^M \ t^{-n/2+j} a_{q,j} + 0(t^{-n/2 +M+1})
$$
where the $a_{q,j} = \int_W \ b_{q,,j,u } $ are constructed out of integrals of locally constructed
forms $b_{q,j,u}$ depending only on finitely many
covariant derivatives of the Hermitian metric $g$ on $TW$ and the connection
form of $D$.
Recall $n = dim \ W = 2N$.

Expressed another way,
with
$$
Tr(e^{- \Delta_{E,K,\partial, q } \, t }) + (\Sigma_{j=1}^{N(q,K)} \ e^{-\nu(q,K,j)\  t })
      - (\Sigma_{j=0}^N \ t^{n/2+j} a_{q,j}) = \mu_M(t)
$$
where $\mu_M(t) = 0 (t^{-n/2+M+1})$.

Another implication of Seeley's analysis is control on the rate of growth
of the real parts of the eigenvalues, which give uniform estimates
$$
Tr(e^{- \Delta^{E,\bar{\partial},K,q} \, t }) \le C e^{-D t} \ for \ t>1 \ ,
$$
for some constants $C>0,D>0$. Also one gets convergence
of the zeta function:
$$
\zeta_{q,K}(s) := \Sigma_j \frac{1}{(\lambda_{q,K,j})^{-s}},  \ for \ Re(s) > n/2 \ .
$$

\vspace{.1in}

These facts may  be used  to give an explicit formula
for the analytic continuation of this zeta function to $s=0$
and its derivative at zero,
as in Cheeger's paper \cite{Cheeger1}.

First one has  for $Re(s) > n/2$ the convergent representation:
$$
\zeta_{q,K}(s) = \frac{1}{\Gamma(s)} \int_0^{\infty} \ t^{s-1} Tr(e^{-\Delta_{E,\bar{\partial},K,q} \, t }) \ dt
$$

Taking $M > n/2+1$ , choosing $\epsilon >0$, replacing $\frac{1}{\Gamma(s)}$
by $\frac{s}{\Gamma(s+1)}$ and substituting, one has  equivalently:
$$
\begin{array}{l}
\zeta_{q,K}(s) = \frac{s}{\Gamma(s+1)} \int_0^{\epsilon} \ t^{s-1} \ \mu_M(t) \ dt \\
     - \frac{s}{\Gamma(s+1)} \int_0^{\epsilon} \ t^{s-1} \ (\Sigma_{j=1}^{N(q,K)} \ e^{-\nu(q,K,j)\  t }) \ dt \\
\hspace{.3in} + \frac{s}{\Gamma(s+1)} \int_0^{\epsilon} \ t^{s-1} \ (\Sigma_{j=0}^M \ t^{-n/2+j} a_{q,j}) \ dt \\
\hspace{.6in}   + \frac{s}{\Gamma(s+1)} \int_{\epsilon}^{\infty} \ t^{s-1} Tr(e^{-\Delta_{E,\bar{\partial},K,q}\, t}) \ dt \\
\end{array}
$$

By $\Gamma(1) =1, \Gamma'(1) = \gamma$, and $M > n/2+1$, the first term has analytic
continuation to $Re(s)>-1$ by this same convergent expression. Its derivative by $s$
is,  $(\frac{1}{\Gamma(s+1)} - \frac{s}{\Gamma(s+1)^2}(\Gamma'(s+1)) )
 \int_0^{\epsilon} \ t^{s-1} \ \mu_N(t) \ dt
+  \frac{s}{\Gamma(s+1)} \int_0^{\epsilon} \ (log(t) t^{s-1}) \ \mu_K(t) \ dt$. Evaluated at
$s=0$ this yields $ \int_0^{\epsilon} \ t^{-1} \ \mu_K(t) \ dt$ which converges.

The second term once $N(q,K)$ is added is similar,
$-\frac{s}{\Gamma(s+1)} \int_0^{\epsilon} \ t^{s-1} \ \newline (\Sigma_{j=1}^{N(q,K)} \
( e^{-\nu(q,K,j) \  t } - 1)) \ dt $ since $e^{-\nu(q,K,j) \  t } - 1 = 0(t)$,
  converges for $Re(s) >-1$. The derivative evaluated at $s=0$ yields
$ -\int_0^{\epsilon} \ t^{-1} \ (\Sigma_{j=1}^{N(q,K)} \
( e^{-\nu(q,K,j) \  t } - 1))  \ dt$ which converges.

The canceling term is  $-N(q,K) \
\frac{s}{\Gamma(s+1)} \int_0^{\epsilon} \ t^{s-1} \cdot  1 \ dt = -N(q,K) \
\frac{s}{\Gamma(s+1)} \frac{\epsilon^s}{s}
= -N(q,K) \ \frac{\epsilon^s}{\Gamma(s+1)}$. This expression  is analytic
for  $Re(s)>-1$ and has derivative at $s=0$ equal to $-N(q,K)(log(\epsilon) - \gamma)$.

In the third term, omitting the constant term $t^0$, i.e. omitting $j = n/2$,  can be integrated,
$\frac{s}{\Gamma(s+1)} \int_0^{\epsilon} \ t^{s-1}
 \ (\Sigma_{j=0, j \neq n/2}^M \ t^{-n/2+j} a_{q,j}) \ dt   = \frac{s}{\Gamma(s+1)} \times
\newline (\Sigma_{j=0, j \neq n/2}^M \
\frac{\epsilon^{-n/2+j+s}}{-n/2+j+s} a_{q,j}) $. It analytically extends to $s=0$ and has derivative
$ (\frac{1}{\Gamma(s+1)} - \frac{s}{\Gamma(s+1)^2}(\Gamma'(s+1)) )  \times(\Sigma_{j=0, j \neq n/2}^M \
\frac{\epsilon^{-n/2+j+s}}{-n/2+j+s} a_{q,j})
+\frac{s}{\Gamma(s+1)} (\Sigma_{j=1}^N \ \frac{\epsilon^{-n/2+j+s}}{-n/2+j+s} a_j) \times \newline log(-n/2+j+s)$.
Evaluated at $s=0$ this is  $(\Sigma_{j=0, j \neq n/2}^M \
\frac{\epsilon^{-n/2+j}}{-n/2+j} a_{q,j})$.

The term $a_{n/2,q}\ t^0$ contributes to the zeta function at zero
the term  \newline $a_{n/2,q}\ (log(\epsilon) - \gamma)$.

The last term is uniformly convergent, extends by this formula to $s=0$ and has at that
point the derivative $\int_{\epsilon}^{\infty} \ t^{-1}\ Tr(e^{-\Delta_{q,K}\, t}) \ dt$,
similar to the first term.

\vspace{.1in}

In toto, the analytic continuation of $\zeta_{q,K}(s)$ to $s=0$ has the derivative at
$s=0$ computed as above. That is,

\begin{thm} \label{thmdiff}
For the notation above,
$$
\begin{array}{l}
\zeta_{q,K}'(0) = \int_0^{\epsilon} \ t^{-1} \ \mu_M(t) \ dt  \\
-\int_0^{\epsilon} \ t^{-1} \ (\Sigma_{j=1}^{N(q,K)} \
( e^{-\nu(q,K,j) \  t } - 1))  \ dt \\
 + (a_{n/2,q}-N(q,K))(log(\epsilon) - \gamma) +(\Sigma_{j=0,j\neq n/2}^M \
\frac{\epsilon^{-n/2+j}}{-n/2+j} a_{q,j}) \\
+ \int_{\epsilon}^{\infty} \ t^{-1}\ Tr(e^{-\Delta_{E,\bar{\partial},K,q}}) \ dt.
\end{array}
$$

Let $0<K<L$ be chosen so that $\Delta_{E,\bar{\partial},q}$ acting on $A^{p,q}(W,E)$ for any $q$ has
no generalized eigenvalue with real part equal to
$K$ or $L$. Let $\lambda_k,k=1,\cdots, T$ be a complete
enumeration, counting with multiplicities,
of the generalized eigenvalues of $\Delta^{\flat}_q$
with real  parts between $K$ and $L$. Denote by
$log(\lambda_k)$ the principle value  of the
eigenvalue $\lambda_k$.

Then
\begin{eqnarray}
\zeta_{q,L}'(0) - \zeta_{q,K}'(0) = - \Sigma_{k=1}^T \ log(\lambda_k) \label{eqn.diffzeta}
\eeqn

\end{thm}

Note that this theorem implies lemma \ref{lemma.3BB}.

\vspace{.1in}

Proof of the equation \ref{eqn.diffzeta}:
  The formula for $\zeta_{q,L}'(0)$ and $ \zeta_{q,K}'(0)$ just  derived
may be subtracted. The terms involving $\mu_M(t)$ and $a_{q,j}$ cancel, yielding
$$
\begin{array}{l}
\zeta_{q,L}'(0) - \zeta_{q,K}'(0) =- \int_0^{\epsilon} \ t^{-1} \ (\Sigma_{k=1}^{T} \
( e^{-\lambda_k \  t } - 1))  \ dt \\
-T(log(\epsilon) - \gamma)
+ \int_{\epsilon}^{\infty} \ t^{-1}\ (\Sigma_{k=1}^T \ e^{-\lambda_k t})  \ dt \\
= \Sigma_{k=1}^T \ [- \int_0^{\epsilon} \ t^{-1} \
( e^{-\lambda_k \  t } - 1))  \ dt -(log(\epsilon) - \gamma) \\
+
\int_{\epsilon}^{\infty} \ t^{-1}\  e^{-\lambda_k t}  \ dt]
= - \Sigma_k \ \log(\lambda_k) \\
\end{array}
$$

To see the last equality, recall that
 for $Re(\lambda)>0$, one has  taking principle values
$\frac{1}{\lambda^s}$ analytic with derivative
$\log(\lambda)(\frac{1}{\lambda^s})$ which evaluates at $s=0$ to $-log(\lambda)$.

On the other hand there is also the convergent integral expression
$\frac{1}{\lambda^s} = \frac{s}{\Gamma(s+1)} \int_0^{\infty} \ t^{s-1} Tr(e^{-\lambda t}) \ dt$
which by the above analysis has analytic continuation to $s=0$ with derivative
given by $[- \int_0^{\epsilon} \ t^{-1} \
( e^{-\lambda_k \  t } - 1))  \ dt -(log(\epsilon) - \gamma) +
\int_{\epsilon}^{\infty} \ t^{-1}\  e^{-\lambda_k t}  \ dt]$. Hence, these are equal.

\vspace{.1in}

Let $g(u), a \le u \le b$ be varying Hermitian metrics on the complex manifold $W$. Suppose that
$K>0$ is not an generalized eigenvalue for any
of the  $\Delta_{E,\bar{\partial},q,g(u)}$ acting on $A^{p,q}(W,E)$ for any $q$ and
any $u$. By the formulas
for the spectral projections $\Pi_{E,K,q,g(u)}$, the  modified
heat kernels, $e^{-\Delta_{E,\bar{\partial},K,q,g(u)} \, t}$,
vary smoothly with  $u$ and one  may address the
question of computing the variation of the analytic continuation
$\zeta_{q,K,g(u)}'(0)$ with the metric. Here one explicitly
records the metric dependence.

The foundational initial formula, as in Ray and Singer \cite{RaySinger1},
is:

\begin{lemma} \label{lemma.AN2}
$$
\begin{array}{l}
d/du\ \Sigma_{q=0}^N \ (-1)^q q \ Tr[ e^{-\Delta_{ E,\bar{\partial},K,q,g(u)} t} ] \\
= \  t \ d/dt[
\Sigma_{q=0}^N \ (-1)^q \ Tr( e^{-\Delta_{E,\bar{\partial},K,q,g(u)} \, t}
 \quad \alpha_{u,p,q} )]
\end{array}
$$
where  $\alpha_{u,p,q} = (\star_{u,p,q})^{-1} \ d/du \ ( \star_{u,p,q})=
(\hat{\star}_{t,p,q})^{-1} \ d/du \ ( \hat{\star}_{t,p,q}) $ acting on
$A^{p,q}(W,E)$.

\end{lemma}

Proof of lemma \ref{lemma.AN2}:
Note that spectral projection $\Pi_{E,K,q,g(u)}=(\Pi_{E,K,q,g(u)})^2 $
has finite rank, and $Q_{E,K,q,g(u)} = 1 -\Pi_{E,K,q,g(u)}$
is also a projection so \newline $(Q_{E,K,q,g(u)})^2 = Q_{E,K,q,g(u)}$.
In particular, $d/du \ Q_{E,K,q,g(u)}  = -d/du \ \Pi_{E,K,q,g(u)}
= - d/du \ (\Pi_{E,K,q,g(u)})^2 =  -d/du \ \Pi_{E,K,q,g(u)}  \ \Pi_{E,K,q,g(u)}
-  \Pi_{E,K,q,g(u)} d/du \ \Pi_{E,K,q,g(u)}$ is of finite rank and consequently
is of trace class.

Differentiating $(Q_{E,K,q,g(u)})^2 = Q_{E,K,q,g(u)}$
 yields  $Q_{E,K,q,g(u)}\  d/du \ (Q_{E,K,q,g(u)}) + d/du \ (Q_{E,K,q,g(u)}) \ Q_{E,K,q,g(u)}
= d/du \ (Q_{E,K,q,g(u)})$. Applying $Q_{E,K,q,g(u)}$ to each side gives
$2 \ Q_{E,K,q,g(u)}\  d/du \ (Q_{E,K,q,g(u)}) \ Q_{E,K,q,g(u)}  \newline
=  \ Q_{E,K,q,g(u)} \  d/du \ (Q_{E,K,q,g(u)}) \ Q_{E,K,q,g(u)}$. Hence,
$$
 Q_{E,K,q,g(u) } \  d/du \ (Q_{E,K,q,g(u)}) \ Q_{E,K,q,g(u)} = 0
$$

In particular,
$d/du \ ( Tr(e^{-\Delta_{E,\bar{\partial},K,q,g(u), t} } )) =
d/du ( Tr( Q_{E,K,q,g(u)} \ e^{- \Delta_{E,\bar{\partial}}  \ t }))
= Tr (  [d/du \  Q_{E,K,q,g(u)}] \ \ e^{- \Delta_{E,\bar{\partial} } \ t})
+ Tr( Q_{E,K,q,g(u)} \  d/du ( e^{- \Delta_{E,\bar{\partial}}  \ t})) \newline
=Tr (  [Q_{E,K,q,g(u)}\  d/du \ (Q_{E,K,q,g(u)}) + d/du \ (Q_{E,K,q,g(u)}) \ Q_{E,K,q,g(u)}
] \ \ e^{- \Delta_{E,\bar{\partial} } \ t})
+ Tr( Q_{E,K,q,g(u)} \  d/du ( e^{- \Delta_{E,\bar{\partial}}  \ t})$.
On the other hand, by the above since $Q_{E,K,q,g(u)}$ commutes with
$e^{ - \Delta_{\bar{\partial},q,g(u)} \ t} $,
$Tr (  [Q_{E,K,q,g(u)}\  d/du \ (Q_{E,K,q,g(u)}) \ e^{- \Delta_{E,\bar{\partial},g(u) } \ t})
= Tr (  [(Q_{E,K,q,g(u)})^2 \  d/du \ (Q_{E,K,q,g(u)}) \ e^{- \Delta_{E,\bar{\partial},g(u) } \ t})
= Tr (  [Q_{E,K,q,g(u)}\ \newline  d/du \ (Q_{E,K,q,g(u)}) \ e^{- \Delta_{E,\bar{\partial},g(u) } \ t} \ Q_{E,K,q,g(u)})
\newline = Tr (  [Q_{E,K,q,g(u)}\  d/du \ (Q_{E,K,q,g(u)}) \  Q_{E,K,q,g(u)} \ e^{- \Delta_{E,\bar{\partial},g(u) } \ t} \ )
=0$ and similarly
$Tr (\   d/du \ (Q_{E,K,q,g(u)}) \ Q_{E,K,q,g(u)}\ \ e^{- \Delta_{E,\bar{\partial},g(u) } \ t})
=0$. These in combination with the above prove the following formula:
$$
d/du \ ( Tr(e^{-\Delta_{E,\bar{\partial},K,q,g(u)} t } ))
= Tr( Q_{E,K,q,g(u)} \  d/du ( e^{- \Delta_{E,\bar{\partial},g(u)}  \ t}))
$$

Now the methods of Ray and Singer in either of \cite{RaySinger1, RaySinger2}
apply verbatim to complete the proof of lemma \ref{lemma.AN2}.

\begin{thm} \label{thm.ANvariation}
Under the assumptions above, the analytic continuation
of $\zeta_{q,K}(s)$ to $s=0$ has the property:
$$
d/du \ [ \Sigma_q \ (-1)^q q \ \zeta_{q,K}'(0)]
= \Sigma_{q=0}^n \ (-1)^q \  Tr( \Pi_{q,K} \  \cdot \alpha ) + a_{n/2,p,q}
$$
where the prime \ $'$ \ denotes the derivative
and $\alpha = \star^{-1} \ d/du(\star)$.

\end{thm}

Following Ray and Singer, a consequence of lemma \ref{lemma.AN2} is
a formula for the  variation of our desired quantity, for $Re(s) >n/2$:
$$
\begin{array}{l}
 d/du [ \Sigma_q \ (-1)^q q \ \zeta_{q,K}(s)] \\
= d/du \ [\Sigma_q (-1)^q q \  \frac{1}{\Gamma(s)}\int_0^{\infty} \ t^{s-1}
Tr(e^{-\Delta_{ E,\bar{\partial},K,q,g(u)} \ t} \ dt] \\
= \frac{1}{\Gamma(s)}\int_0^{\infty} \ t^s \
d/dt[\Sigma_q \ (-1)^q \ Tr( e^{-\Delta_{ E,\bar{\partial},K,q,g(u)} \  t} \cdot \alpha )] \ dt\\
= -\frac{s^2}{\Gamma(s+1)} \Sigma_q \ (-1)^q \ \int_0^{\infty} \ t^{s-1} \
Tr( e^{-\Delta_{ E,\bar{\partial},K,q,g(u)} \  t} \cdot \alpha ) \ dt \\
=-\frac{s^2}{\Gamma(s+1)} \Sigma_q \ (-1)^q \ \int_0^{\infty} \ t^{s-1} \
Tr( e^{- \Delta_{ E,\bar{\partial},q,g(u)}\ , t } \cdot \alpha \\
- \Pi_{E,K,q} \  e^{- \Delta_{ E,\bar{\partial},q,g(u)} \ , t }
\cdot \alpha ) \ dt \\
=-\frac{s^2}{\Gamma(s+1)} \Sigma_q \ (-1)^q \ \int_0^{\infty} \ t^{s-1} \
Tr( e^{- \Delta_{ E,\bar{\partial},q,g(u)} \, t } \cdot \alpha  ) \ dt \\
 +\frac{s^2}{\Gamma(s+1)} \Sigma_q \ (-1)^q \ \int_0^{\infty} \ t^{s-1} \
Tr( - \Pi_{E,K,q} \  (e^{-\Delta_{ E,\bar{\partial},q,g(u)}  \, t } -Id)
\cdot \alpha ) \ dt\\
-\frac{s^2}{\Gamma(s+1)} \Sigma_q \ (-1)^q \ \int_0^{\infty} \ t^{s-1} \
Tr( \Pi_{E,K,q} \
\cdot \alpha ) \ dt
\end{array}
$$

Now $n= dim \ W = 2N$ is even, so
 the asymptotic expansion for \newline $ Tr(e^{- \Delta_{E,\bar{\partial},q}}\ \alpha_{E,p,q,u}) $ for small
time has a term $a_{n/2,p,q,u} \ t^0$   in its expansion about  $t=0$ and the rest positive and
finitely many negative powers of $t$.
 Then inserting this   shows that $\int_{\epsilon}^{\infty} \ t^{s-1} \
Tr(e^{- \Delta_{E,\bar{\partial},q}\  t }\ \alpha_{E,p,q,u})  ) \ dt \newline +
\int_0^{\epsilon} \ t^{s-1} \
[Tr(e^{- \Delta_{E,\bar{\partial},q}\ t }\ \alpha_{E,p,q,u}) - a_{n/2,p,q,u}\ t^0] ) \ dt$ extends as a holomorphic
function to $s=0$, so the derivative of the analytic continuation of the
first term at $s=0$ vanishes. Similar remarks hold for the term involving \newline
$Tr( - \Pi_{q,K,g(u)} \  (e^{- \Delta_{E,\bar{\partial},q}  \, t } -Id)
\cdot \alpha_{E,p,q,u} )$.

 The integral in the last
term has a pole of order one with residue $Tr( \Pi_{q,K,g(u)} \ \cdot \alpha_{E,p,q,u} )$,
so the analytic continuation of this term to zero exists and has derivative
$Tr( \Pi_{q,K,g(u)} \ \cdot \alpha_{E,p,q,u} )$.
Similarly, the integral for $\int_0^{\epsilon} \ t^{s-1} \ a_{n/2,p,q} \ t^0 \ dt$
has a pole of order 1 with residue $a_{n/2,p,q}$, so the analytic continuation
of this term to zero exists and has derivative $a_{n/2,p,q}$. In toto,
this proves the theorem \ref{thm.ANvariation}.

This theorem \ref{thm.ANvariation} is recorded as lemma
\ref{lemma.3A2}.

\vspace{.1in}
The completely parallel analysis for the flat Laplacian $\Delta_E^\flat$, yields
the corresponding results in this setting.
In combination with the analogue of the variation theorem \ref{thm.MAIN},
which is proved by a completely parallel analysis, the independence of
metric result, \ref{thm.MAIN2}, follows.

\begin{flushleft}

Sylvain E. Cappell

Courant Institute, N.Y.U.

251 Mercer Street

New York, NY 10012

email: cappell@courant.nyu.edu

\vspace{.1in}

Edward Y. Miller

Mathematics Department

Polytechnic Institute of New York

Six MetroTech Center

Brooklyn, NY 11201

email: emiller@poly.edu

\end{flushleft}

\end{document}